\documentclass[12pt]{book}
\title{Interval Functions and their Integrals}
\author{Ralph Henstock}
\date{Ph.D.~Thesis, December 1948}

\newcommand{\vt}{\vspace{5pt}\\}
\newcommand{\vs}{\vspace{10pt}}

\newcommand{\R}{\mathbf{R}}

\newcommand{\ve}{\varepsilon}

\newcommand{\var}{{\mathrm{{Var}}}}

\usepackage[pdftex]{graphicx}

\begin{document}
\maketitle

\noindent
{\sf{\large{\textbf{Ralph Henstock's Ph.D.~thesis \\
Preliminary note }}}}

\vs
\noindent
{\sf Ralph Henstock (1923--2007) worked in non-absolute integration, including the Riemann-complete or gauge integral which, independently, Jaroslav Kurzweil also discovered in the 1950's.

As a Cambridge undergraduate Henstock took a course of lectures, by J.C.~Bur\-k\-ill, on the integration of interval functions. Later, under the supervision of Paul Dienes in Birkbeck College, London, he undertook research into the ideas of Burkill (interval function integrands) and Dienes (Stieltjes integrands); and he   presented this thesis   in December 1948.

The thesis contains the germ of Henstock's later work, in terms of overall approach and methods of proof. For example, a notable innovation is a set of axioms for constructing any particular system of integration. 
This  highlights the features held in common by various systems, so that a particular property or theorem can, by a single, common proof, be shown to hold for various kinds of integration.

Within this approach, Henstock's thesis places particular emphasis on various alternative ways of selecting Riemann sums, as the primary distinguishing feature of different systems of integration. This idea was central to his subsequent work and achievement.

Of interest also are those ideas in the thesis which were effectively abandoned in his subsequent work.

In addition to Henstock's own insights at that stage of his work, this thesis provides a good overview of the literature and state of knowledge of integration of the non-Lebesgue kind at that time.  

These are good enough reasons for transcribing the thesis. Another pressing reason is that  the ink and paper of the near 70 years old copy of the thesis in the  Archive in the University of Ulster Library in Coleraine---Henstock's personal, annotated copy---are showing signs of  deterioration. 

\vs

\noindent
Pat Muldowney

\noindent
February 2017
}
\newpage

\noindent
{\large{\textbf{Abstract.}}}

\vs

\noindent
The majority of papers dealing with general interval functions are concerned mainly with their differentiation, and integration is given little space. In this thesis we therefore examine the known results about the integration of interval functions and embed the results in a general theory which uses a family $F$ of ``Riemann successions'' of divisions.

We set up suitable axioms\footnote{The families $F$ and $G$ and the axioms are entirely new.} for $F$ in one-dimensional Cartesian space, and find that all axioms except one (and two special axioms) apply with little change to $n$-dimensional space, and even to an abstract space. The exceptional axiom can with some difficulty be generalised to $n$-dimensional Cartesian space, and becomes quite complicated for the abstract space.

Burkill integrals and norm-limits are defined using $F$, and both reduce to the original extended Burkill integral (Burkill [6]) when the integration uses every possible division of an $n$-dimensional interval. When the Burkill integrals or norm-limits are finite, new additive integrals or limits are defined. The relation between these and the $\sigma$-limit is investigated and the properties of all integrals and limits are examined in some detail. Inequalities are given similar to the false result of Saks [10] (213, Theorem 3) and the true result (214, \textsection 4, Lemma), and a few existence theorems are set out in Chapter 1. Chapter 2 deals with the case $n=2$.

In Chapter 3 we use the theory of Chapters 1 and 2 to examine some special functions. In particular, functions of bounded variation, and functions arising in Kempisty's ``integration around a set'' and my ``density integration'', and Appendix 2 gives a brief note of further generalisations, particularly to the case of an abstract space.

A copy of the author's paper (Henstock [15] \textit{On interval functions and their integrals}) is included in subsidiary matter. 

\vs

\noindent
{\sf{[The latter paper is on-line at:}}

\noindent
\texttt{{\small{{https://academic.oup.com/jlms/article-abstract/s1-21/3/204/845051/}}}}

\noindent
\texttt{{\small{{On-Interval-Functions-and-Their-Integrals?redirectedFrom=fulltext}}}}

\noindent
--{\sf{P.M.]}}

\newpage
 
\noindent
{\large{\textbf{Preface.}}}

\vs
\noindent
On reading papers devoted to the theory of functions of intervals one cannot fail to observe that the majority are concerned mainly with their differentiation, and little space is given to the consideration of their integration. In this thesis we therefore collect together the known results about the integration of interval functions, and embed the results in a general theory which uses a family $F$ of \textit{Riemann successions} of divisions.

Our intervals are in $n$-dimensional Cartesian space, and we set up suitable axioms for the theory. All these axioms, except axioms (vi), (ix), (x), are independent of the dimension, and in fact with little change can be used for an abstract space of objects, as is shown in Appendix 2. This simplifies the theory a good deal, and enables us to take the case $n=1$ alone for results which do not depend on axioms,(vi), (ix),(x). (Axioms (ix) and (x) are rather special.) Two types of integral are defined, Burkill integrals and norm-limits; both kinds are included in the original (extended) integral defined by Burkill in the special case when the integration uses every possible division of an interval. The first integral uses the aggregate of all limit-points of certain sequences of sums given by divisions; and the second uses the limits as $e \rightarrow 0$ of the upper and lower bounds of sums given by divisions with meshes of thickness less than $e>0$.

The theory which uses primes ($n-1$ flats) of divisions in $n$-dimensional space will need an axiom like axiom (vi) for  $n=1$. We give the corresponding axiom (vi)' for $n=2$, dealing with broken lines, and by analogy this can be extended to general $n$. In Appendix 2, after a long discussion, a corresponding axiom for an abstract space is given. The theory which follows from the use of (vi) or (vi)' has a analogue in $n$ dimensions ($n>2$), and a partial analogue in abstract space.

When the Burkill integrals or norm-limits are finite, new integrals are defined, and all except one are additive. The relation between these new integrals and the $\sigma$-limit of Getchell [11] and Hildebrandt [13] is investigated.

Inequalities are given for the integrals and limits, which are useful and simplify subsequent theory, such as is given in Chapter 3 (Special Functions). In that chapter we discuss, among other functions, the functions of bounded variation, and functions arising in Kempisty's ``integration around a set'' and my ``density integration''. 

The present thesis will only be concerned with the ``Riemannian'' integration of interval functions in $n$-dimensional Cartesian space, and will deal with no results on differentiation nor on the general extension of an interval function to produce even a set function for open sets. Differentiation and extension are outside the scope of the thesis. But in Chapter 3 \textsection\textsection 5, 6, we consider the ``Riemannian'' integration of two kinds of interval functions whose integrals sometimes provide a convenient extension of the original interval function.

An integral of an interval function appears unexpectedly in the example given in Appendix 1. But since the theory has no similarity with the main part of the thesis, the example has been put in Appendix 1 instead of in Chapter 3. Appendix 2 gives a brief outline of the way how to extend the results of the thesis to abstract space.

A copy of the author's paper, \textit{On interval functions and their integrals} (Jour.~London Math.~Soc.~21 (1946) 204--209) is also included.

\vs

\noindent
{\sf{[This paper is on-line at:}}

\noindent
\texttt{{\small{{https://academic.oup.com/jlms/article-abstract/s1-21/3/204/845051/}}}}

\noindent
\texttt{{\small{{On-Interval-Functions-and-Their-Integrals?redirectedFrom=fulltext}}}}

\noindent
--{\sf{P.M.]}}

\newpage

\noindent
{\large{\textbf{Contents.}}}









\begin{enumerate}
\item[{Chapter 1}]
\hfill  p.~\pageref{Chapter 1}
\begin{itemize}
\item[1.]
Introduction. \hfill  p.~\pageref{Chapter 1 Introduction}

Examples (i) to (xviii), additive functions, bounded functions, continuous functions.
\item[2.]
Burkill integration.
\hfill  p.~\pageref{Chapter 1 Burkill integration}

The families $F$ and $G$, together with axioms (i) to (v).
\item[3.]
Burkill integration and points of division.
\hfill  p.~\pageref{Burkill integration and points of division}
 
 Axioms (vi) to (viii), and singularities with respect to additivity. Inequalities for sums over non-overlapping interval contained in $\R$.
\item[4.]
$k$-integration and the $\sigma$-limit.
\hfill  p.~\pageref{k-integration and the sigma-limit}

The families $F_k$, $F_{k'}$. The $k$-limits and $k'$-limits. The special axioms (ix) and (x). Inequalities between the several integrals and limits. Inequalities for sums over non-overlapping intervals. The $\sigma$-limit. Upper and lower $\sigma$-limits and their relation to certain upper and lower $k'$-limits.

\item[5.]
Existence theorems.
\hfill  p.~\pageref{Existence theorems}

\end{itemize}

\item[{Chapter 2}]
\hfill  p.~\pageref{Chapter 2}
\begin{itemize}
\item[1.]
Introduction. \hfill  p.~\pageref{Chapter 2 Introduction}

Examples (i) to (ix).
\item[2.]
Burkill integration.
\hfill  p.~\pageref{Chapter 2 Burkill integration}

The ordinary (extended) and the restricted integrals, and the families $H$ and $L$. A result of ``Fubini'' type.

\item[3.]
Burkill integration and lines of division.
\hfill  p.~\pageref{Burkill integration and lines of division}

The common frontier of two non-overlapping finite sums of rectangles $T$. Axiom (vi)' and singularities with respect to additivity. Inequalities for integrals and for sums over non-overlapping rectangles.

\item[4.]
$k$-integration and the $\sigma$-limit.
\hfill  p.~\pageref{Chapter 2 k-integration and the sigma-limit}

$k$-successions, $H_k$, axiom (ix)', and $k$-limits. Inequalities by analogy with Chapter 1.
\end{itemize}

\item[{Chapter 3}]
\hfill  p.~\pageref{Chapter 3}
\begin{itemize}
\item[1.]
Functions of bounded variation. \hfill  p.~\pageref{Functions of bounded variation}

The variations of $g(I)$ and its norm-limits. Properties of the norm-limits. The variation singularities. Decompositions of the norm-limits and $k'$-limits. The upper and lower, positive and negative, variations.

\item[2.]
Absolutely continuous functions.
\hfill  p.~\pageref{Absolutely continuous functions}

The $g(I)$ are also of bounded variation, and the upper and lower norm-limits are absolutely continuous. Absolutely semi-continuous functions. 

\item[3.]
Functions which are monotone on subdivision.
\hfill  p.~\pageref{Functions which are monotone on subdivision}

The upper and lower norm-limits, with special cases when $g(I)$ is continuous or absolutely continuous. The $k'$-limit exists.

\item[4.]
$g(T)=g_1(I_x). g_2(I_y)$ when $T=[I_x,I_y]$ in 2-dimensions.
\hfill  p.~\pageref{product}

\item[5.]
Integration around a set $E$ of points.
\hfill  p.~\pageref{Integration around a set E of points}
 
Kempisty's definitions. Various inequalities. The necessary and sufficient condition for a norm-limit to exist around a set $E$. Absolutely continuous functions, and functions of bounded variation, around a set $E$. 

\item[6.]

Density integration. \hfill p.~\pageref{Density integration}

If $g(I)$ is additive and absolutely continuous the density integral for the measurable set $E$ is the Lebesgue integral $g(E)$ of $g(I)$. In order that the upper and lower density integrals are both finite for every measurable set $E$ it is necessary and sufficient that the norm-limit of $g(I)$ should exist and be absolutely continuous.

\end{itemize}

\item[{Appendix 1:}]
Function Space and Density Integration.
\hfill  p.~\pageref{Appendix 1}

Let $P(f) = \sum_{i=1}^n b_i g(I_i)$ whenever $f=f(x)=
\sum_{i=1}^n b_i c(I_i;x)$, where $b_1, \ldots , b_n$ are real constants, and $I_1, \ldots , I_n$ are disjoint intervals in $(0,1)$. Let
\[
||f||^2 = \int_0^1 |f(x)|^2 dx,\;\;\;\;\;\;||f|| \geq 0.
\]
If $P(f)$ is continuous with respect to $||f||$ then $P(c(E;x))$ is the density integral of $g(I)$ for $E$, and is $g(E)$ (the Lebesgue extension).

\item[{Appendix 2:}]
Further Generalisations\hfill  p.~\pageref{Appendix 2}

To $n$-dimensions. To abstract spaces and many-valued $g(I)$.

\item[{Bibliography.}]
Bibliography and references.
\hfill  p.~\pageref{Bibliography}
\end{enumerate}

\newpage

\noindent
{\large{\textbf{General Definitions.}}}

\vs
\noindent
\textbf{Unless otherwise stated, we assume that the intervals of each set of intervals are non-overlapping.}

\vs
\noindent
\textbf{One dimension.}

\vs
\noindent
We take two fixed numbers $A<B$ and an interval $W$ or $A \leq x \leq B$, and we call $W$ the \textit{fundamental linear interval}. To each pair of numbers $a,b$ such that $A \leq a<b \leq B$ there corresponds four intervals $I$, namely, an \textit{open} interval $(a,b)$ (i.e.~$a<x<b$), a \textit{closed} interval $[a,b]$ (i.e.~$a\leq x\leq b$), and two \textit{half-closed} intervals, $[a,b)$ (i.e.~$a\leq x < b$) and $(a,b]$ (i.e.~$a < x\leq b$). Dienes [14] denotes collectively the four intervals corresponding to $a,b$ by $a\mbox{---}b$. Each has \textit{length} $mI$ equal to $b-a$. Thus $W$ is a closed interval with length $B-A$.

Two intervals \textit{overlap} when they have an interval in common. Two intervals \textit{abut} when they have a common end-point but do not overlap. A finite sum of intervals is denoted by $I^\sigma$. The \textit{norm} of a (non-overlapping) set of intervals is the upper bound $mI$ for intervals $I$ of the set.

The following definition is substantially due to Dienes [14].

Given $R$, an $I^\sigma$ in $W$, then a finite number of points $x_0<x_1< \cdots <x_n$ in $R'$, which include all the end-points of the separate intervals of $R$, will divide $R$ into a finite number of intervals $I_1, \ldots , I_m$. Taking a definite arrangement of brackets (or \textit{bracket convention}) for each interval, the distinct or abutting intervals $I_1, \ldots, I_m$ will be referred to as a \textit{division} $D$ of $R$. Thus there will be $4^m$ divisions corresponding to the \textit{points of division} $x_0<x_1< \cdots <x_n$. The arrangement of brackets at an $x_i$ will be called the \textit{bracket convention} at $x_i$.

The reason\footnote{See definition of $E'$ below. P.M.} for using $R'$ instead of $R$ is to ensure that two abutting intervals of $R$ are treated as one interval in forming a division. If this is not desired, a remark to this effect can be introduced in the definition of the families $F$ and $G$ which are given later. Since $F$ and $G$ are quite general, this will make no difference to the subsequent theory.

The \textit{norm} of a division is the norm of the corresponding set of intervals. Summation for a division $D$ over a set $E$ of points (or of intervals) is denoted by $(E;D)\sum$. If a division $D$ is not in question, or is assumed known, we put $(E)\sum$ for the summation. If also $E$ is assumed known we put $\sum$. This notation is more explicit, and seems to be more useful, than that used by Burkill [6] and others, in which the sum over a mesh-system $\mathcal E$ is denoted by $g(\mathcal E)$.

If $E$ is a set of points, $CE$ is the \textit{complement} of $E$ with respect to $W$, i.e.~ $CE=W-E$. The \textit{interior} of $E$, or the set of interior points of $E$, is denoted by $E^o$, and a point $x$ is \textit{inside}  of $E$ if $x$ is in $E^o$. The \textit{derived set}, or the set of limit-points of $E$, is denoted by $E'$. The \textit{boundary} or \textit{frontier} $FE$ of $E$ is \[E'.(CE)'.\] If the set $E_1$ is contained in the set $E$ we put $E_1 \subset E$. The \textit{characteristic function} of the set $E$ is denoted by $c(E;x)$, so that $c(E;x)=1$ when $x$ is in $E$, $c(E;x)=0$ when $x$ is not in $E$. $(x)$ is the set whose sole member is the point $x$. The \textit{distance} between two sets $E_1, E_2$ denoted by $\rho(E_1;E_2)$, so that
\[
\rho(E_1;E_2) = \mbox{g.l.b.}|x_1-x_2|
\]
for a point $x_1$ in $E_1$ and a point $x_2$ in $E_2$.

\vs
\vs
\noindent
\textbf{Two dimensions.}

\vs
\noindent
The coordinates of a point $P$ are denoted by $x,y$, so that $P=(x,y)$. We take four fixed numbers $A<B$, $C<D$, and a rectangle $W$, or $A\leq x\leq B$, $C\leq y \leq D$, which we write as $[A,B;C,D]$, and we call $W$ the \textit{fundamental rectangle}.

To each set of four numbers $a,b,c,d$ such that $A\leq a<b\leq B$, $C\leq c<d\leq D$, there corresponds $16$ rectangles $T$ according to the $16$ ways of including or not including the \textit{sides} 
\[
a\leq x\leq b,\;\;y=c;\;\;\;a\leq x\leq b,\;\;y=d;\;\;\;x=a,\;\;c\leq y\leq d;\;\;\;x=b,\;\;c\leq y \leq d.
\]
If all sides are included the rectangle is \textit{closed}, and denoted by $[a,b;c,d]$. Thus $W$ is closed. If a side with $x=a$ or $b$, or $y=c$ or $d$, is omitted, the rectangle is denoted by $[a,b;c,d]$ with an index $o$ on the corresponding letter. For example, if $a<x \leq b$, $c \leq y \leq d$, together with the points $(a,c)$, $(a,d)$, then the rectangle is $[a^o,b;c,d]$. The rectangle $[a^o,b^o;c^o,d^o]$ is \textit{open}. The $16$ types of rectangles form our two-dimensional intervals.

The \textit{diameter} $\delta(T)$ of such rectangles is
\[
\sqrt{(b-a)^2+(d-c)^2}.
\]
The \textit{norm} of a set of (non-overlapping) rectangles is the upper bound of $\delta(T)$ for rectangles $T$ of the set. A finite sum of the rectangles is denoted by $T^\sigma$.

Let $V$ be a $T^\sigma$ (in $W$). Divide $V'$ up into a finite number of non-overlapping rectangles $T_1, \ldots ,T_n$, taking a definite arrangement of included and non-included sides around each rectangle $T_i$ ($1\leq i \leq n$). Then $T_1, \ldots , T_n$ will be referred to as a \textit{division} $D$ of $V$. If $V$ is a $T$, and if the lines of division extend right across $T$, we call the division a \textit{restricted division}.

The definitions of division and restricted division follow those introduced by Burkill [6], \textsection 2, but extended to include bracket conventions.

The \textit{norm} of a division is the norm of the corresponding set of rectangles. Taking $E$ to be a set of points, of polygonal lines, or of rectangles $T$, we have definitions of $(E;D)\sum$, $(E)\sum$, $\sum$, as for one-dimensional intervals. If $E$ is a set of points the definitions of $CE$, $E^o$, $E'$, $FE$, $E_1 \subset E$, are as in one dimension. The \textit{characteristic function} of the set $E$ is denoted by $c(E;x,y)$, so that $c(E;x,y)=1$ when $(x,y)$ is in $E$, and $c(E;x,y)=0$ when $(x,y)$ is not in $E$. The set $E$ is \textit{connected} if every two points of $E$ can be connected by a polygonal line lying entirely in $E$.

\chapter{One-dimensional Integration.}\label{Chapter 1}
\section{Introduction.}\label{Chapter 1 Introduction}
Burkill [6] has given the following definition.
   
If rules exist which associate a unique number $g(I)$ with each interval $I$ of a set $S$ of intervals, then $g(I)$ is a \textit{function of intervals } (or an \textit{interval function}) defined in the set $S$.

Let us denote by $\bar E$ the set of all $I^\sigma$ formed by finite sums (not necessarily non-overlapping) of points and the interiors of intervals of $S$. Then if $R_1,R_2$ are in $\bar S$, so is $R_1+R_2$, and $R_1$ is an $I^\sigma$.

In the present chapter we suppose that the intervals of $S$ are contained in the fundamental linear interval $W$. Many examples of functions of intervals occur in the different branches of Analysis. In the following selection the interval function sometimes depends on a further parameter $\xi$ which lies in a range $\rho(I)$ determined by $I$. Putting this interval function as $g(I;\xi)$ we define
\[
g_2(I) = \mbox{l.u.b.} g(I;\xi),\;\;\;\;\;g_2(I) = \mbox{g.l.b.} g(I;\xi),
\]
for $\xi$ in $\rho(I)$. Whenever $g_1(I)$ and $g_2(I)$ can be defined they are (strictly speaking) the interval functions which arise from the given example. 
\begin{enumerate}
\item
Let $a<b$, and $I$ be any of the four intervals $a\mbox{---}b$. Then the \textit{length} of $I$ is $mI=b-a$.
\item
If also $f(x)$ is a bounded real function (of a real point), then the \textit{Stieltjes difference} of $f(x)$ is $S(I)=S(f;I) =f(b)-f(a)$. (Dienes [14].)
\item If at every point $x$ (in $W$) the limits $f(x-0)$, $f(x+0)$ exist, 
then we can define $3^8$ \textit{differences} by taking open intervals $9a,b)$ in 2., and replacing $f(a)$ by $f(a-0),f(a), f(a+0)$, and similarly replacing $f(b)$ by $f(b-0),f(b)$, or $f(b+0)$ (thus obtaining $3^2$ differences for open intervals) and repeating these results for each of the other 3 intervals $a\mbox{---}b$. Only $3^2=9$ of these differences are additive in $W$, in a sense to be defined (Dienes [14]).
\item
The interval function $|S(f;I)|$ is studied for the total variation of functions $f(x)$. (Burkill [6] \textsection 2, number 1.)
\item
We study 
\[
\frac{S(I)}{mI}=\frac{f(b)-f(a)}{b-a}
\]
in the Differential Calculus. (Burkill [6], Introduction.)
\item
Let $t$ be in $W$ and $x=x(t), y=y(t)$ be two continuous functions so that $(x,y)$ traces out a rectifiable curve. Then its length is found by using
\[g(I) = \sqrt{S(x;I)^2+S(y;I)^2}.\]
(Saks [10] Chapter 2, 1.1.)
\item
Hellinger uses 
\[
g(I) = \frac{S(f;I)^2}{S(h;I)},
\]
where $h(x)$ is a bounded strictly increasing function, in his thesis on quadratic forms. See Hobson [3] 257.
\item
Hobson [3] 257 has generalised 7., to become
\[
g(I) = \frac{S(f;I)^{i+1}}{S(h;I)^i}\;\;\;\;\;\;\;\;(i>0).
\]
See also Titchmarsh [19] 384--6 (\textsection 12.44), in which the case $h(x) \equiv x$ is dealt with..
\item
Another type of interval function is the least upper bound $M(I) = M(f;I)$ of $f(x)$ in $I$, where we take bracket conventions into account.
\item
If $f(x)$ is monotone increasing and $h(x)$ is bounded, then for Riemann-Stieltjes integration we consider
\[
g(I) = M(h;I)S(f;I).
\]
(Burkill [6] \textsection 2, number 4.) The upper and lower limits of sums $(R;D) \sum g(I)$, where $R$ is an $I^\sigma$, for every sequence of divisions $D$ of $R$ with norms tending to zero, give the upper and lower Riemann-Stieltjes integrals
\[
\overline{\int}_R h\,df \;\;\;\;\mbox{and }\;\;\;\; 
\underline{\int}_{\;R} h\,df
\]
respectively
When we replace the special $g(I)$ by an arbitrary $g(I)$ we obtain an integration which was first studied by Burkill.
\item
H.L.~Smith [8] has defined a \textit{Stieltjes mean integral} in which he uses
\[
g(I) = \frac 12 \left(h(a) + h(b)\right) S(f;I)
\]
where $h(x)$ is bounded.
\item
Let $h(x)$ be bounded and $f(x-0)$ and $f(x+0)$ exist (finite) for all $x$ in $W$. Then W.H.~Young [2] considers
$g(I;\xi)\; =$
\[
=\;h(a)\left(f(a+0) -f(a)\right)
+h(\xi)\left(f(b-0) -f(a+0)\right)
+h(b)\left(f(b) -f(b-0)\right)
\]
for $I = a$---$b$, $a<\xi<b$.
\item
Hellinger and Radon have considered
\[
g(I) = \frac{M(f;I) S(k;I) S(l;I)}{S(h;I)}
\]
where $h(x)$ is bounded and strictly increasing. (See Lebesgue [16] 296.)
\item
When $f(x)$ is measurable with bounds $A$ and $B$, Lebesgue considers
\[
g(I) = a . \mbox{meas}E\left[a \leq f(x)<b\right]
\]
for $I=[a,b)$. (Burkill [6] \textsection 2, number 2.)
\item
Let $t$ be in $W$ and $x(t),y(t)$ be two continuous functions, so that $(x,y)$ traces out a closed (plane) curve. To find its area we consider 
\[
g(I) = \frac 12 \left(x(a)y(b) -x(b)y(a) \right)
\]
for $I=a$---$b$. (Burkill [7] 321.)
\item
Let $x_1, \ldots , x_n, \ldots$ be a sequence of distinct points and $f(x)$ be such that $\sum_{n=1}^\infty f(x_n)$ is absolutely convergent. Put
\[
g(I) = (I) \sum f(x_n).
\]
Then $g(I)$ is additive, in a sense to be defined.
\item
Let $h(I)$ be additive (in a sense to be defined). Then Dienes [14] considers $g(I;\xi) = f(\xi)h(I)$, taking various $\rho(I)$ including $I^o$.

\item
Let a vector $\mathbf{r}=\mathbf r(x,y)$ be defined at every point of the $(x,y)$ plane and let $\mathbf r$ lie in that plane. Let $A \leq t \leq B$, and in that range let $x(t)$ and $y(t)$ be two real continuous functions of bounded variation, so that $\mathbf{z}(t) \equiv \left(x(t),y(t)\right)$ traces out a simple closed Jordan curve $C$ lying in the $(x,y)$ plane. For $a<b$ let $I=a$---$b$ and $J$ be the line joining $\mathbf{z}(a)$ to $\mathbf{z}(b)$, and take a parameter $\xi$ on $J$. Then for $(x,y)$ on $J$, $\mathbf r$ is a function of $\xi$. Taking $\mathbf N$ as the vector pointing away from $C$, which is perpendicular to $J$ and of the same length, put 
\[
g(I;\xi) = \mathbf{r}(\xi).\mathbf{N}.
\]
We can then consider integration $\int_C \mathbf{r}.\mathbf{dz}$ and Green's Theorem.

\end{enumerate}
An interval function $g(I)$ is \textit{additive in} $S$ if
\begin{enumerate}
\item[(a)]
for every $x,y,z$ with $x<y<z$, and with $(x,y], (y,z), (x,z)$ in $S$ we have
\[
g\left((x,y]\right) + g\left((y,z)\right) = g\left((x,z)\right);
\]
\item[(b)]
for every $x,y,z$ with $x<y<z$, and with $(x,y), [y,z), (x,z)$ in $S$ we have
\[
g\left((x,y)\right) + g\left([y,z)\right) = g\left((x,z)\right);
\]
\end{enumerate}
and if the six corresponding relations using different bracket conventions at $x$ and $z$ are also true.

The interval functions of 1., 2., 16., are additive, but in general the rest are not.

In 3., given any one of the 9 differences for the open interval and supposing $S$ to be all intervals in $W^o$, then we can define $g(I)$ over half-closed intervals by using equations (a) and (b), and then define $g(I)$ over closed intervals by
\begin{enumerate}
\item[(c)]
\[
g\left((x,y)\right) + g\left([y,z]\right) = g\left((x,z]\right);
\]
\end{enumerate}
Hence only 9 of the $3^8$ differences are additive in $W^o$. The non-open intervals of $A$---$x$, $x$---$B$ (with $A<x<B$) and $A$---$B$ are exceptional in that the value of $g(I)$ for them cannot be defined by additivity.

Note that if $S$ contains all the intervals in $W^o$, and no more, then (a), (b), (c) are sufficient for the definition of additivity, since the remaining 5 relations follow easily. 

The following definitions, and results (1.1) to (1.3), are due to Burkill [6].

An interval function $g(I)$ is \textit{bounded} in $S$ if there are numbers $H,K$ such that for all $I$ in $S$ we have $H<g(I)<K$.

Let $\omega (\delta,x)$ be the upper bound of $|g(I)|$ for every $I$ of $S$ contained in an interval with centre $x$ and length $\delta$. Then for fixed $x$, $\omega(\delta,x)$ is monotone decreasing as $\delta \rightarrow 0$, if $\omega(\delta,x)$ exists. Hence $\omega(\delta,x)$ tends to a limit, $\omega(x)$ say, which is called the \textit{oscillation} of $g(I)$ at $x$ (with respect to the set $S$). Naturally this is defined if in every neighbourhood of $x$ there is an interval of $S$. And if so, $\omega(\delta,x)$ is defined for every $\delta$. Then $g(I)$ is \textit{continuous at} $x$ if $\omega(x)=0$, and $g(I)$ is \textit{continuous in a set} $E$ if it is continuous at each point of $E$.

\vs
\noindent
\textbf{(1.1)} \textit{If $g_1(I)$ and $g_2(I)$ are continuous at $x$, so is $g_1(I)+g_2(I)$, provided that $\omega(x)$ for the latter is defined.}

\noindent Let 
$\omega_1(\delta,x)$, $\omega_2(\delta,x)$, $\omega_3(\delta,x)$
be the oscillation functions of $g_1, g_2, g=g_1+g_2$. Then, for each $\delta>0$,
\[
0 \leq \omega(\delta,x) \leq \omega_1(\delta,x)+\omega_2(\delta,x)\,\;\;\;\rightarrow 0
\]
as $\delta \rightarrow 0$. Hence $\omega(x)=0$.

\vs
\noindent
\textbf{(1.2)} \textit{If $g_1(I)$ is continuous at $x$,  $g_2(I)$ is bounded, and $S$ is the same for both, then $g_1g_2$ is continuous at $x$.}

\noindent For some $M$, $0 \leq |g_1g_2| \leq M|g_1|$ so that
\[
0\leq \omega(\delta,x) \leq M\omega_1(\delta,x),\;\;\;\;\;\;\;\;
0\leq \omega(x) \leq M\omega_1(x).
\]

\noindent
\textbf{(1.3)} \textit{If at each point $x$ of a closed interval $J$, $\omega (x) \leq k$, then given $\ve>0$ we can find $\delta = \delta(\ve)$ so that $|g(I)| <k+\ve$ for every $I\subset J$ with $mI<\delta$.}

\noindent
Suppose false. Then there is a sequence of intervals $I_1,I_2, \ldots$ in $J$ such that $|g(I_n)| \geq k+\ve$ and $mI_n \rightarrow 0$ as $n \rightarrow \infty$. The centres of intervals $I_n$ have at least one limit-point, $\xi$ say, and since $J$ is closed, $\xi$ is in $J$. 
But $\omega(\xi) \geq k+\ve$, contradicting $\omega(x) \leq k$ in $J$.

\vs
\noindent
\textbf{Corollary:}\textit{ The uniformity of continuity theorem for functions of intervals. }
 (Take $\omega(x) =0, k=0$.)

The interval functions of all the examples, except posibly 5., 7., 8., 13., are bounded. The interval functions of 1., 6., 15., 18., are continuous, but in general the rest are not. If, however, $f(x)$ is continuous, the interval functions of 2., 3., 4., 10., 11., 12., are also continuous.

\section{Burkill integration} \label{Chapter 1 Burkill integration}
A \textit{Riemann succession} of divisions of $R$, an $I^\sigma$, is a sequence of divisions of $R$ in which norm$(D_n) \rightarrow 0$ as $n \rightarrow \infty$. (Dienes [14].)

In integration we often consider a family $F$ of Riemann successions $\{D_n\}$ of divisions formed from intervals of $S$, supposing that $F$ has the following properties.
\begin{enumerate}
\item[(i)]
If $R$ is in $\bar S$ there is at least one Riemann succession of divisions of $R$ which is in $F$.
\item[(ii)]
If $\{D_n\}$ is a Riemann succession in $F$ and $n_1,n_2, \ldots$ is any sequence of integers tending to infinity then $\{D_{n_i}\}$ is in $F$.
\item[(iii)]
If $R_1, R_2$ are non-overlapping and in $\bar S$, and if $\{D_{1,n}\}, \{D_{2,n}\}$ are in $F$, where $\{D_{i,n}\}$ is a Riemann succession of divisions of $R_i$ ($i=1,2$), then $\{D_n\}$ is in $F$, where $D_n$ is the division of $R_1+R_2$ formed from the intervals of the divisions $D_{1,n}$ and $D_{2,n}$, i.e.~$D_n = D_{1,n}+D_{2,n}$.
\end{enumerate}
We often suppose that $F$ also has the following properties.
\begin{enumerate}
\item[(iv)]
Let $\{D_n^{(i)}\}$ ($i=1,2, \ldots$) be a set of Riemann successions in $F$, of divisions of $R$ in $\bar S$. Then if
\[
D_1^{(1)},\;\;\;D_2^{(1)},\;\;\;D_1^{(2)},\;\;\;D_3^{(1)},\;\;\;D_2^{(2)},\;\;\;D_1^{(3)},\;\;\;D_4^{(1)},\;\;\;\ldots
\]
is a Riemann succession it is in $F$.

This may be called an axiom of ``closure''. It is needed in order that limits using certain Riemann successions may be the same as certain norm-limits, as will be seen later.
\item[(v)]
Let $R$ (in $\bar S$) be the sum of non-overlapping sets $R_1, \dots , R_m$, each in $\bar S$, with the property that \textbf{every} division $D$ of $R$ can be divided up to form $m$ divisions $D^{(1)}, \ldots , D^{(m)}$ such that $D^{(i)}$ is a division of $R_i$ ($i=1, \ldots , m$). Then if the Riemann succession $\{D_n\}$ of divisions of $R$ is in $F$, the Riemann successions 
\[
\left\{D_n^{(1)}\right\}, \ldots , \left\{D_n^{(m)}\right\}
\]
are also in $F$, where $D_n^{(i)}$ is over $R_i$ ($i=1, \ldots , m$) and $D_n =D_n^{(1)} + \cdots + D_n^{(m)}$.
\end{enumerate}
Since our intervals are in one dimension, this axiom can in general only apply when $R'$ is the sum of $q\geq m$ closed disjoint intervals. In two or more dimensions this would not necessarily be so.

Denote by $G$ the set of divisions which occur in the Riemann successions of $F$. If $S$ consists of all intervals in $W$, denote $S$ by $S_1$. If $S =S_1$ and $F$ comprises all Riemann successions over all $I^\sigma$ in $W$, denote $F$ by $F_1$. The corresponding $G$ is $G_1$.

\vs
\noindent
\textbf{(2.01)} \textit{If $R_1, R_2$ are two non-overlapping $I^\sigma$ in $\bar S$, and if $D^{(i)}$ in $G$ is a division of $R_i$ ($i=1,2$), then $D^{(1)} + D^{(2)}$ is in $G$.}

\noindent
For, by hypothesis and (ii), we can take $D^{(1)}, D^{(2)}$ as the first divisions of two Riemann successions, and then (iii) gives the result.

\vs
\noindent
\textbf{(2.02)} \textit{ If $R$ is in $\bar S$, and if $D$ in $G$ is a division of $R$, and if $E$ is an $I^\sigma$ formed from a subset of the set of intervals comprising $D$, then there is an $I^\sigma$ in $\bar S$, $R_1$ say, such that $R_1 = R-E'$.}

\noindent
For let $R_2$ be the $I^\sigma$ formed from those intervals of $D$ which are not in $E$. Then $R_2'$ has the required properties.

\noindent
\textbf{(2.03)} \textit{    From (iv) and (ii), if If $R_1, R_2$ are two non-overlapping $I^\sigma$ in $\bar S$, and if $D^{(i)}$ in $G$ is a division of $R_i$ ($i=1,2$), then $D^{(1)} + D^{(2)}$ is in $G$.}

\vs
\noindent When $F$ exists and $S$ is not vacuous we may deduce the following.

\vs
\noindent
\textbf{(2.04)}  \textit{If 
$I$ is in $\bar S$ then $I$ contains intervals $J$ of $S$ with arbitrarily small norm. (This follows from (i).) Also $\omega(x)$ exists for all $x$ inside $I$.}

\vs
Following Burkill [6], but extending his definition to the family $F$, we define the \textit{upper Burkill integral}
$(F)\overline{\int}_R \,g(I)$ of $g(I)$ over $R$ in $\bar S$, with respect to the family $F$, to be the upper bound of
\[
\overline{\lim}_{n \rightarrow \infty} (R;D_n)\sum g(I)
\]
for every Riemann succession $\{D_n\}$ in $F$, of divisions of $R$. The \textit{lower Burkill integral }
$(F)\underline{\int}_{\,R}\, g(I)$  is the lower bound of 
\[
\underline{\lim}_{n \rightarrow \infty} (R;D_n)\sum g(I).
\]
If for $R$ the upper and lower Burkill integrals are finite and equal, we say that $g(I)$ \textit{is integrable over $R$ in $\bar S$, with respect to the family $F$, in the sense of Burkill}, writing the common value as $(F)\int_R g(I)$.

The above definition uses RIemann successions of divisions. Another definition, given below, produces norm-limits in the sense of Moore and Smith [4]. These integrals could just as well have been called Burkill integrands, but for distinction we will call them norm-limits. (See Getchell [11] and Hildebrandt [13].)

The \textit{upper norm-limit $(N;G)\overline{\int}_R \,g(I)$ of $g(I)$ over $R$ in $\bar S$, with respect to the set $G$ of divisions}, is the greatest lower bound of numbers $d(e)$ for all $e>0$, where $d(e)$ is the least upper bound of $(R;D)\sum g(I)$ for all $D$ in $G$ such that $D$ is a division of $R$ and norm$(D) < e$.

As $e \rightarrow 0$, $d(e)$ is monotone decreasing, and so the lower bound is also the limit. A similar definition holds for the \textit{lower norm-limit} $(N;G)\underline{\int}_{\,R}\,g(I)$. If for $R$ the upper and lower norm-limits are finite and equal, we say that \textit{the norm-limit of $g(I)$ exists over $R$ in $\bar S$, with respect to the family $G$ of divisions, in the sense of Moore and Smith}, and we write the common value as $(N;G)\int_R g(I)$.

By (i) the upper and lower Burkill integrals and the upper and lower norm-limits all exist, though some may be infinite.

\vs
\noindent
\textbf{(2.05)}  \textit{ $(F)\overline{\int}_R g(I) \geq (F)\underline{\int}_{\,R}\, g(I)$. If $(F)\overline{\int}_R g(I)= (F)\underline{\int}_R g(I)$ then for every $\{D_n\}$ in $F$, of divisions of $R$,
\[
\lim_{n \rightarrow \infty} (R;D_n)\sum g(I) = (F)\overline{\int}_R g(I).
\]}
These results are obvious from the definitions.

\vs
\noindent
\textbf{(2.06)} \textit{ $(N;G)\overline{\int}_R g(I) \geq (N;G)\underline{\int}_{\,R}\, g(I)$. If $(N;G)\overline{\int}_R g(I)$ exists then given $\ve>0$ there is a $\delta>0$ such that if $D$ in $G$ is a division of $R$ with norm$(D)<\delta$, then 
\[
\left|(R;D)\sum g(I) - (N;G)\int_R g(I) \right|<\ve.
\]
If $(N;G) \overline{\int}_R g(I) = -\infty$ then given $A>0$ there is a $\delta>0$ such that if $D$ in $G$ is a division of $R$ with norm$(D)<\delta$, then
\[
(R;D) \sum g(I) < -A.
\]
If $(N;G) \underline{\int}_{\,R}\, g(I) = +\infty$ then given $A>0$ there is a $\delta>0$ such that if $D$ in $G$ is a division of $R$ with norm$(D)<\delta$, then
\[
(R;D) \sum g(I) >A.
\]}
These results are obvious from the definitions.

\vs
\noindent
\textbf{(2.07)}  \textit{
Let $\{D_n\}$ in $F$ be of divisions of $R$ in $\bar S$. If $L$ is a limit-point of the sums $(R;D_n) \sum g(I)$ as $n \rightarrow \infty$, then there is in $F$ a $\{D'_n\}$ of divisions of $R$ such that
\[
\lim_{n \rightarrow \infty}(R;D_n)\sum g(I) =L.
\]}
(By (ii) and elementary analysis.)

\vs
\noindent
\textbf{Corollary:}
\textit{ We can take $L$ equal to 
\[
\overline{\lim}_{n \rightarrow \infty} (R;D_n) \sum g(I)
\;\;\;\;\;\;\mbox{ or }\;\;\;\;\;\;
\underline{\lim}_{n \rightarrow \infty} (R;D_n) \sum g(I)
\]
}
\noindent
\textbf{(2.08) (a)}
\textit{ 
If $R$ is in $\bar S$ and $(F)\overline{\int}_R g(I) > -\infty$ then given $A<(F)\overline{\int}_R g(I)$ there is in $F$ a Riemann succession $\{D_n\}$ of divisions of $R$ with
\[
\lim_{n\rightarrow \infty} (R;D_n) \sum g(I) >A.
\]}
\textbf{(2.08) (b)}  
\textit{ 
If $F$ satisfies (iv) then it contains a Riemann succession $\{D_n\}$ of divisions of $R$ with
\[
\lim_{n \rightarrow \infty} (R;D_n)\sum g(I) =
(F)\overline{\int}_R g(I).
\]}
For (a), by definition, there is in $F$ a Riemann succession  $\{D'_n\}$ with
\[
\overline{\lim}_{n\rightarrow \infty} (R;D'_n) \sum g(I)>A.
\]
Then, (2.07, Corollary) gives the result. For (b) we can suppose $(F)\overline{\int}_R \,g(I)>-\infty$, since otherwise (2.05) will give the result. Let $\{D_n^{(i)}\}$ be the $\{D_n\}$ of (a) for $A=A_i$, where 
\[
A_1<A_2< \cdots \;\;\;\mbox{ and }\;\;\;
A_i \rightarrow (F)\overline{\int}_R g(I)
\]
as $i \rightarrow \infty$, and such that norm$\left(D_n^{(i)}\right)<(i+n)^{-1}$. This can be arranged by (ii). Then the sequence of (iv) is a Riemann succession and so by (iv) is in $F$. Denote this succession by $\{D'_n\}$. Then for all $i>0$, 
\[
\overline{\lim}_{n\rightarrow \infty} (R;D'_n) \sum g(I)\geq A_i, \;\;\;\;\mbox{ so }\;\;\;\;
\overline{\lim}_{n\rightarrow \infty} (R;D'_n) \sum g(I)
\geq
(F)\overline{\int}_R \,g(I),
\]
and by definition we have equality. Using (2.07, Corollary) we have the result.

Similar results hold for lower integrals.

\vs

\noindent
\textbf{Corollary:}
 In a similar way we can prove that 
\textit{if $F$ satisfies (iv), and $R$ is fixed, then the set of limit-points of all sequences of sums $(R;D_n)\sum g(I)$, for $\{D_n\}$ in $F$, is closed.}

(2.08) is stated in such a way as to include the case 
$(F)\overline{\int}_R\, g(I)=+\infty$.

\vs
\noindent
\textbf{(2.09) (a)}
\textit{
Let $(N;G)\overline{\int}_R\, g(I) < +\infty$. Then given  
$A>(N;G)\overline{\int}_R\, g(I)$
there is an $e>0$ such that every division of $R$ and in $G$, with norm$(D)<e$ has $(R;D) \sum g(I)<A$.}

\vs
\noindent
\textbf{(2.09) (b)}
\textit{
Let $(N;G)\overline{\int}_R\, g(I) >-\infty$
Then given $e>0$ and 
$A<(N;G)\overline{\int}_R\, g(I)$
there is a  division $D$ of $R$ and in $G$, with norm$(D)<e$ and \[(R;D) \sum g(I)<A.\]}

\noindent
For (a) we take $e>0$ so that $d(e) <A$. Then the result follows. For (b) we have that $d(e) >A$ for all $e>0$. The result follows by definition of $d(e)$.

Similarly for lower norm-limits.

Next we have a connection between the two integrations.

\vs
\noindent
\textbf{(2.10) (a)}
\[
(N;G)\underline{\int}_{\,R}\, g(I)
 \leq
 (F)\underline{\int}_{\,R}\, g(I) \leq
 (F)\overline{\int}_{\,R}\, g(I)
 \leq  
 (N;G)\overline{\int}_{\,R} \, g(I).
 \]

\noindent
\textbf{(2.10) (b)}
\textit{If $F$ satisfies (iv) then}
\[
(N;G)\underline{\int}_{\,R}\, g(I)
 =
 (F)\underline{\int}_{\,R}\, g(I),
 \;\;\;\;\;\;\;\;
 (N;G)\overline{\int}_{\,R}\, g(I)
 =  
 (F)\overline{\int}_{\,R}\, g(I).
 \]
For (a) we can obviously suppose that $(F)\overline{\int}_{\,R}\, g(I)>-\infty$. Then, by (2.08a), given $A<(F)\overline{\int}_{\,R}\, g(I)$, there is a Riemann succession $\{D_n\}$ in $F$, of divisions of $R$, such that
\[
\lim_{n\rightarrow \infty} (R;D_n) \sum g(I) >A.
\]
Hence there is in $G$ a division $D$ of $R$ with arbitrarily small norm, such that $(R;D) \sum g(I) >A$, so that $d(e) >A$ for all $e>0$. Hence \[(N;G)\overline{\int}_{R}\, g(I)\geq A. \]
Since $A$ is arbitrary, we have the third inequality of (a).
Similarly for lower integrals. 

For (b) suppose that for some $A$ we have
\[
(N;G)\overline{\int}_{R} g(I) >A> (F)\overline{\int}_{R} g(I).
\]
Then for every $e>0$, $d(e) >A$, so that for each integer $n$ there is a division $D_n$ of $R$, and in $G$, such that
\[
(R;D_n) \sum g(I) >A \;\;\;\mbox{ and } \;
\mbox{norm}(D_n) < \frac 1 n.
\]
By (2.03), which uses (iv) and (ii), $\{D_n\}$ is in $F$. Hence $(F)\overline{\int}_{R} \,g(I)\geq A$. This gives a contradiction. Hence, from (a), 
\[
(N;G)\overline{\int}_{R} g(I) = (F)\overline{\int}_{R} g(I).
\]
Similarly for lower integrals.

The use of the norm-limit is therefore equivalent to assuming that $F$ satisfies (iv) as well as (i), (ii).

\vs

\noindent
\textbf{Corollary: }
\textit{If $F$ satisfies (iv) then}
\[
(N;G)\underline{\int}_{\,R}\, g(I)
 =
 (F^*)\underline{\int}_{\,R}\, g(I),
 \;\;\;\;\;\;\;\;
 (N;G)\overline{\int}_{\,R}\, g(I)
 =  
 (F^*)\overline{\int}_{\,R}\, g(I),
 \]
\textit{where $F*$ is the ``closure'' of $F$, in an obvious sense.}

\vs
\noindent
\textbf{(2.11) }
\textit{If $(F)\overline{\int}_{\,R}\, g(I)=-\infty$
and either}

\noindent
\textbf{(b) }
\textit{$F$ satisfies (iv) and $(F)\overline{\int}_{\,R}\, g(I)=+\infty$ or}

\noindent
\textbf{(c) }
\textit{there is in $F$ a Riemann succession $\{D_n^{(2)}\}$ of divisions of $R_2$ such that
\[
\lim_{n \rightarrow \infty} (R_2; D_n^{(2)})\sum g(I) = +\infty,
\]
then $(F)\overline{\int}_{\,R_1 +R_2}\, g(I) = +\infty$. Similarly for lower integrals.}

\noindent
(a) By (2.08a), and given $A_i <(F)\overline{\int}_{\,R}\, g(I)$, we can find a Riemann succession $\{D_n^{(i)}\}$, in $F$ and of divisions of $R_i$, such that
\[
\lim_{n\rightarrow \infty} (R_i;D_n^{(i)}) \sum g(I) >A_i \;\;\;\;\;\;(i=1,2).
\]
Then $\;\;\;\lim_{n\rightarrow \infty} (R_1+R_2;D_n^{(1)}+D_n^{(2)}) \sum g(I)$ $=$
\[
=\;\;\lim_{n\rightarrow \infty} (R_1;D_n^{(1)})\sum g(I) +
\lim_{n\rightarrow \infty} (R_2;D_n^{(2)})\sum g(I) \;\;>\;\;A_1+A_2.
\]
By (2.08b), (b) implies (c). We therefore assume the latter.

\noindent
(c) if $\{D_n^{(1)}$ is a Riemann succession in $F$, of divisions of $R_1$, then for each integer $i$ there is an integer $n_i\geq i$ such that
\[
(R_2;D_{n_i}^{(2)})\sum g(I) > i - (R_1;D_i^{(1)})\sum g(I).
\]
Using (ii) and then (iii), $\{D_i\}$ is in $F$, where 
$D_i=D_i^{(1)}+D_{n_i}^{(2)}$. And
\[
(R_1+R_2;D_i) \sum g(I)\;\; =\;\;
(R_1;D_i^{(1)})\sum g(I)+(R_2;D_{n_i}^{(2)})\sum g(I)\;\; >\;\; i .
\]
Hence $(F)\overline{\int}_{\,R_1 +R_2}\, g(I)=+\infty$.

\vs
\noindent
\textbf{Corollary:} By (2.10) Corollary, \textit{We have similar results to (a) and (b) for norm-limits.}

\vs
\noindent
\textbf{(2.12} \textit{If $R_1, R_2$ are non-overlapping and in $\bar S$, and if}

\vs
\noindent
\textbf{(2.12) (a)} \textit{$(F)\overline{\int}_{\,R_1 +R_2}\, g(I)$ is finite then either $(F)\overline{\int}_{\,R_i }\, g(I)<+\infty$ ($i=1,2$) or one of $(F)\overline{\int}_{\,R_i}\, g(I)$ ($i=1,2$) is $+\infty$ and the other is $-\infty$. The second alternative is false if $F$ satisfies (iv).} (Saks [10] 213, Theorem 1, for $F=F_1$.)

\vs
\noindent
\textbf{(2.12) (b)} \textit{If both Burkill integrals over $R_1+R_2$ are finite then both Burkill integrals over $R_1$ and over $R_2$ are also finite.} (Saks [10] 213 Theorem 4 ($1^o$) for $F=F_1$.)

\noindent
(a) follows from (2.11), and (b) follows from (a). For if we have (say) $(F)\overline{\int}_{\,R_2}\, g(I)=-\infty$, then by the analogue of (2.11c) for lower integrals we would have $(F)\overline{\int}_{\,R_1+R_2}\, g(I)=-\infty$, contrary to hypothesis.

\vs
\noindent
\textbf{Corollary:} \textit{By (2.10 Corollary) we have similar results for norm-limits to those in the case when $F$ satisfies (iv).}

\vs
\noindent
\textbf{(2.13)} \textit{Let $R$ be as in axiom (v). Then if $F$ satisfies (v),}

\noindent
\textbf{(a)}
\[
(F)\overline{\int}_{\,R}\, g(I) = \sum_{i=1}^n (F)\overline{\int}_{\,R_i}\, g(I)
\]
\textit{whenever the right-hand side is finite.}

\noindent
\textbf{(b)}
\textit{If one or more of the $(F)\overline{\int}_{\,R_i}\, g(I)=-\infty$, and for Riemann successions in $F$, of divisions of the rest of the $I^\sigma$, the corresponding upper limits are always less than $+\infty$, then $(F)\overline{\int}_{\,R_i}\, g(I)=-\infty$.}

\noindent
\textbf{(c)}
\textit{If one or more of $(F)\overline{\int}_{\,R_i}\, g(I)=+\infty$, and if (b) does not hold, then \[(F)\overline{\int}_{\,R}\, g(I)=+\infty.\]}
(c) follows from (2.11c). For (a) we note that, by (2.11a),
\[
(F)\overline{\int}_{\,R}\, g(I) \geq \sum_{i=1}^n 
(F)\overline{\int}_{\,R_i}\, g(I).
\]
We therefore prove the reverse inequality. Given $A<(F)\overline{\int}_{\,R}\, g(I)$, let $\{D_n\}$ in $F$ and over $R$ be such that 
\[
\lim_{n \rightarrow \infty}(R;D_n)\sum g(I) >A.
\]
Such a $\{D_n\}$ exits by (2.08a). By (v) we may split up $\{D_n\}$ into $\{D_n^{(i)}\}$ in $F$ and over $R_i$ for $i=1, \ldots m$. Then
\begin{eqnarray*}
(R;D_n)\sum g(I) &=& \sum_{i=1}^m (R_i;D_n^{(i)})\sum g(I) \vt
A&<&  \sum_{i=1}^m \overline{\lim}_{n \rightarrow \infty} (R_i;D_n^{(i)})\sum g(I)  \vt
& \leq & 
\sum_{i=1}^m 
(F)\overline{\int}_{\,R_i}\, g(I).
\end{eqnarray*}
Hence $ 
(F)\overline{\int}_{\,R}\, g(I)$ is finite and $ 
(F)\overline{\int}_{\,R}\, g(I) \leq \sum_{i=1}^m  
(F)\overline{\int}_{\,R_i}\, g(I)$. Hence (a).
For (b) we have (say) for any $\{D_n\}$ in $F$ and over $R$, which by (v) gives $\{D_n^{(i)}\}$ in $F$ and over $R_i$ ($i=1, \ldots , m$),
\begin{eqnarray*}
 {\lim}_{n \rightarrow \infty} (R_1;D_n^{(1)})\sum g(I)
 &=& -\infty, \vt
  \overline{\lim}_{n \rightarrow \infty} (R_i;D_n^{(i)})\sum g(I)
 &<&+\infty \;\;\;\;\;\;(i=2, \ldots, m). \;\;\;\;\;\;\;\;\mbox{Then}\vt
  \overline{\lim}_{n \rightarrow \infty} (R;D_n)\sum g(I)
&\leq & 
\sum_{i=1}^m  \overline{\lim}_{n \rightarrow \infty} (R_i;D_n^{(i)})\sum g(I) \;\;\;=\;\;\;-\infty.
\end{eqnarray*}
Hence (b).

\vs
\noindent
\textbf{Corollary.}
\textit{Similar results hold for lower integrals and for norm-limits. Note that (b) does not hold if any $ 
(N;G)\overline{\int}_{\,R_i}\, g(I) =+\infty$.}

The following are further simple results for Burkill integrals, and so for norm-limits.

\vs
\noindent
\textbf{(2.14)}
\textit{If $g(I) \geq 0$ then $ 
(F)\underline{\int}_{\,R}\, g(I)\geq 0.$ (Burkill [6] 2.5, for $F=F_1$.}

\vs
\noindent
\textbf{(2.15)}
\textit{If $H.mI \geq g(I) \geq K.mI$ then
\[
H.mR \geq (F)\overline{\int}_{\,R}\, g(I)
\geq (F)\underline{\int}_{\,R}\, g(I) \geq K.mR.
\]
}
(Burkill [6] 2.5, for $F=F_1$.)

\vs
\noindent
\textbf{(2.16)} \textit{If $c$ is a constant and}

\noindent
\textbf{(a)} \textit{$c>0$ then 
$(F)\overline{\int}_{\,R}\, cg(I)=c(F)\overline{\int}_{\,R}\, g(I)$, $(F)\underline{\int}_{\,R}\,c g(I)
=c (F)\underline{\int}_{\,R}\, g(I)$.}

\noindent
\textbf{(b)}
\textit{$c<0$ then 
$(F)\overline{\int}_{\,R}\, c g(I) = c(F)\underline{\int}_{\,R}\, g(I)$, $(F)\underline{\int}_{\,R}\, cg(I)
= c(F)\overline{\int}_{\,R}\, g(I).$
}

\noindent
(Burkill [6] 2.3, for $F=F_1$.)

\vs
\noindent
\textbf{(2.17)}
\textit{Let $g(I) =g_1(I) + g_2(I)$. Then
\begin{eqnarray*}
(F)\underline{\int}_{\,R}\, g_1(I) + (F)\underline{\int}_{\,R}\, g_2(I) &\leq &
(F)\underline{\int}_{\,R}\, g(I) \vt
&\leq & (F)\underline{\int}_{\,R}\, g_1(I)
+ (F)\overline{\int}_{\,R}\, g_2(I) \vt
&\leq &
(F)\overline{\int}_{\,R}\, g(I) \vt
&\leq & (F)\overline{\int}_{\,R}\, g_1(I) +
(F)\overline{\int}_{\,R}\, g_2(I)
\end{eqnarray*}
whenever this has meaning.
} 
(Burkill [6] 2.6, for $F=F_1$.)

\noindent
For if $\{D_n\}$ in $F$ is over $R$ then
\[
(R;D_n)\sum g(I) = (R;D_n)\sum g_1(I) + (R;D_n)\sum g_2(I).
\]
By elementary analysis this gives
\begin{eqnarray*}
&&\underline{\lim}_{n\rightarrow \infty} (R;D_n)\sum g_1(I)
+\underline{\lim}_{n\rightarrow \infty} (R;D_n)\sum g_2(I)
\vt
&\leq &
\underline{\lim}_{n\rightarrow \infty} (R;D_n)\sum g(I) \vt
&\leq &
\underline{\lim}_{n\rightarrow \infty} (R;D_n)\sum g_1(I) +
\overline{\lim}_{n\rightarrow \infty} (R;D_n)\sum g_2(I) \vt
&\leq & \overline{\lim}_{n\rightarrow \infty} (R;D_n)\sum g(I) \vt
&\leq &
\overline{\lim}_{n\rightarrow \infty} (R;D_n)\sum g_1(I) +\overline{\lim}_{n\rightarrow \infty} (R;D_n)\sum g_2(I).
\end{eqnarray*}
Hence the first inequality in the result, taking the lower bound of the right-hand side. From the second limit inequality we have
\[
(F)\underline{\int}_{\,R}\, g(I)
\leq \underline{\lim}_{n\rightarrow \infty} (R;D_n)\sum g_1(I) + (F)\overline{\int}_{\,R}\, g_2(I).
\]
Hence the second inequality in the result, taking the lower bound of the right-hand side. Similarly for the rest.

\vs
\noindent
\textbf{Corollary:} \textit{If $g_1, g_2$ are Burkill-integrable so is $g_1+g_2$, its integral being the sum of the other two integrals.}

\vs

From inequalities such as those of Schwarz, H\"{older}, Minkowski we may deduce corresponding inequalities for Burkill integrals. Thus Burkill [6] 2.7 gives the following (for $F=F_1$).

\vs

\noindent
\textbf{(2.18)} \textit{Schwarz' inequality.}

\noindent
\textbf{(a)} 
\[
\left( (F)\overline{\int}_{\,R}\, g_1(I)g_2(I)\right)^2
\leq (F)\overline{\int}_{\,R}\, g_1^2(I)\;\;\times \;\;
(F)\overline{\int}_{\,R}\, g_2^2(I).
\]
\textbf{(b)}
\[
\left( (F)\underline{\int}_{\,R}\, \left|g_1(I)g_2(I)\right|\right)^2
\leq (F)\underline{\int}_{\,R}\, g_1^2(I)\;\;\times \;\;
(F)\overline{\int}_{\,R}\, g_2^2(I).
\]
For by Schwarz' inequality
\[
\left( (R;D)\sum g_1(I)g_2(I)\right)^2
\leq (R;D)\sum g_1^2(I)\;\;\times \;\;
(R;D)\sum g_2^2(I).
\;\;\;\;\mbox{ Hence}
\]

\begin{eqnarray*}
\overline{\lim}_{n\rightarrow \infty}\left( (R;D_n)\sum g_1(I)g_2(I)\right)^2
&\leq &
\overline{\lim}_{n\rightarrow \infty} 
\left(\left((R;D_n)\sum g_1^2(I)\right)
\left((R;D_n)\sum g_2^2(I)\right)\right) \vt
&\leq &
\left(\overline{\lim}_{n\rightarrow \infty} 
(R;D_n)\sum g_1^2(I)\right)
\left(\overline{\lim}_{n\rightarrow \infty} 
(R;D_n)\sum g_2^2(I) \right)
\vt
& \leq & \left((F)\overline{\int}_{\,R}\, g_1^2(I)\right)
\left((F)\overline{\int}_{\,R}\, g_2^2(I)\right).
\end{eqnarray*}
Hence the result (a). Also
\begin{eqnarray*}
\underline{\lim}_{n\rightarrow \infty}\left( (R;D_n)\sum \left|g_1(I)g_2(I)\right|\right)^2
&\leq &
\underline{\lim}_{n\rightarrow \infty} 
\left(\left((R;D_n)\sum g_1^2(I)\right)
\left((R;D_n)\sum g_2^2(I)\right)\right) \vt
&\leq &
\left(\underline{\lim}_{n\rightarrow \infty} 
(R;D_n)\sum g_1^2(I)\right)
\left(\overline{\lim}_{n\rightarrow \infty} 
(R;D_n)\sum g_2^2(I) \right)
\vt
& \leq & \left(\underline{\lim}_{n\rightarrow \infty} 
\left((R;D_n)\sum g_1^2(I)\right)
\right)
\left((F)\overline{\int}_{\,R}\, g_2^2(I)\right)
\end{eqnarray*}
by elementary analysis. Hence
\[
\left( (F)\overline{\int}_{\,R}\, \left|g_1(I)g_2(I)\right|\right)^2
\leq 
\left(\underline{\lim}_{n\rightarrow \infty} 
\left((R;D_n)\sum g_1^2(I)\right)
\right)
\left((F)\overline{\int}_{\,R}\, g_2^2(I)\right).
\]
Hence result (b).
Similarly from H\"older's inequality
\[
\sum |a_ib_i| \leq \left(\sum |a_i|^p \right) ^{\frac 1p}
\left(\sum |b_i|^q \right) ^{\frac 1q}
\]
for $p>1$, $q>1$, $p+q=pq$, we may derive

\vs
\vs
\vs
\noindent
\textbf{(2.19)}

\noindent
\textbf{(2.19)(a)}
\[
 (F)\overline{\int}_{\,R}\, g_1(I)g_2(I)
\leq 
\left( (F)\overline{\int}_{\,R}\, g_1^p(I) \right)^{\frac 1p}
\left( (F)\overline{\int}_{\,R}\, g_2^q(I) \right)^{\frac 1q},
\]
\textbf{(2.19)(b)}
\[
 (F)\underline{\int}_{\,R}\, \left|g_1(I)g_2(I)\right|
\leq 
\left( (F)\underline{\int}_{\,R}\,\left| g_1^p(I)\right| \right)^{\frac 1p}
\left( (F)\overline{\int}_{\,R}\,\left| g_2^q(I)\right| \right)^{\frac 1q},
\]
\textit{for $p>1$, $q>1$, $p+q=pq$.} For $p=q=2$ this is (2.18). And from Minkowski's inequality
\[
\left(\sum |a_i+b_i|^p\right)^\frac 1p \leq \left(\sum |a_i|^p \right) ^{\frac 1p} +
\left(\sum |b_i|^p \right) ^{\frac 1p}
\]
we may derive

\vs
\noindent
\textbf{(2.20)}

\noindent
\textbf{(a)}
\[
 \left((F)\overline{\int}_{\,R}\, \left|g_1(I)+g_2(I)\right|^p\right)^\frac 1p
\leq 
\left( (F)\overline{\int}_{\,R}\,\left| g_1(I)\right|^p \right)^{\frac 1p}+
\left( (F)\overline{\int}_{\,R}\,\left| g_2(I)\right|^p \right)^{\frac 1p},
\]
\textbf{(b)}
\[
 \left((F)\underline{\int}_{\,R}\, \left|g_1(I)+g_2(I)\right|^p\right)^\frac 1p
\leq 
\left( (F)\underline{\int}_{\,R}\,\left| g_1(I)\right|^p \right)^{\frac 1p}+
\left( (F)\overline{\int}_{\,R}\,\left| g_2(I)\right|^p \right)^{\frac 1p},
\]
\section{Burkill integration and points of division}
\label{Burkill integration and points of division}
In this section we consider points of division and bracket conventions at these points.

Saks [10] 211 defines the following functions when $F=F_1$ and $y$ is in $W^o$.
\begin{eqnarray*}
A(y) &=& \max \left( \overline{\lim} \left(g(I) -g(I_1)-g(I_2)\right);0\right), \vt
a(y) &=& \min \left( \underline{\lim} \left(g(I) -g(I_1)-g(I_2)\right);0\right),
\end{eqnarray*}
where $x<y<z$, $I=x$---$z$, $I_1 = x$---$y$, $I_2=y$---$z$, and the $\overline{\lim}$, $\underline{\lim}$
are taken as $x$ and $z$ tend independently to $y$. He then states that if 
$a<y<b$, $J=a$---$b$, $J_1 = a$---$y$, $J_2=y$---$b$, then

\noindent
(A)
\[
(F_1) \overline{\int}_J g(I)
= (F_1) \overline{\int}_{J_1} g(I) +(F_1) \overline{\int}_{J_2} g(I)
+A(y).
\]
That (A) is false may be seen by considering the following $g(I)$.
Let $f(x)$ be continuous in $[0,1]$ and linear in the intervals $[1-2^{-n}, 1-2^{-n-1}]$ with value 0 at $1-2^{-2n}$, 1 at $1-2^{-2n-1}$ ($n=0,1, \ldots $).
When $0 \leq a<b<1$ and $I=a$---$b$ we put $g(I) = f(b)-f(a)$. When $I'$ is $[1-2^{-2n}, 1+2^{-2n}]$ we put $g(I) =1$. Otherwise $g(I)=0$ in $[0,2]$. Then $A(1)=1$ and
\[
(F_1)\overline{\int}_{(0,1)} g(I) =
\overline{\lim}_{x \rightarrow 1-0}\left(f(x)-f(0)\right) =1,\;\;\;\;\;\;(F_1)\overline{\int}_{(1,2)} g(I)=0.
\]
In a division $D$ of $(0,2)$ let the greatest division-point less than 1 be $x$, and let $y$ be the nearest division-point greater tan $x$. Then
\[
\left((0,2);D\right) \sum g(I) = f(x) -f(0) +g(x\mbox{---}y),
\]
and $g(x$---$y)=0$ unless $f(x)=0$, when sometimes $g(x$---$y)=1$. Hence \[(F_1)\overline{\int}_{(0,2)} g(I)=1.\]
Hence (A) is false for $a=0$, $y=1$, $b=2$, and this $g(I)$, which is independent of bracket conventions. However, we need to retain a result like that of (A) (i.e.~Saks [10] 213, Theorem 3), and for symmetry, when $F=F_1$ we will use
\[
\overline{\lim} \left|g(I)-g(I_1)-g(I_2) \right|
\]
instead of $A(y)$ and $a(y)$.

Returning to the more general $F$ and $S$, supposing them to exist and be not vacuous, we may deduce the following result.

\vs
\noindent
\textbf{(3.01)}
\textit{If $I$ is in $\bar S$} then $I$ contains intervals $J$ of arbitrarily small norm, which are of the following types.

\noindent
\textbf{(a)}
\textit{Each $J$ has left (right)-hand end-point that of $I$.}

\noindent
\textbf{(b)}
\textit{For $x$ in $I^o$ then $x$ is in $J^o$ or else there are non-overlapping pairs of the $J$ with $x$ as common end-point.}

If a point $x$ is a point of division of all but a finite number of the divisions in a Riemann succession, then $x$ is a \textit{permanent point of division} of the Riemann succession. (Dienes [14].)

If a point $x$ inside $R$, an $I^\sigma$ of $\bar S$, is not a point of division of any but a finite number of the division of $R$ in a Riemann succession, then we can say that $x$ is a \textit{permanent interior point } of the Riemann succession.

\vs
\noindent
\textbf{(3.02)}

\noindent
\textbf{(a)}
\textit{If $x$ is a common end-point of two non-overlapping intervals $I_1, I_2$ in $S$ then $x$ is a permanent point of division for some Riemann succession in $F$.}

\noindent
\textbf{(b)}
\textit{If $x$ is inside $R$ an $I^\sigma$ of $\bar S$, and is never an end-point of any interval of $S$, then $x$ is a permanent interior point of every Riemann succession in $F$ of divisions of $R$.}

\noindent
(b) is obvious; (a) follows from (i) and (iii).

\vs
Let $x$ be inside an interval $I$ of $\bar S$, and $\{D_n\}$ in $F$ be over $I$. Then $x$ \textit{has the property} (p) \textit{with respect to} $\{D_n\}$ if the following are satisfied.

\noindent
(a) $I=I_1+I_2 + (x)$ where $x$ is the common end-point of $I_1, I_2$, each of them being in $\bar S$.

\noindent
(b) $\{D_{1,n}\}$ over $I_1$ and $\{D_{2,n}\}$ over $I_2$
are both in $F$, where $\{D_{1,n}\}$ and $\{D_{2,n}\}$ are defined in (c).

\noindent
(c)
If $x$ is inside an interval $J_n$ of $D_n$ then each interval of $D_n$, except $J_n$, occurs in $D_{1,n}+D_{2,n}$ with the same bracket convention at its ends, and there is one remaining interval $J_{1,n}$ 
of $D_{1,n}$ and one remaining interval $J_{2,n}$ 
of $D_{2,n}$ so that $J'_{1,n}+J'_{2,n}=J'_n$. If $x$ is a point of division of $D_n$ then each interval of $D_n$ occurs in $D_{1,n} + D_{2,n}$ with the same bracket conventions at its ends.

\vs

We sometimes suppose that $F$ has the property (vi).

\vs
\noindent
\textbf{(vi)} Every point $x$ satisfying condition (a) of property (p), has the property (p) with respect to every $\{D_n\}$ over $I$ and in $F$.

If $x$ has the property (p) for every $\{D_n\}$ in $F$ and over $I$, for which $x$ is a permanent interior point, we define
\[
B(x;\{D_n\}) = \overline{\lim}_{n \rightarrow \infty} \left| g(J_n) - g(J_{1,n}) - g(J_{2,n})\right|.
\]
If $x$ has the property (p) for every $\{D_n\}$ in $F$ and over $I$, for which $x$ is a permanent interior point, we define
\[
b(x;I;F) = \mbox{l.u.b.}B(x;\{D_n\})\;\;\;\mbox{ for all such }\;\;\{D_n\}.
\]
Then $b(x;I;F) \geq 0$. If it is $>0$ we say that \textit{$x$ is a singularity in $I$ of $g(I)$ with respect to additivity and respect to $F$.}

We now consider similar definitions for the norm-limit method. In the definition of the property (p) we replace $F$ by $G,$  $ \{D_n\}$ over $I$  by a single division $D$ over $I$, $J_n$ by $J$, $J_{1,n}$ by $J_1$, and $J_{2,n}$ by $J_2$. We then have a property (P) with respect to $D$, and put
\[
C(x;G;e) = \mbox{l.u.b.}\left| g(J) - g(J_1) - g(J_2)\right|
\]
for every triad $J, J_1, J_2$ of intervals from a division $D$ (over $I$ and in $G$) with norm less than $e>0$, and with $x$ inside the interval $J$ of $D$. If no such $D$ exists we put $C(x;G;e)=0$. Then
\[
c(x;I;G) =\lim_{\ve\rightarrow 0} C(x;G;e).
\]
Obviously $c(x;I;G) \geq 0$ If it is $>0$ we say that $x$ \textit{is a singularity in $I$ of $g(I)$ with respect to additivity and with respect to $G$.}
We \textbf{cannot} put
\[
c(x;I;G) = \overline{\lim}\left| g(J) -g(J_1)-g(J_2)\right| \;\;\mbox{ as }\; \max\{mJ; mJ_1; mJ_2\} \rightarrow 0
\]
when we have a general $F$, since there may be divisions $D$ in $G$ with arbitrarily small $J$, but with norms greater than some positive number. However, since in Henstock [15] we have $F=F_1$, the results for the $D(x,y,z)$ and $\sigma(y)$ of that paper are not invalidated.

\vs
\noindent
\textbf{(3.03)}

\noindent
\textbf{(a)}
\textit{If $x$ has the property (p) with respect to all $\{D_n\}$ over $I$ and in $F$ then
\[
B(x;\{D_n\}) \leq c(x;I;G),\;\;\;\;\;\;b(x;I;F)=c(x;I;G).
\]}

\noindent\textbf{(b)}
\textit{If $F$ also satisfies (iv) then $b(x;I;F) = c(x;I;G)$.}

\noindent
\textbf{(c)}
\textit{If $F=F_1$ then
\[
c(x;I;G) \leq \sigma(x) \leq 2c(x;I;G)
\]
where $\sigma(x)$ is the function defined in Henstock [15].}

\noindent
\textbf{For (a)}, norm$(D_n) < e$ when $n>n_0(e)$, so that if $x$ is a permanent interior point of $\{D_n\}$, and $n>n_0$,
\begin{eqnarray*}
\left| g(J_n) - g(J_{1,n})-g(J_{2,n}\right| & \leq & C(x;G;e);\;\;\;\;\mbox{ hence}\vt
B(x;\{D_n\}) & \leq & C(x;G;e),
\end{eqnarray*}
and this for every $e>0$, so that
\[
B(x;\{D_n\}) \leq c(x;I;G),\;\;\;\;\;\;b(x;I;F) \leq c(x;I;G)
\]
\textbf{In (b)}, given $e>0$, we have $C(x;G;e) \geq cx;I;G)$ so that there is a division $D$ over $I$ and in $G$, and with norm less than $e$, such that for $\delta>0$ we have
\[
\left|g(J) - g(J_1) - g(J_2)\right| > C(x;G;e) - \delta 
\geq c(x;I;G) - \delta.
\]
For $e = \delta = n^{-1}$ let $D$ be $D_n$. Then by (2.03), which uses (ii) and (iv), $\{D_n\}$ is in $F$.
\[
B(x;\{D_n\}) = \overline{\lim}_{n \rightarrow \infty} \left| g(J_n)-g(J_{1,n}) - g(J_{2,n}) \right| \geq \overline{\lim}_{n\rightarrow \infty} \left(c(x;I;G) - \frac 1n \right) = c(x;I;G).
\]
Hence $b(x;I;F) \geq c(x;I;G)$, so that (b) follows from (a).

\noindent
\textbf{For (c).} In the paper Henstock [15] the following definitions are given. For $x<y<z$ let $D(x;y;z)$ be the maximum of 
\begin{eqnarray*}
\left| g(x\mbox{---}z)- g(x\mbox{---}y) - g(y\mbox{---}z)\right| && (2^6\mbox{ alternatives) and}\vt
\left|g(x\mbox{---}y)
+ g(y\mbox{---}z) - g(x\mbox{---}y - g(y\mbox{---}z)) \right| && \left(\frac{15 \times 16}{2} \mbox{ non-trivial alternatives}\right).
\end{eqnarray*}
We put \[
\sigma(y) = \overline{\lim}_{x,z \rightarrow y} D(x;y;z).
\]
For $\{D_n\}$ let $J_n=x_n$---$z_n$. Then 
\[
B(y;\{D_n\}) \leq \overline{\lim}_{n \rightarrow \infty} D(x_n;y;z_n)\;\;\;\mbox{ so that  } \;\;\; b(y;I;F) \leq \sigma(y).
\]
But there are $x_n, z_n \rightarrow y$ such that 
\[
\sigma(y) = \lim_{n \rightarrow \infty} D(x_n;y;z_n)
\leq \overline{\lim}_{n \rightarrow \infty} 2 \max 
\left| g(x_n\mbox{---}z_n) - g(x_n\mbox{---}y_n) - g(y_n\mbox{---}z_n)\right|.
\]
Since $F=F_1$ we can find a $\{D'_n\}$ with $J_n=x_n$---$z_n$, and $J_n, J_{1,n}, J_{2,n}$ so that they give the maximum value to $|g(x_n\mbox{---}z_n) - g(x_n\mbox{---}y) - g(y\mbox{---}z_n)|$. Hence
\[
\sigma(y) \leq 2B(y;\{D'_n\}),\;\;\;\;\;\;\sigma(y) \leq 2b(y;I;F).
\]
Then (b) completes the proof of (c).

\vs
\noindent
\textbf{Corollary.}
\textit{When $F=F_1$, if $y$ is a singularity in the $\sigma(y)$ sense, $y$ is a singularity in the $b(y;I;F)=c(y;I;G)$ sense, and conversely.}

To connect these new functions with the integrals we have the following results. We put 
\[
\mbox{osc}(g;R;F) = (F) \overline{\int}_R g(I) - (F)\underline{\int}_{\,R} g(I) 
\]
whenever the right-hand side exists. Thus the two Burkill integrals cannot be infinite of the same sign.
Similarly
\[
\mbox{osc}(g;R;G) = (N;G) \overline{\int}_R g(I) - (N;G)\underline{\int}_{\,R} g(I) .
\]
\textbf{(3.04)}
\[
B(x;\{D_n\}) \leq \mbox{osc}(g;I;F),\;\;\;\;
b(x;I;F) \leq \mbox{osc}(g;I;F),\;\;\;\;
c(x;I;G) \leq \mbox{osc}(g;I;G).
\]
For $|g(J_n) - g(J_{1,n}) - g(J_{2,n})|
=|(I;D_n) \sum g(I) -(I;D_{1,n}+D_{2,n})\sum g(I)|$. Hence
\begin{eqnarray*}
B(x;\{D_n\}) &=& \overline{\lim}_{n \rightarrow \infty} \left| (I;D_n) \sum g(I) - (I;D_{1,n}+D_{2,n})\sum g(I)\right| \vt
&\leq & (F) \overline{\int}_I g(I) - (F)\underline{\int}_{\,I} g(I) \;\;\;\leq \;\;\;\mbox{osc}(g;I;F).
\end{eqnarray*}
Hence also $b(x;I;F) \leq \mbox{osc}(g;I;F)$. Similarly
\[
\left|g(J) - g(J_{1}) - g(J_{2})\right| =
\left|(I;D) \sum g(I) -(I;D_{1}+D_{2})\sum g(I)\right|.
\]
Then by (2.09a) and a similar result for lower limits, given $\ve>0$ there is a $\delta>0$ such that
\[
C(x;G;\delta) \leq (N;G) \overline{\int}_I g(I) - (N;G)\underline{\int}_{\,I} g(I) +2\ve.
\]
Hence $c(x;I;G) \leq \mbox{osc} N(g;I;G)$.

\vs
\noindent
\textbf{(3.05)}
\textit{There is a $g(I)$ with finite Burkill integrals such that if $x'<y<z'$ and $(x',z')$ is in $\bar S$ then 
\[
b(y;(x',z');F)=0=c(y;(x',z');G), \;\;\;\;\mbox{ but }\;\;\;
\overline{\lim}_{x,z\rightarrow y}\mbox{osc}(g;(x,z);F) \geq 2
\]
for $x<y<z$ and $(x,z) $ in $\bar S$.}

Let $x'<y<z'$ where $(x',z')$ is in $\bar S$. Then by (3.01b) there is a sequence $(x_1,z_1), \ldots , (x_m,z_m), \ldots$ of intervals of $\bar S$ and in $(x',z')$ such that \[
x_1< \cdots <x_m< \cdots,\;\;\;\;\;\;z_1> \cdots > z_m> \cdots ,
\]
 and $x_m, z_m \rightarrow y$ as $m \rightarrow \infty$.
 
For each $(x_m,z_m)$ choose a Riemann succession $\{D_{mn}\}$ in $F$. Let $m_i$ be the $i$th term in the sequence $1,1,2,1,2,3,1,2,3,4,1,\ldots$. We proceed by induction to define a sequence $n_1, \ldots , n_i, \ldots $ and numbers $u_0=z_1, u_1, u_2, \ldots$ as follows. Assuming that $n_1, \ldots , n_{i-1}, u_0, u_1, \ldots , u_{i-1}$ have been defined, $n_i$ is the smallest integer greater than $n_{i-1}$ such that $D_{m_i, n_i}$ has the point $u_i$ of division nearest to, but greater than $y$, with $u_i<u_{i-1}$.

We now put $f(x) = (-1)^{p_i}$ for $u_{i-1}>x\geq u_i$ and $f(z_1)=-1$, where $p_i$ is the $i$th term of the sequence $1,2,1,3,2,1,4,3,2,1,5,\ldots$. We then\footnote{If $I$ is an interval such as $[a,b)$ then $S(f;I)$ is the Stieltjes increment of $f$, $f([a,b)) = f(b)-f(a)$.---P.M.} put $g(I)=S(f;I)$ 
for $I'$ in $(y,z_i]$, $g(I) =0$ otherwise in $(x',z')$. Then
\[
\left(\left(x_{m_i},z_{m_i}\right);D_{m_i,n_i}\right) \sum g(I) =
f\left(z_{m_i}\right)-f\left(u_{i}\right) =
f\left(z_{m_i}\right) + (-1)^{p_i+1},
\]
and by construction of $p_i$, and for $m=1,2, \ldots$,
\begin{eqnarray*}
\overline{\lim}_{n \rightarrow \infty}
\left(\left(x_{m},z_{m}\right);D_{m,n}\right) \sum g(I)
& \geq & f(z_m) +1, \vt
\underline{\lim}_{n \rightarrow \infty}
\left(\left(x_{m},z_{m}\right);D_{m,n}\right) \sum g(I)
& \leq & f(z_m) -1, \vt
\mbox{osc}(g;(x_m,z_m);F) &\geq & 2, \vt
\underline{\lim}_{m \rightarrow \infty} \mbox{osc}(g;(x_m,z_m);F) &\geq & 2
\end{eqnarray*}
But $b(y;(x',z');F)=0=c(y;(x',z');G)$ since $g(I)=0$ when $y$ is in $I'$. Hence the results. This a ``singularity'' in an interval is no guarantee of a singularity of $g(I)$. 

Note that if $S$ satisfies (vii) below, then 
\[
b(y;F) =0=c(y;F)\;\;\;\;\mbox{ and }\;\;\;
\lim_{x,z \rightarrow y} \mbox{osc}(g;(x,z);F) \geq 2.
\]
At this point we introduce two new axioms.

\noindent
\textbf{(vii)}
If $R_1, R_2$ are in $\bar S$ and $R_1 \subset R_2$, then $R_2-R'_1$ is in $\bar S$.

This may be called an axiom of subtraction. $R_1$ is closed in order that $R_2-R'_1$ contains no isolated points.

\vs
\noindent
\textbf{(viii)}
If $I$ is in $S$ there is a division $D$ of $G$ and over $I$ which consists of $I$ alone.

\vs
\noindent
\textbf{(3.06)}

\noindent
\textbf{(a)}
\textit{If $S$ satisfies (vii) and if $R_1,R_2$ are in $\bar S$, so are $R_1 -R_2'$, $R_2 - R_1'$, $R_1$, $R_2$.}

\noindent
\textbf{(b)}
Also \textit{if $R_1 \subset R_2$, then osc$(g;R_1;F) \leq \mbox{osc}(g;R_2;F)$ whenever both sides exist, and osc$N(g;R_1;G) \leq \mbox{osc}N(g;R_2;G)$ whenever both sides exist. }

\noindent
\textbf{(c)}
\textit{Also $b(x;I;F) \equiv b(x;F)$ and $c(x;I;F) \equiv c(x;F)$ are independent of $I$.}

(a): By construction $R_1+R_2$ is in $\bar S$, and contains $R_1$. Hence by (vii), $R_1+R_2 -R_1' = R_2 - R_1'$ is in $\bar S$. Similarly $R_1-R_2'$ is in $\bar S$. And $R_1.R_2 \subset R$, so that substituting $R_2 - R_1'$ for $R_1$ in (vii), we have that $R_1.R_2$ is in $\bar S$.

(b): (2.11a) and (vii) give 
\[
\mbox{osc}(g;R_1;F) \leq \mbox{osc}(g;R_1+(R_2-R_1');F) = \mbox{osc}(g;R_2;F).
\]
Similarly fot osc$N(g;R;F)$.

(c): Let $I_1 \subset I_2$, where $I_1, I_2$ are in $\bar S$ and have a common end-point, and let $y$ be the other end-point of $I_1$. Then by (vii), $I_2-I_1'$ is in $\bar S$, which shows that $y$ satisfies (a) of property (p). Then (vi) implies that if $\{D_n\}$ is in $F$ and over $I_2$ we have a $\{D_n'\}$ in $F$ and over $I_1$, all of whose interval, save at most one in each division, are the intervals of the corresponding divisions of $\{D_n\}$.

Again, if $\{D_n'\}$ in $F$ is over $I_1$ there is by (vii) and (i) a $\{D_n''\}$ in $F$ and over $I_2-I_1'$, and so by (iii) a $\{D_n\}$ in $F$ and over $I_2$, such that all the intervals of $D_n'$ are intervals of $D_n$, for $n=1,2, \ldots$.

These results imply that if $x$ is inside $I_1$, then $b(x;I_1;F) = b(x;I_2;F)$. Similarly if $I_1'$ lies inside $I_2$.
We then have two points $y$ to consider. Now let $I_1,I_2$ be general intervals of $\bar S$. Then $I_1, I_2 \subset I_1+I_2$, so that
\[
b(x;I_1;F) = b(x;I_1+I_2;F) = b(x;I_2;F)
\]
for any point $x$ inside both $I_1$ and $I_2$. Thus $b(x;I;F)$ is independent of $I$. Similarly for $c(x;I;G)$.

\vs
\noindent
\textbf{Corollary.}
\textit{By (b) and (c) we can prove the note after (3.05):}

\vs
\noindent
\textbf{(3.07)}
\textit{If $F$ satisfies (vi), (vii), and if the upper and lower Burkill integrals of $g(I)$ over $R$ in $S$, with respect to the family $F$, are finite, then $g(I)$ has at most an enumerable number of singularities $y_1, \ldots ,y_n, \ldots $ in $R^o$, and
\[
\sum_{n=1}^\infty b(y_n;F) \leq \mbox{osc}(g;R;F).
\]
A similar result holds for norm-limits.}

For $n$ distinct points $y_1, \ldots , y_n$ inside $R$ there are $n$ overlapping intervals $I_1, \ldots , I_n$ of $\bar S$ and in $R$ such that $y_i$ is inside $I_i$ ($i=1, \ldots ,n$) (by (3.01)). From (3.04), using (2.12b) to show that all the osc$(g;I_i;F)$ are finite, and using (3.06c) for $b(y_i;F)$, we have
\[
\sum_{i=1}^n b(y_i;F) \leq \sum_{i=1}^n \mbox{osc}(g;I_i;F),\;\;< \mbox{osc}(g;R;F)
\]
by (vii) and (3.06b). The result then follows by an argument due originally to Cantor. For there are at most $m-1$ of the $y_i$ with $b(y_i;F) > \mbox{osc}(g;R;f)/m$, so that we may enumerate the $y_i$ by the size of unequal $b(y_i;F)$, and the position of $y_i$ for equal $b(y_i;F)$. And $\sum_{i=1}^\infty b(y_i;F) \leq \mbox{osc}(g;R;F)$.
For $F=F_1$, similar results have been given by Saks [10] 213 and Henstock [15] 205.

We now consider some very useful inequalities, assuming that the lower and upper Burkill integrals are finite. After (3.08) it is necessary to use norm-limits.

\vs
\noindent
\textbf{(3.08)}
\textit{Let $R_1, R_2$ in $\bar S$ be non-overlapping with sum $R_3$. Then }

\noindent
\textbf{(a)}
\[
(F)\overline{\int}_{R_3} g(I) \leq
(F)\overline{\int}_{R_1} g(I)
(F)\overline{\int}_{R_2} g(I)
+(FR_1.FR_2)\sum b(y;F),
\]
\textbf{(b)}
\[
(F)\underline{\int}_{\,R_3} g(I) \geq
(F)\underline{\int}_{\,R_1} g(I)
(F)\underline{\int}_{\,R_2} g(I)
-(FR_1.FR_2)\sum b(y;F).
\]
\textit{Similar results hold for norm-limits.}

Each point of $FR_1.FR_2$ satisfies (a) of property (p), so that by (vi) each has that property. Hence if $\{D_n\}$ in $F$ is over $R_3$ we obtain Riemann successions of divisions over the separate intervals of $R_3'$ (which successions are in $F$, by (v)), and then using (vi) and recombining by using (iii), we obtain $\{D_{i,n}\}$ in $F$ over $R_i$ ($i=1,2$). Then putting
\[
d_n(y) = \left| g(J_n) - g(J_{1,n})-g(J_{2,n}) \right|
\]
where the $J$'s are from $D_n$ for $y$,

\begin{eqnarray*}
(R_3;D_n)\sum g(I) &\leq & (R_1;D_{1,n})\sum g(I)+ (R_2;D_{2,n})\sum g(I)\;\; + \vt
&&\;\;\;\;\;\; + \;\;(FR_1.FR_2)\sum d_n(y), \vt
\overline{\lim}_{n\rightarrow \infty}(R_3;D_n)\sum g(I) &\leq & \overline{\lim}_{n\rightarrow \infty}(R_1;D_{1,n})\sum g(I) \;\; +\vt
&& \;\;\;\;\;\;+\;\;\overline{\lim}_{n\rightarrow \infty}(R_2;D_{2,n})\sum g(I)\;\;+\vt
&&\;\;\;\;\;\;+\;\; (FR_1.FR_2)\sum B(y;\{D_n\}),
\end{eqnarray*}
where we have extended the notation $B(y;\{D_n\})$ slightly.
The result (a) now follows. Similarly for (b) and for norm-limits.

\vs
\noindent
\textbf{Corollary:}
\textit{Using (2.11) also, we have}

\noindent
\textbf{(a)}
\[
(F)\overline{\int}_R (F)\overline{\int}_I g(I)
\leq (F)\overline{\int}_R g(I) \leq (F)\overline{\int}_R (F)\overline{\int}_I g(I) + (R^o)\sum b(y;F).
\]
\noindent
\textbf{(b)}
\textit{If $b(y;F) =0$ everywhere then the upper and lower integrals of $g(I)$ are additive.}

\noindent
\textbf{(c)}
\textit{The $b(y;F)$ for $(F)\overline{\int}_R g(I)$ is not greater than the $b(y;F)$ for $g(I)$.}

\vs
\noindent
\textbf{(3.09)}
\textit{Let the upper and lower norm-limits of $g(I)$ over $R$ in $\bar S$ with respect to $G$ be finite. Then given $\ve>0$ we can find $\delta>0$ such that every finite non-overlapping set $S_2$ of intervals of $S$ in $R$ with norm$(S_2<\delta$ gives
\[
(S_2) \sum g(I) < \ve + (N;G) \overline{\int}_{S_2}g(I) + (FS_2.R^o) \sum c(y;G).
\]
Similarly for lower limits. We have used $S_2$ both for the finite set of intervals and for the corresponding $I^\sigma$, and we suppose that $G$ satisfies (viii).}

Let $\delta>0$ be such that every division $D$ in $G$ and over $R$, with norm$(D)<\delta$, gives
\[
(R;D) \sum g(I) < (N;G) \overline{\int}_{R} g(I) + \frac \ve 2.
\]
This is possible by (2.09a). If $S_2$ has norm $<\delta$, and $R_1 = R-S_2'$, then $R_1$ is in $\bar S$ by (vii), and we can take a division $D_1$ of $R_1$ with norm $<\delta$, such that
\[
(R_1;D_1) \sum g(I) > (N;G) \overline{\int}_{R_1} g(I) - \frac \ve 2,
\]
by (2.09b). Now each interval of $S_2$ is in $S$, so that by (viii) and (iii), the intervals of $S_2$ form a division $D_2$ in $G$ over $S_2$. By (iii) again, $D_1+D_2$ is a division in $G$ over $R$, and norm$(D_1+D_2) < \delta$. Hence we have
\begin{eqnarray*}
(S_2)\sum g(I) &=& (R;D_1+D_2) \sum g(I) - (R_1;D_1) \sum g(I) \vt
&<& (N;G) \overline{\int}_{R} g(I) +\frac \ve 2 -(N;G) \overline{\int}_{R_1} g(I) + \frac \ve 2, \vt
&<& (N;G) \overline{\int}_{S_2} g(I) +  \ve  +(FS_2.R^o) \sum c(y;G)
\end{eqnarray*}
by (3.08a) for norm-limits.

Note that the effect of (viii) and (iii) together, is to ensure that every division (using intervals of $S$) of an $R$ in $\bar S$, is in $G$. This is a rather wide assumption, but it is necessary for (3.09). Axiom (v) and part of (vi) are then unnecessary for $G$.

\vs
\noindent
\textbf{(3.10)}
\textit{Let the norm-limit of $g(I)$ exist over $R$ in $\bar S$. Then given $\ve>0$ we can find $\delta>0$ such that for every finite family $S_2$ of non-overlapping intervals in $R$ with norm $<\delta$, we have}
\[
\left| (S_2) \sum \left(g(I) - (N;G) \int_I g(I) \right)\right| < \ve\;\;\;\;\mbox{ and }\;\;\;\;
 (S_2) \sum \left|\left(g(I) - (N;G) \int_I g(I) \right)\right| <2 \ve.
 \]
Since by (3.04) $c(y;G) =0$ everywhere in $R^o$ we obtain the first result from (3.09). But note that we can use (2.11) and omit (vi) and all mention of points, and still prove the results. For the second result take first the terms with \[g(I) - (N;G) \int_I g(I) \geq 0,\] and then the rest.

Saks (17) 167 proves the first result directly, when $F=F_1$.

\vs
\noindent
\textbf{(3.11)}
\textit{Let $\delta$ be as in (3.10). If for a division $D$ in $G$ over $R_1$ in $R$ we have norm$(D)<\delta$ and
\[
(R_1;D)\sum g(I) > (N;G) \overline{\int}_{R_1} g(I) - \ve 
\]
then every partial set $S_2$ of the set $S_3$ on intervals of $D$ is such that
\[
(S_2)\sum g(I) > -2\ve - (F(R_1-S_2).R^o) \sum c(y;G) + (N;G) \overline{\int}_{S_2} g(I).
\]
Similarly for lower integrals.
}

\begin{eqnarray*}
(S_2)\sum g(I) &=& (R_1;D) \sum g(I) - (S_3 - S_2)\sum g(I),
\;\;\;\;\;\mbox{
and by (3.09)},\vt
&>& (N;G) \overline{\int}_{R_1} g(I) - \ve -\ve -
(N;G) \overline{\int}_{R_1-S_2} g(I) \;\;- \vt
&&\;\;\;\;\;\;\;\;\;\;\;\;\;\;\;\;\;\;\;\;\;\;\;\;\;\;\;\;\;\;\;\;-\;\; (F(R_1-S_2).R^o) \sum c(y;G) \vt
&\geq &  (N;G) \overline{\int}_{S_2} g(I) - 2\ve
- (F(R_1-S_2).R^o) \sum c(y;G) 
\end{eqnarray*}
by (2.11). Hence the result.

\noindent
\textbf{(3.12)}
\textit{Let $x<y<z$, $I=x$---$z$, $I_1=x$---$y$, $I_2=y$---$z$ where $I_1, I_2$ are in $S$. Then, as $x,z$ tend independently to $y$,
\begin{eqnarray*}
\overline{\lim} \left(g(J) -(N;G) \overline{\int}_J g(I)\right) &=& 0 \;\mbox{ for } J=I_1, I_2;\;\mbox{ and either for }\;J=I,\;\mbox{ or else} \vt
\overline{\lim} \left(g(I_1)+g(I_2) -(N;G) \overline{\int}_I g(I)\right) &=& 0.\;\mbox{ Similarly}\vt
\underline{\lim} \left(g(J) -(N;G) \underline{\int}_{\,J} g(I)\right) &=& 0 \;\mbox{ for } J=I_1, I_2;\;\mbox{ and either for }\;J=I,\;\mbox{ or else} \vt
\underline{\lim} \left(g(I_1)+g(I_2) -(N;G) \overline{\int}_I g(I)\right) &=& 0.
\end{eqnarray*}
The alternative forms for $I$ are equal if $c(y;G)=0$ since $y$ satisfies (a) of property (p).}

If $J=I_1$, put $R_1$ as $(ay)$ in $\bar S$. If $J=I_2$ put $R_1$ as $(y,b)$ in $\bar S$. And if $J=I$ put $R_1$ as $(a,b)$. Then $F(R_1-J)$ does not include $y$, so that from (3.09), 
\[
\overline{\lim} \left(g(J) -(N;G) \overline{\int}_J g(I)\right) \leq 
(F(R_1-J).R^o) \sum c(s;G) \leq(R_1')\sum c(s;G) - c(y;G).
\]
Now $(R_1')\sum c(s;G)$ is convergent, and so tends to $c(y;G)$ as $a,b \rightarrow y$. Hence
$\overline{\lim} \left(g(J) -(N;G) \overline{\int}_J\, g(I)\right) \leq 0$. Similarly \[\overline{\lim} \left(g(I_1)+g(I_2) -(N;G) \overline{\int}_I\, g(I)\right) \leq 0.\]
We now use (3.11) in a similar proof of the opposite inequalities. When $J=I$ the special divisions in (3.11) determine which inequality is true. Note that if $I$ is not in $S$ we have the second of the alternative forms.

\noindent
\textbf{(3.13)}
\textit{If $g(I)$ is continuous, so is $(N;G)\int_I g(I)$ when it exists.}

\noindent
(Saks [17] 167, for $G=G_1$.)
Use (3.10). The following example shows that the upper and lower norm-limits need not be continuous even if $g(I)$ is continuous. Let
\[
f\left(\frac 1{2n}\right) =0,\;\;\;\;\;\;f\left(\frac 1{2n+1}\right) =1,
\]
and between these two values let $f(x)$ be linear. Let $g(I) =S(f;I)$ (\textit{the Stieltjes version of $f$---P.M.}) when $I=x$---$y$ with $x>0$ and $0<y-x<x^3$, and otherwise $g(I)=0$. Then it is easily seen that $g(I)$ is continuous, but 
\[
(N;G)\overline{\int}_{(0,x)}\,g(I) =f(x),
\]
which oscillates between 0 and 1.

The following result does not assume axioms (vi), (vii), but needs (viii).

\vs
\noindent
\textbf{(3.14)}
\textit{If the upper norm-limit of $g(I)$ over $R$ in $\bar S$, with respect to $G$, is finite, then given $\ve>0$ there is a $\delta>0$ such that if $D$ over $R$ and in $G$ is formed of $I_1, \ldots ,I_n$ and norm$(D)<\delta$ the for $s=1,2, \ldots ,n$,
\[
\sum_{i=1}^s (N;G)\overline{\int}_{I_i}\,g(I) +
\sum_{i=s+1}^n g(I_i) < (N;G)\overline{\int}_{R}\,g(I)+\ve.
\]
Similarly for lower norm-limits.
}

\noindent
(Saks [10] \textsection 4, Lemma, for $G=G_1$.)
Given $\ve>0$ there is a $\delta>0$ such that if $D$ in $G$ is over $R$ with norm$(D)<\delta$ then

\noindent
\textbf{(1)}
\[
(R;D) \sum g(I) < (N;G)\overline{\int}_{R}\,g(I) + \frac \ve 2.
\]
Now let $I_1, \ldots ,I_n$ form such a $D$, and $s$ be one of $1, \ldots ,n$. For each $I_i$ take a $D_i$ over $I_i$ and in $G$ with norm $<\delta$ and such that

\noindent
\textbf{(2)}
\[
(I_i;D_i) \sum g(I) > (N;G)\overline{\int}_{I_i}\,g(I)
-\frac \ve{2n}.
\]
By (viii) and (iii), $D_{(s)}$ is in $G$, where $D_{(s)}$ is formed of $D_1, \ldots , D_s, I_{s+1}, \ldots ,I_n$. Then by (1) and (2),
\begin{eqnarray*}
\sum_{i=1}^s (N;G)\overline{\int}_{I_i}\,g(I)+
\sum_{i=s+1}^n g(I_i) &<& 
(R;D_{(s)} \sum g(I) + \frac \ve 2 \vt
&<& (N;G)\overline{\int}_{R}\,g(I)+ \ve.
\end{eqnarray*}
Hence the result. Note that none of the $(N;G)\overline{\int}_{R}\,g(I)$ can be $+\infty$ (by (2.12) (Corollary) and (a)), and we have assumed that none is $-\infty$, for otherwise the result would be trivial.

\section{$k$-integration and the $\sigma$-limit}
\label{k-integration and the sigma-limit}
When $(F)\overline{\int}_{R}\,g(I)$ and $(F)\underline{\int}_{R}\,g(I)$ are finite for some $R$ in $\bar S$, $(F)\overline{\int}_{I}\,g(I)$ and $(F)\underline{\int}_{I}\,g(I)$ are interval functions over those $I$ in $\bar S$ which lie in $R$; and when the integrals $(F){\int}_{I}\,g(I)$ exist, (2.11a) shows that they are additive (in the sense of \textsection 1).
But in general, if osc$(g;R;F)>0$ the upper and lower integrals are not additive. For take $F=F_1$ and put $g(-a$---$a)=1$ for every $a>0$, and otherwise $g(I)=0$. Then $(F_1)\overline{\int}_{I}\,g(I)=1$ when the origin is inside $I$, and $=0$ when the origin is outside or at an end of $I$. Thus for $x=-1$, $y=0$, $z=+1$ in the definition of an additive function, (a) is false.

A further difficulty arises from the fact that when $g(I)$ is additive it does not necessarily follow that the upper and lower Burkill integrals over $I$ in $S$ are equal to $g(I)$. Suppose namely, that $g(I)$ is ``1 for the origin and 0 elsewhere.', i.e.~$g(a$---$b)=1$ whenever (1) $a<0<b$, or (2) $a=0<b$ with bracket $[$ at $a$, or (3) $b=0>a$ with bracket $]$ at $b$, and $g(I)=0$ otherwise. Then $g(I)$ is additive. But $(F)\underline{\int}_{I}\,g(I)=0 \neq g(I)$ when the origin is inside $I$, since the lowest limit-point is obtained by taking brackets $)$ $($ at the origin. 

The preceding two paragraphs are essentially those in Henstock [15] \textsection 2 (pages 206--207). It was then pointed out that in the theory of Riemann-Stieltjes integration, Pollard [5] suggested the use of permanent points, and Dienes [14] suggested the use of bracket conventions. We will make slight extensions of the definitions of $k$-conventions and $k$-successions of the paper Henstock [15].

We suppose that $F$ obeys axioms (i) to (vii), save possibly (iv) and (v). Then the function $b(y;F)$ and the singularities $y_1, \ldots , y_n, \ldots $ can be defined when the $(F)\overline{\int}_{R}\,g(I)$ and $(F)\underline{\int}_{\,R}\,g(I)$ are finite.

If in a Riemann succession in $F$, each singular point $y_i$ is a fixed point of the divisions, with a specified set of bracket conventions (the $k$-convention) after a finite number of divisions depending on $y_i$ then we call the succession a $k$-\textit{succession}.

The family of all $k$-successions in $F$ is denoted by $F_k$. When $F=F_1$ then $F_k$ is denoted by $F_{1k}$. The family $F_k$ does not necessarily exist, but we will suppose that $F_k$ exists and satisfies axiom (i). Then we have

\noindent  
\textbf{(1)}
\[
-\infty < (F)\underline{\int}_{\,R}\,g(I) \leq(F_k)\underline{\int}_{\,R}\,g(I)\leq 
(F_k)\overline{\int}_{R}\,g(I) \leq(F)\overline{\int}_{R}\,g(I) < + \infty.
\]
We now suppose that the upper and lower norm-limits are finite. The function $c(y;G)$ can also be defined, giving the singularities $y_1', y_2', \ldots $ where $c(y;G)>0$. Since by (3.03a), $b(y;F) \leq c(y;G)$, the set $y_1, y_2, \ldots$ is included in $y_1', y_2', \ldots $. we can quite easily have $b(y,F) =0$, $c(y;G)>0$, so that some $y_i'$ are sometimes not in $y_1, y_2, \ldots$.

For example, let the Riemann succession $\{D_n\}$ of divisions in $F$ be any Riemann succession over $R$ in $[-1,1]$, except that it includes for some integer $m$, the following arrangements about the origin. We have
\[
-\frac 1{m^i}\mbox{---}\frac 1{m^i},\;\mbox{ or }\;
-\frac 1{m^i}\mbox{---}0,\;\mbox{ or }\;
0\mbox{---}\frac 1{m^i},\;\;\;
(i=i_1, \ldots , i_n, \ldots\;\mbox{ tending to }\;+\infty).
\]
If $m$ is not a power of a smaller integer put

\[
g\left( \frac{-1}{m}\mbox{---}\frac 1m\right)=1,\;\;\;\;\;\;(m=1,2,3,5, \ldots).
\]
Otherwise put $g(I)=0$. Then
\[
B\left(0;\{D_n\}\right)
=\overline{\lim}_{n \rightarrow \infty} \left| 
g\left(-\frac 1{m^i}\mbox{---}\frac 1{m^i}\right)
- g\left(-\frac 1{m^i}\mbox{---}0\right)
- g\left(0\mbox{---}\frac 1{m^i}\right)
\right| =0.
\]  
Hence $b(0;F) =0$. But $C(0;G;e)=1$ so that $c(0;G)=1$.

At each $y_i'$ we take a set of bracket conventions, $k'$-\textit{convention}, such that if $y_i'$ is a $y_j$ the $k'$  and $k$-conventions there are the same. If in a Riemann succession in $F$ each singular point $y_i'$ is a fixed point of the divisions with correct $k'$-convention after a finite number of divisions depending on $y_i'$ then we call the succession a $k'$-\textit{succession}. The family of all $k'$-successions in $F$ is denoted by $F_{k'}$. Then $F_{k'} \subset F_k$. When $F=F_1$ then $b(y;F_1) =c(y;G_1)$, so that the $F_{k'}$ is $F_{1k}$. We will suppose that $F_{k'}$ exists and satisfies axiom (i). Then

\noindent
\textbf{(2)}
\[
(F_k)\underline{\int}_{\,R} g(I)
\leq (F_{k'})\underline{\int}_{\,R} g(I)
\leq (F_{k'})\overline{\int}_{\,R} g(I)
\leq (F_k)\overline{\int}_{\,R} g(I).
\]
Since axioms (ii), (iii) are true for $F$ they are true for $F_k$ and $F_{k'}$. Axiom (iv) in general is never true for an $F_k$ or an $F_{k'}$, since there are usually Riemann successions formed from divisions in $G$, which are not $k$-successions. In the present theory we are dealing with two parameters, norm$(D)$ and the number of permanent points included with the correct $k$-convention (or $k'$-convention) in the division $D$, so that the corresponding ``closure'' axioms are as follows.

\vs
\noindent
\textbf{(ix)} Let $\{D_n^{(i)}\}$ ($i=1,2, \ldots $) be a set of $k$-successions in $F_k$, of divisions of $R$ in $\bar S$. Then if
\[
D_1^{(1)},\;\;D_2^{(1)},\;\;D_1^{(2)},\;\;D_3^{(1)},\;\;D_2^{(2)},\;\;D_1^{(3)},\;\;D_4^{(1)},\dots 
\]
is a $k$-succession it is in $F_k$.

\vs
\noindent
\textbf{(x)}
For (x) replace $k$ by $k'$ throughout (ix).

\vs

If (iv) is true for $F$ then (ix) is true for $F_k$. If (ix) is true for $F_k$ then (x) is true for $F_{k'}$. These results are obvious. Also if (v) is true for $F$ it is true for $F_k$ and $F_{k'}$.

When the upper and lower norm-limits are finite, limits corresponding to them may be defined as follows. The \textit{upper $k$-limit $(k;G)\overline{\int}_R g(I)$ of $g(I)$ over $R$ in $\bar S$, with respect to the set $G$ of divisions}
is the greatest lower bound of numbers $dk(e;n)$ for all $e>0$ and integers $n$, where $dk(e;n)$ is the least upper bound of $(R;D)\sum g(I)$ for all $D$ in $G$ such that $D$ is a division of $R$ with (at least) the first $n$ singularities $y_1, y_2, \ldots$ included with the correct $k$-convention, and with norm$(D)<e$. As $e \rightarrow 0$ and $n \rightarrow \infty$, $dk(e;n)$ is monotone decreasing, and so the lower bound is also the limit.

A similar definition holds for the \textit{lower $k$-limit}
 $(k;G)\underline{\int}_{\,R} g(I)$. If for $R$ the upper and lower $k$-limits are equal, we say that the \textit{$k$-limit of $g(I)$ exists over $R$ in $\bar S$, with respect to the family $G$ of divisions}, and we write the common value as $(k;G) \int_R g(I)$.
 
 Replacing $k$ by $k'$ and $y_1, y_2, \ldots $ by $y'_1, y'_2, \ldots $, we obtain the definitions of the \textit{upper and lower $k'$-limits and the $k'$-limit}; respectively
 \[
 (k';G)\overline{\int}_{\,R} g(I),\;\;\;\;\;\;\;\;(k';G)\underline{\int}_{\,R} g(I), \;\;\;\;\;\;\;\;(k';G){\int}_{\,R} g(I),
 \]
\textbf{(3)}
\begin{eqnarray*}
&&-\infty <   (N;G)\underline{\int}_{\,R} g(I) 
\leq  (k;G)\underline{\int}_{\,R} g(I)
\leq  (k';G)\underline{\int}_{\,R} g(I) \leq \vt
&&\;\;\;\;\;\;\;\leq  (k';G)\overline{\int}_{R}\, g(I)
\leq  (k;G)\overline{\int}_{R}\, g(I)
\leq  (N;G)\overline{\int}_{R}\, g(I)
<\infty.
\end{eqnarray*}
The theorems (2.05), (2.07), (2.08a) are true for $F=F_k$ or $F_{k'}$. Te remaining part of (2.06) is true in the forms

\vs
\noindent
\textbf{(4.01)}

\noindent
\textbf{(b)}
\textit{If $(k;G){\int}_{\,R} g(I)$ exists the given $\ve>0$ there are a $\delta>0$ and an integer $n$ such that if $D$ in $G$ is a division of $R$ with norm$(D)<\delta$ and the first $n$ singularities $y_1,y_2, \ldots$ included with the correct $k$-convention, then}
\[
\left|(R;D)\sum g(I) - (k;G){\int}_{R}\, g(I)\right| < \ve.
\]
\textbf{(b)}
\textit{We can replace $k$ by $k'$; $y_1,y_2, \ldots $by $y_1', y_2', \ldots $, in (a).}

Theorem (2.08b) is true for $F_k, F_{k'}$ in the form

\vs
\noindent
\textbf{(4.02)}

\noindent
\textbf{(a)}
\textit{If $R$ is in $\bar S$ and $F_k$ satisfies (ix) then $F_k$ contains $\{D_n\}$ over $R$ with}
\[
\lim_{n \rightarrow \infty}
(R;D_n)\sum g(I) = (F_k)\overline{\int}_R \,g(I).
\]
\textbf{(b)}
\textit{We can replace $k$ by $k'$, (ix) by (x), in (a).}

Theorem (2.09) is altered in much the same way as (2.06) was altered to become (4.01). Theorem (2.10) becomes

\noindent
\textbf{(4)}
\[
(k;G)\underline{\int}_{\,R} g(I) 
\leq   (F_{k})\underline{\int}_{\,R} g(I)
\leq  (F_k)\overline{\int}_{\,R} g(I) 
\leq  (k;G)\overline{\int}_{R}\, g(I);
\]

\noindent
\textbf{(5)}
\[
(k';G)\underline{\int}_{\,R} g(I) 
\leq   (F_{k'})\underline{\int}_{\,R} g(I)
\leq  (F_k')\overline{\int}_{\,R} g(I) 
\leq  (k';G)\overline{\int}_{R}\, g(I).
\]
\noindent
\textbf{(6)}
\textit{If $F_{k'}$ satisfies (x) then}
\[
(k';G)\underline{\int}_{\,R} g(I) 
=   (F_{k'})\underline{\int}_{\,R} g(I)
=  (F_k')\overline{\int}_{\,R} g(I) 
=  (k';G)\overline{\int}_{R}\, g(I).
\]
\noindent
\textbf{(7)}
\textit{If $F_{k'}$ satisfies (x) then
}
\[
(k;G)\underline{\int}_{\,R} g(I) 
=   (F_{k})\underline{\int}_{\,R} g(I),\;\;\;\;\;\;
 (F_k)\overline{\int}_{\,R} g(I) 
=  (k;G)\overline{\int}_{R}\, g(I).
\]
Theorem (2.11a) is true for $F=F_k$ or $F_{k'}$. Of course, the parts (b) and (c) do not apply here. If $F$ satisfies (v) then so do $F_k$ and $F_{k'}$, as already noted, and Theorem (2.13a) is true for $F_k$ and $F_{k'}$ also. Similarly Theorems (2.14) to (2.20) are true for $F=F_k$ or $F_{k'}$.

The theorems mentioned in the last paragraph are true also for $k$ and $k'$-limits.

If the $k$-convention and $k'$-convention are all conventions then (vi) is true for $F_k$ and $F_{k'}$ since it is true for $F$. Otherwise we need in general to postulate that $F_k$ and $F_{k'}$ satisfy (vi). Since $F$ satisfies (vii) so do $F_k$ and $F_{k'}$. Obviously
\[
b(x;F_{k'})\leq b(x;F_{k}) \leq b(x;F)=0
\]
if $x$ is not a $y_i$ we may define $ck(x;I;G)$ corresponding to $c(x;I;G)$, as 
\[
\lim_{e\rightarrow 0,\;n \rightarrow \infty} Ck(x;G;e;n)\;\;\mbox{ where }\;\;Ck(x;G;e;n) = \mbox{l.u.b.}\left|g(j)-g(J_1)-g(J_2)\right|
\]
(as in property (p)) for every triad $J,J_1,J_2$ of intervals from a division $D$ (over $I$ and in $G$) with norm$(D)<e$ with (at least) the first $n$ singularities $y_1,y_2, \ldots $ included with the correct $k$-convention, and with $x$ inside the interval $J$ of $D$.

Replacing $k$ by $k'$ and $y_1,y_2, \ldots$ by $y_1', y_2', \ldots$, we obtain the definition of $ck'(x;I;G)$. When $x$ is not a $y_i'$ we have 
\[
ck'(x;I;G) \leq ck(x;I;G) \leq c(x;I;G)=0
\;\mbox{ so that }\; ck'(x;I;G)=0=ck(x;I;G).
\]
When $x$ is not a $y_i$ we have
$ck(x;I;G) \leq c(x;I;G)$, but the latter is not necessarily 0. In fact, for the example given in which $b(0;F)=0$, $c(0;G)=1$ we have $ck(0;I;G)=1$ also.

Since $b(x;F_k) =0=b(x;F_{k'})$ when $x$ is not a $y_i$, and  $ck'(x;I;G)=0$ when $x$ is not a $y_i'$, the analogues of (3.08) and (2.11a) for $F=F_k$ and $F_{k'}$, and for $k'$-limits, give (4.03), as in (3.08 Corollary (b)).

\vs
\noindent
\textbf{(4.03)}
\textit{The upper and lower Burkill integrals with respect to $F_k$ and with respect to $F_{k'}$, and the upper and lower $k'$-limits are each additive for $R$ in $\bar S$.}
\noindent (Henstock [15] 207, Theorem (3.1), for $F_k=F_{1k}$.)

\vs
A division in $G$ over $R$, with norm $<e$, and with the correct $k$-convention at each of the $y_1, \ldots , y_n$, is called a \textit{$k_{ne}$-division of $R$}.

A division in $G$ over $R$, with norm $<e$, and with the correct $k'$-convention at each of the $y_1', \ldots , y_n'$ (at least), is called a \textit{$k'_{ne}$-division of $R$}. Similarly for a $k'_{ne}$-\textit{set} of a finite number of intervals. 
 (Henstock [15] 208, for $G=G_1$.)

\noindent
\textbf{(8)}
\textit{If the $k'$-convention at each $y_i'$ which is not a $y_j$, is all conventions, then}
\[
(k';G)\underline{\int}_{\,R} g(I) 
\leq   (F_{k})\underline{\int}_{\,R} g(I)
\leq  (F_k)\overline{\int}_{\,R} g(I) 
\leq (k';G)\overline{\int}_{R}\, g(I).
\]
Let $y_1', \ldots , y_n'$ divide $R$ up into $R_1, \ldots , R_q$, all of which are in $\bar S$ by repeated applications of (vi) and (vii). Given $e, e'>0$ there is by (2.08a) for $F_k$, a $k_{ne}$-division $D_i$ of $R_i$ such that
\[
(R_i;D_i) \sum g(I) >(F_k) \overline{\int}_
{R_i} g(I) - \frac{e'}q,\; \;\;\;(i=1, \ldots ,q).
\]
Then by (4.03) the division $D=\sum_{i=1}^q D_i$ of $R$ gives
\[
(R;D) \sum g(I) >(F_k) \overline{\int}_
{R} g(I) - {e'}.
\]
But $D_i$ is in $G$, so that by (2.01), $D$ is in $G$. Thus $D$ is a $k'_{ne}$-division of $R$, so that by (2.01), $D$ is in $G$. Thus $D$ is a $k'_{ne}$-division of $R$, so that
\begin{eqnarray*}
dk'(e;n) &>&(F_k) \overline{\int}_{R} g(I) - {e'}, \vt
dk'(e;n) & \geq & (F_k) \overline{\int}_{R} g(I) ,\vt
(k';G) \overline{\int}_{R} g(I) 
&\geq &
(F_k) \overline{\int}_{R} \,g(I).
\end{eqnarray*}
Similarly for lower $k'$-limits.

\vs
\noindent
\textbf{(4.04)}
\textit{If, keeping the same $k'$-convention at each of the $y_i'$, we add more permanent points, each of which satisfy (a) of property (p), we do not alter the $k'$-limits.}

\noindent
(Henstock [15] 208, Theorem (3.2), for $G=G_1$.)

Let the first $n$ extra points divide $R$ up into $R_1, \ldots ,R_q$. Then the proof proceeds as in (8), using 
$(k';G) \overline{\int}_{R_i}\, g(I) $ instead of $(F_k) \overline{\int}_{R_i} \,g(I)$, and $k'_{ne}$-divisions instead of $k_{ne}$-divisions. Thus the new upper integral is not less than $(k';G) \overline{\int}_{R}\, g(I) $.

The opposite inequality is obvious. Hence the results. The conventions at the extra points do not matter since $c(y;G)=0$. 

At this point it would be useful to collect together the results of inequalities (1) to (8). Though we only give results for upper limits and integrals, there are corresponding results for lower limits and integrals.

We denote $(A) \overline{\int}_{R}\, g(I) $ or $(A;G) \overline{\int}_{R}\, g(I) $ symbolically by $A$, and $\leq$ by $\rightarrow$. Then (1) to (8) give

\vs
\noindent
\textbf{(4.05)}
\[
\begin{array}{lllllll}
   &   &F_k&\rightarrow &F &    &    \vt
   & \nearrow &   &   &   & \searrow &   \vt
F_{k'}&   &   & \searrow &   &   & N\;\;(<+\infty) \vt
   & \searrow &   &   &    & \nearrow &   \vt
   &   & k'& \rightarrow & k &   &      
\end{array}
\]

\noindent
\textbf{(4.05)}
\textit{If the $k'$-convention at each $y_i'$ which is not a $y_i$, \textbf{all} conventions then}
\[
\begin{array}{lllllllll}
   &   &   &   &k' & \rightarrow & k&   &   \vt
   &   &   &\nearrow &   &   &   &\searrow &   \vt
 F_{k'} & \rightarrow & F_k &   &   &   &   &   &N \vt
    &   &   &   &\searrow & &  \nearrow &   &   \vt
 &&&&&F&&&
\end{array}
\]
\textbf{(4.07)}
\textit{If $F_{k'}$ satisfies (x) then}
\[
\begin{array}{lllllllll}
   &   &  &   &   &   & F &   & \vt
   &&&&& \nearrow &&\searrow& \vt
F_{k'}&=&k'&\rightarrow &F_k&&&&N \vt
&&&&&\searrow && \nearrow &\vt
&&&&&&k&&
\end{array}
\]
\textbf{(4.08)}
\textit{If $F_K$ satisfies (ix) then}
\[
F_{k'} = k' \rightarrow F_k = k \rightarrow F \rightarrow N
\]
\textbf{(4.09)}
\textit{If $F_{k'}$ satisfies (x) and if the $k'$-convention at each $y_i'$ which is not a $y_j$, is \textbf{all} conventions then}
\[
{\sf\mbox{\textit{{\sf{(Chart missing in thesis---P.M.)}}}}}
\]
\textbf{(4.10)}
\textit{If $f$ satisfies (ix) and if the $k'$-convention at each $y_i'$ which is not a $y_j$ is \textbf{all} conventions then
\[
F_{k'}=F_k=k'=k \rightarrow F \rightarrow N.\]}

The following result is almost obvious.

\noindent
\textbf{(4.11)}
\textit{If $G$ is the family of divisions of every $R$ in $W$, which only use the conventions $)[$ or $[($, with a fixed convention at each boundary point of $R$; if $g(I)$ is additive; and if $I_1$ has the correct fixed convention at its ends; then $(N;G)\int_{I_1}g(I)$ exists and is equal to $g(I_1)$.} (Henstock [15] 208, Theorem (3.3).)

\noindent If $D$ is n $G$ and $I_1$ has the correct fixed conventions at its ends then \[(I_1;D)\sum g(I) = g(I_1).\] Hence the result.

The following are analogues of the results (3.09), (3.11) and (3.12).

\vs
\noindent
\textbf{(4.12)}
\textit{Given $\delta>0$ there is an $e>0$ and an integer $n$ so that every $k'_{ne}$-set $S_2$ in $R$ gives
\[
(S_2) \sum g(I) < (k';G)\overline{\int}_{S_2} \,g(I) + \delta,
\]
whenever $G$ satisfies (viii). Similarly for lower integrals.}

\noindent (Henstock [15] 208, Theorem (3.4), for $G=G_1$.)

\vs
\noindent
\textbf{(4.13)}
\textit{Let $n, e$ be as in (4.12). If for a $k'_{ne}$-division $D$ in $G$ over $R_1$ in $R$ we have
\[
(R_1;D)\sum g(I) > (k';G)\overline{\int}_{S_2} \,g(I) -\delta
\]
then every  partial set $S_2$ of the set $S_3$ of intervals of $D$ is such that 
\[
(S_2)\sum g(I) > (k';G)\overline{\int}_{S_2} \,g(I) -2\delta
\]
Similarly for lower integrals.}
(Henstock [15] 209, Theorem (3.5), for $G=G_1$.)

\vs
\noindent
\textbf{(4.14)}
\textit{
Let $x<y<z$, $I=x$---$z$, $I_1=x$---$y$, $I_2=y$---$z$, where $I_1,I_2$ are in $S$. Then as $x,z$ tend independently to $y$,
\[
\overline{\lim}\left(g(J) -(k';G)\overline{\int}_J\,g(I)\right) =0
\]
for $J=I_1, I_2$, and $I$ if $I$ is in $S$, provided that if $y$ is a singularity $y_i'$ then $J=I_1$ or $I_2$ alone, with the correct $k'$-convention at $y$.}

\noindent (Henstock [15] 209, Theorem (3.6), for $G=G_1$.)

Let $I^k$ denote an interval $I_1$ with $I_1'=I'$ and with the correct $k$-convention (i.e.~one of the fixed set of bracket conventions) at an end if a $y_i$ occurs there. A similar definition can be given of $I^{k'}$.

\vs
\noindent
\textbf{(4.15)}
\textit{Let $h(I) = g(I^k)$. Then for $F=F_1$, $(F_k)\overline{\int}_R\,h(I)\leq (F_k)\overline{\int}_R\,g(I)$
and we can have $(F_k)\overline{\int}_R\,h(I) < (F_k)\overline{\int}_R\,g(I)$ for some $g$.}

\noindent
Let $x_r=2^{-r}$, $a(I) =1$ when $I=[0,x_r]$, ($r=1,2,\ldots$), $a(I) =0$ otherwise, and 
\[
g(I) = (I) \sum x_r + a(I).
\]
Let the $k$-convention be $)[$ at $0$ and $x_1,x_2, \ldots$. Then $a(I^k)=0$ and
\[
(F_{1k})\overline{\int}_{90,1)}\,h(I) = \sum_{r=1}^\infty \frac 1{2^r} =1,\;\;\;\mbox{ whereas }\;\;\;(F_{1k})\overline{\int}_{(0,1)}\,g(I)=2.
\]
Thus we could define a new set of integrals by using $g(I^k)$ and $g(I^{k'})$, but their relations with $(F_k)\overline{\int}_R\,g(I)$ etc.~are not interesting, and so have been disregarded.

The case where we take \textbf{all} conventions at each $y_i'$ can be connected with the $\sigma$-limit of Getchell [11] and Hildebrandt [13] over $R$ in $\bar S$.

Given $R$ in $\bar S$, if there is a finite number $A$ such that for every $\ve>0$ we can find a division $D$ in $G$ over $R$ such that $|(R;D')\sum g(I) -A|<\ve$ whenever $D'$ in $G$ and over $R$ contains all the division points of $D$, then we call $A$ the $\sigma$-\textit{limit} of $g(I)$ over $R$ and with respect to $G$, and write $A$ as $(\sigma;G)\int_R\,g(I)$.

Getchell [11] and Hildebrandt [13] take the simple case $G=G_1$. If $(k';G){\int}_R\,g(I)$ exists when we take all conventions at each $y_i'$, then by (4.01b), for each $\ve>0$ we can find $\delta>0$, an integer $n$, and such a $D$ (in $G$ and over $R$) by including in $D$ the first $n$ singularities $y_1', \ldots , y_n'$ and then adding other points so that norm$(D)<\delta$. Hence in this case $(\sigma;G)\int_R g(I)$ exists and is equal to $(k;G){\int}_R\,g(I)$.

If $(\sigma;G)\int_R\,g(I)$ exists and the upper and lower norm-limits of $g(I)$ over $R$, with respect to $G$ are finite, then we define $\{y_i'\}$ and $(k';G)\overline{\int}_R\,g(I)$, $(k';G)\underline{\int}_R\,g(I)$,
where $k'$ is \textbf{all} conventions. Then given $\ve_1, \ve_2, \ldots $ tending to 0 let the corresponding divisions $D$ (in $G$ and over $R$) in the definition of the $\sigma$-limit be $D_1, D_2, \ldots $, and take an arbitrary sequence $n_1, n_2, \ldots $ of integers tending to $+\infty$. We add $y_1', \dots , y_{n_i}'$ to $D_i$ to form $D_i'$, which by (vi) and (iii) is in $G$, and then use $\{D_i'\}$ instead of $\{D_i\}$ without altering the $\sigma$-limit $A$, so that by (4.04) we see that
\[
(k';G)\overline{\int}_R\,g(I) = A =
(k';G)\underline{\int}_R\,g(I),\;\;\;\;\mbox{ i.e.,}
\]
\textbf{(4.16)}
\textit{Let $k'$ be \textbf{all} conventions, taking finite upper and lower norm-limits of $g(I)$ over $R$ in $\bar S$, with respect to $G$. Then the existence of the $\sigma$-limit of $g(I)$ over $R$ with respect to $G$ implies and is implied by the existence of the corresponding $k'$-limit, and the two limits are equal.}

We now consider suitable definitions of upper and lower $\sigma$-limits. Let $D$ be over $R$ and in $G$, and let $G.D$ be the family of all divisions in $G$ which use all the division-points of $D$ (at least) which lie in the $I^\sigma$ which the divisions cover. By (i) and (iii), $G.D$ has divisions over $R$ of arbitrarily small norm. Then the \textit{upper $\sigma$-limit} $(\sigma;G)\overline{\int}_R\,g(I)$, \textit{of $g(I)$ over $R$ in $\bar S$, with respect to $G$} is the greatest lower bound for all $D$ (in $G$ and over $R$) of $(N;G.D)\overline{\int}_R\,g(I)$. The \textit{lower $\sigma$-limit} $(\sigma;G)\underline{\int}_R\,g(I)$ is the least upper bound of
$(N;G.D)\underline{\int}_R\,g(I)$.

\vs
\noindent
\textbf{(4.17)}
\textit{There is a set $s_1, s_2, \ldots $ of points such that if they are used to produce upper and lower ``$s$-limits'' over $R$ in the same way as the upper and lower $k'$-limits are produced (for $k'$=all conventions) then the upper and lower $s$-limits are equal respectively to the upper and lower $\sigma$-limits.}

\noindent
For there is a sequence $\{D_n\}$ of divisions over $R$ and in $G$, such that
\begin{eqnarray*}
\lim_{n\rightarrow \infty}(N;G.D_n)\overline{\int}_R\,g(I) &=&(\sigma;G)\overline{\int}_R\,g(I)\;\;\mbox{and }\vt
\lim_{n\rightarrow \infty}(N;G.D_n)\underline{\int}_R\,g(I) &=&(\sigma;G)\underline{\int}_R\,g(I).
\end{eqnarray*}
Let $D'_n$ be the division consisting of all the division-points of $D_1, \ldots , D_n$. Then by (vi) and (iii), $D_n'$ is in $G$. And $\{D_n'\}$ can replace $\{D_n\}$. Then we take $s_1, s_2, \ldots$ to be first the division-points of $D_1'$ in order of increasing value, followed by those of $D_2'$ which are not in $D_1'$, in order of increasing magnitude, and so on. Let $D_n''$ be given by $s_1, \ldots ,s_n$. If $s_n$ is not in $D_i'$ then
\[
(\sigma;G)\overline{\int}_R\,g(I)
\leq (N;G.D_n'')\overline{\int}_R\,g(I)
\leq (N;G.D_i')\overline{\int}_R\,g(I).
\]
Since $i$ can tend to $+\infty$ as $n$ does so, $\{D_n''\}$ can replace $\{D_n'\}$ and so $\{D_n\}$. Then we have an $e>0$ such that (as for the norm-limit))
\[
ds(e,n) - \frac 1m <(N;G.D_n'')\overline{\int}_R\,g(I)
\leq ds(e,n)
\]
since $ds(e,n)$ is the $d(e)$ for $G.D_n''$. Hence
\begin{eqnarray*}
(s;G)\overline{\int}_\R\,g(I) &=&
\lim_{e \rightarrow 0, \;n \rightarrow \infty} ds(e,n) \vt
&=& \lim_{n\rightarrow \infty} (N;G,D_n'')
\overline{\int}_\R\,g(I) \vt &=&
(\sigma;G)\overline{\int}_\R\,g(I).
\end{eqnarray*}
Similarly for the lower ``$s$-limit''.

\vs
\noindent
\textbf{Corollary.}
\textit{We may add any enumerable set of points to $\{s_n\}$ without altering the ``$s$-limits''.}

\vs
\noindent
\textbf{(4.18)}
\textit{The $\sigma$-limit (of $g(I)$ over $R$ in $\bar S$, with respect to $G$) being equal to $A$, implies and is implied by the upper and lower $\sigma$-limits being equal to $A$.}
\noindent
Let \[(\sigma;G)\overline{\int}_R\,g(I) =A=
(\sigma;G)\underline{\int}_R\,g(I).\]
Then for an ``$s_{ne}$-division'' $D$ of $R$ with $n>n_0(\delta)$, $e < e_0(\delta)$, we have by (4.17),
\[
A-\delta < (R;D')\sum g(I) < A+\delta
\]
when $D'$ has all the division points of $D$. I.E.~the $\sigma$-limit exists equal to $A$. 

\noindent
If the $\sigma$-limit exists equal to $A$, let $\ve_n\rightarrow 0$ as $n \rightarrow \infty$, and let $\{D_{1,n}\}$ be the corresponding divisions in the definition of the $\sigma$-limit. Let $D_{2,n}$ be formed from the division-points of $D_{1,n}$ and $s_1, \ldots ,s_n$. Then we can replace $\{D_{1,n}\}$ by $\{D_{2,n}\}$, so that
\[
A-\ve_n \leq (N;G.D_{2,n})\underline{\int}_R\,g(I) \leq
(N;G.D_{2,n})\overline{\int}_R\,g(I) \leq A+\delta.
\]
Hence 
\[(\sigma;G)\overline{\int}_R\,g(I) =A=
(\sigma;G)\underline{\int}_R\,g(I).\]
Thus the upper and lower $\sigma$-limits are rightly named. They are of little use unless for some $D$ over $R$ and in $G$, 
\[
(N;G.D)\underline{\int}_R\,g(I) \;\;\;\mbox{ and }\;\;\;
(N;G.D)\overline{\int}_R\,g(I)
\]
are both finite. In this case we can define $c(y;G.D)$ and the corresponding $\{y_i'\}$, and so 
\[
(k';G.D)\underline{\int}_R\,g(I) \;\;\;\mbox{ and }\;\;\;
(k';G.D)\overline{\int}_R\,g(I)
\]
We then have

\vs
\noindent
\textbf{(4.19)}
\[
(\sigma;G)\underline{\int}_R\,g(I) =
(k';G.D)\underline{\int}_R\,g(I),\;\;\;\;\;\;
(\sigma;G)\overline{\int}_R\,g(I) =
(k';G.D)\overline{\int}_R\,g(I)
\]
\textit{where $k'$ is all conventions.}

\noindent
Take $s_1, s_2, \ldots $ equal to $y_1', y_2', \ldots $ and the points of division of $D$. By (4.04), any other points are useless. We then have the results.

The definitions of the upper and lower $\sigma$-limits give no construction to find $\{D_n\}$ and so $\{s_n\}$, so that the definitions are not as good as those of the upper and lower $k'$-limits with respect to $G.D$ whenever the norm-limits with respect to $G.D$ are \textit{known to be finite}. We can also use sets of bracket conventions for the upper and lower $k'$-limits, and can consider $k$-limits and integrals with respect to $F_k$ and $F_{k'}$. Thus the idea of upper and lower $k'$-limits is usually better than the idea of upper and lower $\sigma$-limits.

There are analogues of (2.06), (2.09). (2.11) becomes

\vs
\noindent
\textbf{(4.20)}
\textit{The upper and lower $\sigma$-limits are additive.}

\noindent
Let $R =R_1+R_2$ where $R_1,R_2$ are non-overlapping and in $\bar S$, and let $s_1,s_2, \ldots $ be $FR_1.FR_2$ together with the sequences for $R,R_1,R_2$. Then $\{s_n\}$ by (4.17 \& Corollary) will do for $R, R_1,R_2$. Then by (iii) and (vi) repeatedly, the divisions in $R_1,R_2$ are independent and form a division of $R$; and when $FR_1.FR_2$ is included, each division of $R$ forms a division of $R_1$ and a division of $R_2$. We thus have
\[
(\sigma;G)\overline{\int}_R\,g(I) =
(\sigma;G)\overline{\int}_{R_1}\,g(I) +
(\sigma;G)\overline{\int}_{R_2}\,g(I) 
\]
with the convention here that $+\infty-\infty = -\infty+\infty =+\infty$.

A similar result holds for the lower $\sigma$-limit, with
$+\infty-\infty = -\infty+\infty =+\infty$.

There are also analogues of (2.14) to (2.20).

\vs
\noindent
\textbf{(4.21)}
\[
(\sigma;G)\overline{\int}_{R }\,g(I) 
\geq (N;G)\overline{\int}_{R}\,g(I)
- (R^o)\sum c(s_i;G)
\]
\textit{whenever the right-hand side has meaning. Similarly for lower $\sigma$-limits.}

\noindent
Let $s_1, \ldots ,s_n$ divide $R$ up into $R_1, \ldots, R_q$. Then there is an $s_{ne}$-division $D_i$ of $R_i$ such that
\[
(R_i;D_i)\sum g(I) > (N;G) \overline{\int}_{R_i }\,g(I) ,\;\;\;\;\;(i=1, \ldots ,q)
\]
Put $D=\sum_{i=1}^q D_i$ so that by (iii), $D$ is in $G$, and therefore an $s_{ne}$-division of $R$. Then for $n>n_0(\delta)$, $e< e_0(\delta)$,
\begin{eqnarray*}
(\sigma;G)\overline{\int}_{R }\,g(I) + \delta &>& (R;D) \sum g(I) \vt
&>& \sum_{i=1}^q (N;G)\overline{\int}_{R_i }\,g(I),\;\;\;\;\;\mbox{ and by (3.08)},\vt
&\geq& (N;G)\overline{\int}_{R }\,g(I) -(R^o)\sum c(y;G). 
\end{eqnarray*}
Hence the result since $\delta>0$ is arbitrary.

\vs
\noindent
\textbf{Corollary.}
\textit{If $k'$ is \textbf{all} conventions then}
\[
(k';G)\overline{\int}_{R }\,g(I)\geq
(N;G)\overline{\int}_{R }\,g(I) -(R^o)\sum c(y;G).
\]  

\section{Existence theorems.}\label{Existence theorems}
\textbf{(5.1)}
\textit{The necessary and sufficient condition that $(F){\int}_{R }\,g(I)$ exists, is that for each different $\{D_n\}$ of $F$ and over $R$, and given $\ve>0$, there is an $n_0=n_0(\ve;\{D_n\};\{D_n'\})$ such that for $m,n>n_0$,
}
\[
\left|(R;D_m)\sum g(I) - (R;D_n')\sum g(I)\right|<\ve.
\]
The necessity is obvious. For sufficiency let $n \rightarrow \infty$. Then 
\[
\left|(R;D_m)\sum g(I) - \overline{\lim}_{n \rightarrow \infty}(R;D_n')\sum g(I)\right|\leq \ve,
\]
so that if $l>n_0$,
\[
\left|(R;D_m)\sum g(I) - (R;D_l)\sum g(I)\right|\leq 2\ve.
\]
Hence by the General Principle of Convergence,
$
\lim_{m \rightarrow \infty} (R;D_m) \sum g(I)$ exists and is finite. Similarly
$
\lim_{n \rightarrow \infty} (R;D_n') \sum g(I)$ exists and is finite. And then the two limits are equal.

\vs
\noindent
\textbf{(5.2)}
\textit{The necessary and sufficient condition that $(N;G)\int_{R }\,g(I)$ exists, is that given $\ve>0$ there is a $\delta>0$ such that each pair $D_1,D_2$ of divisions of $R$ with norm less than $\delta$ are such that
\[
\left|(R;D_1)\sum g(I) - (R;D_2)\sum g(I)\right| <\ve.
\]}
\textbf{(5.3)}
\textit{When the upper and lower Burkill integrals of $g(I)$ over $R$ in $\bar S$, with respect to $F$, are finite, the necessary and sufficient condition that $(k;G)\int_{R }\,g(I)$ exists, is that given $\delta>0$ there is an $e>0$ and an integer $n$ such that every pair $D_1,D_2$ of $k_{ne}$-divisions of $R$ satisfies
\[
\left|(R;D_1)\sum g(I) - (R;D_2)\sum g(I)\right| <\delta.
\]}

\noindent
\textbf{(5.4)}
\textit{When the upper and lower norm limits of $g(I)$ over $R$ in $\bar S$, with respect to $G$, are finite, the necessary and sufficient condition that $(k';G)\int_{R }\,g(I)$ exists, is that given $\delta>0$ there is an $e>0$ and an integer $n$ such that every pair $D_1,D_2$ of $k'_{ne}$-divisions of $R$ satisfies
\[
\left|(R;D_1)\sum g(I) - (R;D_2)\sum g(I)\right| <\delta.
\]}

\noindent
\textbf{(5.5)}
\textit{The necessary and sufficient condition that the $\sigma$-limit of $g(I)$ exists over $R$ in $\bar S$, with respect to $G$,
is that for some $D$ (over $R$ and in $G$), the $k'$-limit of $g(I)$ exists over $R$ in $\bar S$, with respect to $G.D$.}

\vs
\noindent
\textbf{(5.6)}
\textit{The necessary and sufficient condition that the norm-limit  exists over $R$ in $\bar S$, with respect to $G$,
are
that the $\sigma$-limit  exists over $R$ and with respect to $G$, and $c(y;G)=0$ in $R^o$.} 
\noindent
(Getchell [11] 415 for $G=G_1$.) From (4.21).

There are similar results in two dimensions.

\chapter{Two-dimensional Integration.}\label{Chapter 2}

\section{Introduction.}\label{Chapter 2 Introduction}
The definition of Burkill [6], of a function of intervals $I$, given in Chapter 1, \textsection 1, is equally the definition of a function of two-dimensional intervals (i.e.~of rectangles) $T$ whose sides are parallel to the coordinate axes, for we can substitute $T$ for $I$, $U$ for $S$, and ``rectangle'' for ``interval''. 

We denote by $\bar U$ the set of all $T^\sigma$ formed by finite (not necessarily overlapping) sums of points, lines parallel to the axes, and the interiors of rectangles of $U$. Thus if $V_1,V_2$ are in $\bar U$ so are $V_1^o, V_1^1, V_1+V_2$, and $V_1$ is a $T^\sigma$. The rectangles of $U$ are supposed to be contained in $W$.

The following is a selection of functions of rectangles $T=[a,b;c,d]$ for some bracket convention. Denote such $T$ collectively by $[a$---$b;c$---$d]$.
\begin{enumerate}
\item[(i)]
The \textit{diameter}, $\delta(T)$, the \textit{area} $m(T)$, and the \textit{parameter of regularity}, $p(T)$, of $T$, i.e.~if the shorter and longer sides are respectively of lengths $a$ and $b$, 
\[
p(T) =\frac ab;\;\;\;\;\;\;\;\;\;m(T) =ab;\;\;\;\;\;\;\;\;\;\delta(T) = \sqrt{a^2+b^2}.
\]
\item[(ii)]
The two-dimensional \textit{Stieltjes difference} $S(T)=S(f;T)$ of the real finite function $f(x,y)$, where
\[
S(T)=f(a,c)-f(a,d) -f(b,c) + f(b,d).
\]
\item[(iii)]
The difference $f(b,d) - f(a,c)$.
\item[(iv)]
In the theory of the total variation of $f(x,y)$ we study $|g(T)|$ where $g(T)$ has the form (ii) or (iii).
\item[(v)]
Let $(x_1,y_1), \ldots , (x_n,y_n), \ldots $ be a sequence of distinct points and $f(x,y)$ be such that $\sum_{n=1}^\infty f(x_n,y_n)$ is absolutely convergent. Put \[g(T) = (T)\sum f(x_n,y_n).\]
Then $g(T)$ is additive.
\item[(vi)]
Let $f(x,y)$ be real and continuous in $T$. We can then define
\begin{eqnarray*}
H_1(f;T) &=& \int_a^b \left(f(x,d) - f(x,c)\right) dx, \vt
H_2(f;T) &=& \int_c^d \left(f(b,y) - f(a,y)\right) dy.
\end{eqnarray*}
If $f(z) =u(x,y) +\iota v(x,y)$ where $z=x+\iota y$, and $u$ and $v$ are real and continuous in $T$, we can define
\[
H(f;T) = -\left(H_1(u;T) + H_2(v,T)\right)
+ \iota\left(H_2(u;T) - H_1(v,T)\right),
\]
the curvilinear integral of $f(z)$ following the contour of $T$. (Saks [17] 195.)
\item[(vii)]
If $f(x,y)$ is real in $T$ put $g_1(I_y) =f(x,y_2) - f(x,y_1)$ for $I_y=y_1$---$y_2$, and $g_2(I_x) =f(x_2,y)-f(x_1,y)$ for $I_x=x_1$---$x_2$, and then
\begin{eqnarray*}
V_x[f;T]&=& (F_1)\int_{(c,d)}|g_1(I_y)|, \vt
V_y[f;T]&=& (F_1)\int_{(a,b)}|g_2(I_x)|, \vt
V[f;T]&=& \int_a^b V_x[f;T]\,dx + \int_c^d V_y[f;T]\, dy.
\end{eqnarray*}
$V[f;T$ is the \textit{Tonelli variation} of $f(x,y)$ in $T$. (Tonelli [9],  Saks [177 169.)
\item[(viii)]
For $f(x,y)$ real and continuous in $T$ put
\begin{eqnarray*}
G_1(f;T) &=& \int_a^b |f(x,d) - f(x,c)|\,dx, \vt
G_2(f;T) &=& \int_c^d |f(b,y) - f(a,y)|\,dxy, \vt
G^2(f;T) &=& G_1^2(f;T) + G_2^2(f;T) + (b-a)^2(d-c)^2.
\end{eqnarray*}
These are the expressions of Z.~de Geocze [1]. Then the integral of $G(f;T)$ ($>0$) over a rectangle $T_0$ is equal to the area of the surface $z=f(x,Y)$ on $T_0$. (Saks [17] 171--180, particularly 179.)
\item[(ix)]
Let $x=x(u,v)$ and $y=y(u,v)$ be real and continuous for $(u,v)$ in $T_0$. Then as $(u,v)$ traces out $FT \subset T_0^o$, $(x,y)$ traces out a curve $C(T)$ in the $(x,y)$ plane. Putting
\[
x(a,c) = x_1,\;\;\;\;\;x(b,c) =x_2,\;\;\;\;\;x(b,d) = x_3,\;\;\;\;\;x(a,d) =x_4,
\]
and similarly for $y$, we consider
\[
\Delta(x,y) =\frac 12 \left(
x_1y_2-x_2y_1+x_2y_3-x_3y_2+x_3y_4-x_4y_3+x_4y_1-x_1y_4\right),
\]
in finding the area enclosed by $C(T)$. (Burkill [7] 313.) Similarly for a curved surface
\[
x=x(u,v),\;\;\;\;\;y=y(u,v),\;\;\;\;\;z=z(u,v),
\]
using $\;\;\Delta(y,z),\;\; \Delta(z,x),\;\; \Delta(x,y)$. (Burkill [7] 314.)
\end{enumerate}

\section{Burkill integration.}\label{Chapter 2 Burkill integration}  
A \textit{Riemann succession} of divisions of $V$, a $T^\sigma$, is a sequence $\{D_n\}$ of divisions of $V$ in which norm$(D_n) \rightarrow 0$. If $V$ is a $T$, a \textit{restricted Riemann succession} of divisions of $V$ is a Riemann succession of restricted divisions of $V$.

To give some distinction between the one- and the two-dimensional cases we will use the following symbols.
\begin{itemize}
\item
$T$ for a rectangle replaces $I$ for an interval.
\item
$V$ for a $T^\sigma$ replaces $R$ for an $I^\sigma$.
\item
The family $H$  replaces the family $F$.
\item
The family $L$ of divisions replaces the family $G$ of divisions.
\end{itemize}
We keep unchanged the symbols $D$ (division), $E$ (general set), $W$ (the fundamental one- or two-dimensional interval).

With this replacement we replace axioms (1) to (viii) (without (vi)), respectively by (i)' to (viii)' (without (vi)'). Axiom (vi) is exclusively an axiom for one-dimension, and so must be modified, as will be seen.

We see at once that a family $H$ of restricted Riemann successions of divisions (denoted by $H_r$) will not obey (i)' since restricted Riemann successions are defined only over rectangles. We might extend the definition of a restricted division to $V$, a $T^\sigma$, by prolonging the sides of the component rectangles of $V'$ indefinitely both ways and using these in a mesh of lines extended indefinitely both ways, and parallel to the axes. 

By taking these segments of the lines of the mesh which lie in $V'$, and taking suitable bracket conventions along these lines we obtain a division of $V$. With this definition we can conside an $H$ of restricted Riemann successions which satisfies (i)' and (ii)'. But obviously (iii)' will not be satisfied. It is possible for $H_r$ to satisfy (iv)', (v)', (vii)', (viii)'.

Kempisty [12] and [18] supposes that the rectangles $T$ of $U$ have parameter of regularity $(p(T) \geq 1/2$; or for $0<\ve<1$, $p(T) \geq \ve$.

If $U$ consists of all rectangles $T$ in $W$, denote $U$ by $U_0$. If $U$ consists of all $T$ in $W$ with $p(T) geq \ve$ denote $U$ by $U_\ve$. Since we suppose that $0<\ve<1$, $U_1$ is not defined. As pointed out by Kempisty [12] 13, if we took $p(T)=1$, we would only be able to define a division over a $T$ with commensurable sides, or a finite sum of such $T$.

If $U=U_s$ and $H$ consists of all Riemann successions of divisions using rectangles of $U_s$, and over all $T^\sigma$ in $W$, denote $H$ by $H_s$ and the corresponding $L$ by $L_s$. If $U=U_s$ and $H$ consists of all Riemann successions of restricted divisions using rectangles of $U_s$, and over all $T^\sigma$ in $W$, denote $H$ by $H_{rs}$ and the corresponding $L$ by $L_{rs}$. This is for $0\leq s<1$. If $H=H_r$ denote the corresponding $L$ by $L_r$.

For a general $H$ results corresponding to those in Chapter 1, \textsection 1, now follow. For $H_r$ results corresponding to (2.02), ..., (2.10), (2.14), ... , (2.20) hold.

The Burkill integral with respect to the family $H_{r0}$ (or the equivalent norm-limit with respect to $L_{r0}$) is defined by Burkill [6] \textsection 2 and taken as standard, whilst the Burkill integral with respect to $H_0$ (or the equivalent norm-limit with respect to $L_{r0}$) is defined and called the ``extended integral''. But since the family $H_{r0}$ cannot satisfy (iii)', it would seem best to take the extended integral as standard, and to call the other integral the ``restricted integral''. However, our general theory covers both integrals.

\vs
\noindent
\textbf{(2.1)}
\textit{Let $H_r$ be the family of all restricted Riemann successions in $H$. Then $(H_r)\overline{\int}_V\,g(T) \leq (H)\overline{\int}_V\,g(T)$, and there are an $H$ and a $g(T)$ so that
\[(H_r)\overline{\int}_V\,g(T) < (H)\overline{\int}_V\,g(T).\] Similarly for lower integrals and for upper and lower norm-limits.
}

\noindent
The first result is obvious. For the second, take, for example, $H=H_0$ and the following $g(T)$.

Let $M$ be the unit square with two vertices $(0,0)$ and $(1,1)$. The points $P(1/3,1/2)$ and $Q(2/3,1/2)$ lie in $M$. For squares $T$, centre $P$ and side $2^{-2n}$, and squares $T$, centre $Q$ and side $2^{-2n-1}$, put $g(T)=1$ ($n=1,2,\ldots $). Otherwise put $g(T)=0$. Then $(H_0)\overline{\int}_M\,g(T)=2$, but $(H_{r0})\overline{\int}_M\,g(T)=1$ since in a restricted division we can have at most one of the special squares.

By using axioms in Chapter 1 and working from them, we have been enabled to generalise the results of \textsection 2 of that Chapter to the case of two dimensions. We could equally well have generalised the results to the case of any finite number of dimensions, or even of an abstract space. However, points are introduced in axiom (vi), so that care is necessary to generalise the results of Chapter 1, \textsection\textsection 3,4 to the case of $n$ ($>1$). The case $n=2$ will be considered in \textsection\textsection 3,4  of the present Chapter.

There now follows a theorem of Fubini type, in which we disregard bracket conventions.

\vs
\noindent
\textbf{(2.2)}
\textit{Let $T_1=[\alpha$---$\beta; \gamma$---$\delta]$ be a rectangle in $\bar U$ and put $I_{1,x} =\alpha$---$\beta$, $I_{1,y}= \gamma$---$\delta$. For  $T=[a$---$b; c$---$d]$ 
in $T_1$
put
$I_{x} =a$---$b$, $I_{y}= c$---$d$, and $g(T) = g(I_x,I_y)$.
Let $F_x$ be the family obtained from an $H_r$ by replacing each $T$ by $I_x$, and similarly for $F_y$. Then let $H$ be the family obtained by using any $\{D_n\}$ from $F_x$ to give lines parallel to the $y$-axis, and then in each $D_n$ using any division $D$ (in $G_y$) in each ``column'' such that norm$(D_n)<\ve_n$, and $\ve_n \rightarrow 0$ as $n \rightarrow \infty$. Then}
\begin{eqnarray*}
(H)\underline{\int}_{T_1} \, g(T) & \leq &
(F_x)\underline{\int}_{I_{1,x}} \,\,(N;G_y)\underline{\int}_{I_{1,y}}\;g(I_x,I_y) \vt
&\leq &
(F_x)\overline{\int}_{I_{1,x}} \,\,(N;G_y)\overline{\int}_{I_{1,y}}\;g(I_x,I_y) \vt
&\leq &\overline{\int}_{T_1}\,g(T).
\end{eqnarray*}
For given $\ve>0$ there is a $\{D_n\}$ in $F_x$ such that
$(F_x)\overline{\int}_{I_{1,x}} \,(N;G_y)\overline{\int}_{I_{1,y}}\;g(I_x,I_y) <$
\begin{eqnarray*}
&<&
\lim_{n\rightarrow \infty} (I_{1,x};D_n) \sum \;(N;G_y)\overline{\int}_{I_{1,y}}\;g(I_x,I_y) +  \ve \vt
&<&
\overline{\lim}_{n\rightarrow \infty} (I_{1,x};D_n)\sum \; (I_{1,y};D_n) \sum g(I_x,I_y) +  2\ve
,\vt
&&\;\;\;\;\;\;\mbox{ where $D=D(n;I_x)$ and norm}(D)<\frac 1n,  \vt
&<&
\overline{\lim}_{n\rightarrow \infty} (T_1;D_n'')\sum g(T) + 2\ve, \vt
&&\mbox{where $D_n''$ is the division of $T_1$ consisting of $D_n$ along the $x$-axis}\vt
&&\mbox{ and $D=D(n;I_x)$ along the ``columns''}\vt
&\leq & (H) \overline{\int}_{T_1}\,g(T) + 2\ve.
\end{eqnarray*}
Since $\ve>0$ is arbitrary, we have one result. Similarly for the lower integrals.

\vs
\noindent
\textbf{Corollary.}
\textit{If $g(T)$ is Burkill-integrable over $T_1$ with respect to $H$ then}
\[
(H)\int_{T_1}g(T)=
(F_x){\int}_{I_{1,x}} \,(N;G_y)\overline{\int}_{I_{1,y}}\;g(I_x,I_y) 
= (F_x){\int}_{I_{1,x}} \,(N;G_y)\underline{\int}_{I_{1,y}}\;g(I_x,I_y) .
\]

\section{Burkill integration and lines of division.} 
\label{Burkill integration and lines of division}
Let each of the (closed) non-overlapping lines $l_1, \ldots , l_m$ be parallel to one of the axes, with the property that for $i=1, \ldots ,m-1$, $l_i$ and $l_{i+1}$ have only a common end-point, and $l_i,l_j$ have no common point if $|i-j|>1$. Then the set $l_1, \ldots ,l_m$ or $\{l\}$, with \textbf{end-points} that end-point of $l_1$ and that of $l_m$ which are not end-points of $l_2, l_{m-1}$ respectively. Then $l_1, \ldots , l_m$ cannot  cross each other.

If $V_1$ and $V_2$ are non-overlapping $T^\sigma$ then $V_1'$ and $V_2'$ meet in an at most finite number of broken lines. Let $(V_1')^o$ and $(V_2')^o$ be connected and let $F(V_1').F(V_2')$ be split up into non-overlapping non-crossing broken lines $\{l\}_1, \ldots , \{l\}_m$ such that each $\{l\}_i$ lies inside $V=V_1'+V_2'$, \textbf{except that the end-points of $\{l\}_i$ lie on $FV$.} (This may not always be possible.)

Suppose that (say) $\{l\}_1$ and $\{l\}_2$ have a common point $x$. If $x$ is not an end-point of $\{l\}_1$ it is inside $V$, and so not an end-point of $\{l\}_2$. Then since $\{l\}_1$ and $\{l\}_2$ are non-overlapping and non-crossing, the two lines of of $\{l\}_1$ which meet at $x$ are perpendicular, and similarly for $\{l\}_2$. We thus have the arrangement of Fig.~1, page \pageref{Figures}. 

Since $x$ is not in $V_1')^o + (V_2')^o$, connection of $(V_1')^o$ through $x$ is excluded, and similarly for $(V_2')^o$. Then by topological considerations, if $(V_1')^o$ is connected, $V_2')^o$ cannot be, and vice versa.

Thus a common point $x$ of $\{l\}_1$ and $\{l\}_2$ must be a common end-point, and so must lie in $FV$. Then the two lines, say $l_1,l_2$ of $\{l\}_1, \{l\}_2$, which have $x$ as a common end-point, cannot be in line. For $l_1, l_2$ do not lie along $FV$, so that since $V$ is a $T^\sigma$, the two sides, $s_1,s_2$ say, of $V$ meeting at $x$ must each be perpendicular to $l_1$ and $l_2$. Hence if $l_1,l_2$ are in line we have the arrangement of Fig.~2, page \pageref{Figures}.

\begin{figure}[!t]
\begin{center}
\includegraphics[width=6in]{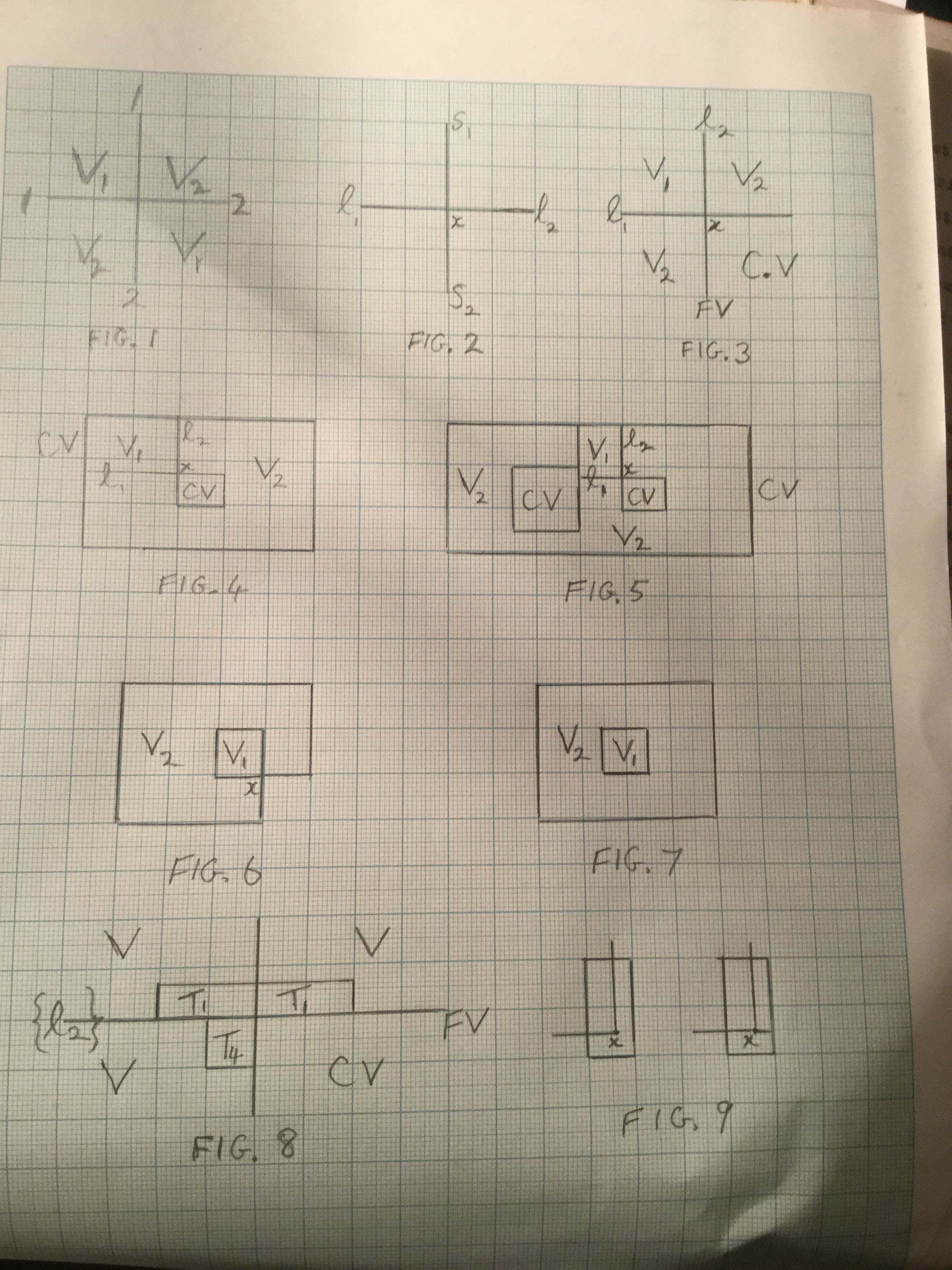}
\caption{Figures}
\label{Figures}
\end{center}
\end{figure}

But apart from end-points, $l_1$ and $l_2$ lie in $V^o$, so that this arrangement is impossible since each side of $V$ is parallel to one of the two perpendicular axes. Hence $l_1,l_2$ are perpendicular, with the arrangement of Fig.~3, page \pageref{Figures},  or the similar arrangement with $V_1$ and $V_2$ interchanged.
 
Since $(V_2')^o$ is connected there is a polygonal line through $x$ and in $V$ which surrounds part of $CV$ and cannot be shrunk to a point. This arrangement is like Fig.~4, Fig.~5 of page \pageref{Figures}.  Thus $V$ is multiply-connected.

When $(V_1')^o$ and $(V_2')^o$ are each connected and are non-overlapping $T^\sigma$ then either $F(V_1').F(V_2')$ can be split up into $\{l\}_1, \ldots , \{l\}_m$ as above, or else it contains a closed circuit. For example, Fig.~6, Fig.~7 of page \pageref{Figures}. 

If the non-overlapping, non-crossing broken lines $\{l\}_1, \{l\}_2$ have only two common points, the two end-points, then the whole forms a closed circuit and is a \textit{closed broken line} $\{cl\}$.

If the non-overlapping, non-crossing broken lines $\{l\}_1,\ldots , \{l\}_n$ are such that for $i=1, \ldots n-1$, $\{l\}_{i}$ and $\{l\}_{i+1}$ have only one common point, and end-point $x_i$, and $\{l\}_{i}$, $\{l\}_{j}$ are disjoint for $|i-j|>1$, and $\{l\}_{i}$, $\{l\}_{i+1}$ are perpendicular ``near'' $x_i$, then $x_1, \ldots , x_{n-1}$ are \textit{special corner points} and $\{l\}_{1}, \ldots , \{l\}_{n}$ form a \textit{multiple broken line} $\{ml\}$ with \textit{end-points} that of $\{l\}_{1}$ which is not $x_1$, and that of $\{l\}_{n}$ which is not $x_{n-1}$.
Note that
$\{ml\}$ is itself a broken line in which the points $x_1, \ldots ,x_{n-1}$ are special.

\vs
\noindent
\textbf{(3.01)}
\textit{If $V_1,V_2$ are non-overlapping $T^\sigma$ then $F(V_1').F(V_2')$ consists of an at most finite number of $\{l\}$, $\{cl\}$, $\{ml\}$. If $(V_1')^o$ and $(V_2')^o$ are each connected then either (a) $F(V_1').F(V_2')$ is a $\{cl\}$, or (b) there are disjoint} 
\[\{l\}_{1},\dots ,\{l\}_{p},\{ml\}_{1},\ldots ,\{ml\}_{q}\] \textit{such that}
\begin{enumerate}
\item[($\alpha$)]
$\sum_{i=1}^p\{l\}_{i} + \sum_{j=1}^q \{ml\}_{j}
= F(V_1').F(V_2')$,
\item[($\beta$)]
\textit{Each of $\{l\}_{i},\{ml\}_{j}$ lies inside $V=V_1'+V_2'$, except that the end-points of $\{l\}_{i}$ and of $\{ml\}_{j}$, and the special corner points of $\{ml\}_{j}$ all lie on $FV$.
\item[($\gamma$)]
Through each special corner point there is a polygonal line in $V$ which cannot be shrunk to a point.}
\end{enumerate}
If $V$ is a $T^\sigma$ then $(V')^o$ is a finite sum of disjoint open connected sets (i.e.~domains) $T^\sigma$.
A finite non-overlapping set $T_1, \ldots , T_n$ of rectangles in $U$ is a \textit{cluster} $\{T\}_i$ round an $\{l\}$ if
\begin{enumerate}
\item[(a)]
for $i=1, \ldots ,n$, the set $\{l\}.T_i'$ is not null,
\item[(b)]
$\{l\}$ is inside $\sum_{i=1}^n T_i'$ except that the end-points of $\{l\}$ are on $F\left(\sum_{i=1}^n T_i'\right)$.
\end{enumerate}
Similarly for an $\{ml\}$.

$\{l\}_j$ and $\{l\}_{j+1}$ define a quadrant $Q_j$ ``near'' $x_j$, and the quadrant at $x_j$ opposite to $Q_j$ is free of rectangles $T_i$ ``near'' $x_j$.

Similarly for a $\{c\}$ (which has no end-points, but sometimes special corner points. See Fig.~6, page \pageref{Figures}. 

Supposing the non-empty $H,U$ to exist, we may deduce from (i)' the following results.

\vs
\noindent
\textbf{(3.02)}
\textit{If $V$ is in $\bar U$ then $V$ contains rectangles $J$ of arbitrarily small norm, which are of the following types:}
\begin{enumerate}
\item[\textbf{(a)}]
\textit{Along each side of $V'$ there is a set of non-overlapping rectangles $J$ such that each point of the side is in at least one of the $J'$.}
\item[\textbf{(b)}]
\textit{If $\{l\}$ is a broken line (closed broken line, multiple broken line) in $V^o$ except that its end-points and special corner points (if any) are on $F(V')$, there is a set $J$ forming a cluster point round $\{l\}$.}
\end{enumerate}
If a broken line (closed broken line, multiple broken line)
$\{l\}$ forms part of the lines of division of all but a finite number of the divisions in a Riemann succession, then $\{l\}$ is a \textit{permanent broken line (closed broken line, multiple broken line) of division} of the Riemann succession.

Let $V$ be in $\bar U$, let $(V')^o$ be connected, and let the disjoint \[\{ml\}_1, \ldots ,\{ml\}_q, \{l\}_1, \ldots \{l\}_p\]
\textit{have the property (1) with respect to $\{F_n\}$} if we have the following.

\begin{enumerate}
\item[\textbf{(a)}]
There are non-overlapping $V_1,V_2$ in $\bar U$ such that $(V_1')^o$ and $V_2')^o$ are connected, and $V_1'+V_2'=V'$, and
\[
F(V_1')F(V_2') = \sum_{i=1}^p \{l\}_i +  \sum_{j=1}^q \{ml\}_j.
\]
\item[\textbf{(b)}]
$\{D_{i,n}\}$ over $V_i$ ($i=1,2$) are in $H$, where the 
$\{D_{i,n}\}$ are defined in \textbf{(c)}.
\item[\textbf{(c)}]
If a rectangle of $D_n$ lies in $V_1$ or $V_2$ then it occurs with the same bracket convention in $D_{1,n} + D_{2,n}$. If a rectangle $T$ of $D_n$ does not lie in $V_1$ nor in $V_2$ then there is a finite number on non-overlapping rectangles $T_i$ such that $\sum T_i' =T'$, and each $T_i$ lies in $V_1$ or in $V_2$, occurs in $D_{1,n}+D_{2,n}$, and has a segment of $F(V_1')F(V_2')$ on its frontier.
\end{enumerate}
Similarly for a closed broken line.

\vs
\noindent
\textbf{(vi)'}
Every set of disjoint $\{l\}_1, \ldots \{l\}_p,\{ml\}_1, \ldots ,\{ml\}_q$, and every closed broken line, all satisfying (a) of property (1), have the property (1) with respect to every $\{D_n\}$ over $V$ and in $H$.

If $\{l\}_1, \ldots \{l\}_p,\{ml\}_1, \ldots ,\{ml\}_q$ have the property (1) for $\{D_n\}$ over $V$ and in $H$, we define $B(\{l\}_1, \ldots \{l\}_p,\{ml\}_1, \ldots ,\{ml\}_q;\{D\}_n) \;=$
\[
=\; \overline{\lim}_{n\rightarrow \infty} \left|
(V;D_n)\sum g(T) -(V;D_{1,n}+D_{2,n})\sum g(T)\right|.
\]
If $\{l\}_1, \ldots \{l\}_p,\{ml\}_1, \ldots ,\{ml\}_q$ have  the property (1) for all $\{D_n\}$ over $V$ and in $H$, we define
\[
b\left(\{l\}_1, \ldots  ,\{ml\}_q ;V;H\right)
= \mbox{l.u.b.}
B\left(\{l\}_1, \ldots ,\{ml\}_q;\{D_n\}\right)
\]
for all such $\{D_n\}$. There are similar definitions of $B(\{cl\};\{D_n\})$, $b(\{cl\};V;H)$.

\vs
\noindent
\textbf{(3.03)}
\textit{If $H$ satisfies (vi)' then there are non-overlapping $V_{1,1},\ldots ,V_{1,p}, V_{2,1}, \ldots , V_{2,p}$ in $\bar U$ such that  $\sum_{i=1}^p V_{1,i}' + \sum_{i=1}^q V_{2,i}' =V'$ and}
\[
b\left(\{l\}_1, \ldots \{l\}_p,\{ml\}_1, \ldots ,\{ml\}_q;V;H\right) =
\sum_{i=1}^p b\left(\{l\}_i;V_{1,i};H\right) +
\sum_{i=1}^q b\left(\{ml\}_i;V_{2,i};H\right).
\]
For $\{l\}_1, \ldots  ,\{ml\}_q$ are disjoint and closed and so at a positive distance apart, so that by (i)' there is in $L$ a division $D$ over $V$ with the following properties. Some of its division lines form broken lines separating 
$\{l\}_1, \ldots ,\{ml\}_q$ from each other and divide $V'$ up into $V_{1,1}, \ldots ,V_{2,q}$ in $\bar U$ such that $V_{1,i}^o$ contains $\{I\}_i$ ($i=1, \ldots ,p$) and $V_{2,i}^o$ contains $\{ml\}_i$ ($i=1, \ldots ,q$) (apart from end-points and special corner points).

Then by (vi)' applied for the broken lines from $D$, there is for each $\{D_n\}$ over $V$ and in $H$, a $\{D_{i,n}'\}$ over $V_{1,i}$ and in $H$, and a $\{D_{j,n}''\}$ over $V_{2,j}$ and in $H$, such that the clusters round $\{l\}_i, \{ml\}_j$ from $D_n$ when $n>n_0$, are unaltered when we use $\{D_{i,n}'\}$, $\{D_{j,n}''\}$. (Note that $n_0$ depends on norm$(D_{n_0})$. Hence
\[
\sum_{i=1}^p B\left(\{l\}_i; \{D_{i,n}'\}\right)
+ \sum_{j=1}^q B\left(\{ml\}_i; \{D_{j,n}''\}\right)
=B\left(\{l\}_1, \ldots ,\{ml\}_q;\{D\}_n\right).
\]
Further, if any $\{D_{i,n}'\}$ over $V_{1,i}$ ($i=1, \ldots ,p$) and any $\{D_{i,n}''\}$ over $V_{2,i}$ ($i=1, \ldots ,q$) are all in $H$, then by (iii)', $\{D_n'\}$ is in $H$, where $D_n'=\sum_{i=1}^p D_{i,n}'+ \sum_{j=1}^q D_{j,n}''$. Hence
\[
B\left(\{l\}_1, \ldots ,\{ml\}_q;\{D'_n\}\right)
=
\sum_{i=1}^p B\left( \{l\}_i ;\{D_{i,n}'\}\right)+ \sum_{j=1}^q B\left(\{ml\}_j ;\{D_{j,n}''\}\right).
\]
Hence the result.

Obviously all the $B$'s and $b$'s are non-negative. If $b(\{l\};V;H)>0$ we say that $\{l\}$ \textit{is a singularity in $V$ of $g(T)$ with respect to additivity and with respect to $H$}. Similarly for $\{cl\}$ and $\{ml\}$.

For the norm-limit method we obtain a property (L) with respect to $D$ by replacing $H$ by $L$; $\{D_n\}$ over $V$ by $D$ over $V$; $\{D_{i,n}\}$ by $D_i$ ($i=1,2$); and putting
\[
C\left(\{l\}_1, \ldots ,\{ml\}_q;V;L;\ve\right)
=\mbox{l.u.b.} \left|(V;D) \sum g(T) -(V;D_1+D_2)\sum g(T)\right|
\]
for every $D$ over $V$ and in $L$, with norm less than $\ve>0$, we define
\[
c\left(\{l\}_1, \ldots ,\{ml\}_q;V;L\right)
=
\lim_{\ve \rightarrow 0}C\left(\{l\}_1, \ldots ,\{ml\}_q;V;L;\ve\right).
\]
Similarly for $\{cl\}$. We now have a theorem similar to (3.03).

\vs
\noindent
\textbf{(3.04)}
\textit{If $H$ satisfies (vi)' then there are non-overlapping $V_{1,1},\ldots ,V_{1,p}, V_{2,1}, \ldots , V_{2,p}$ in $\bar U$ such that  $\sum_{i=1}^p V_{1,i}' + \sum_{j=1}^q V_{2,j}' =V'$ and}
\[
c\left(\{l\}_1, \ldots  \ldots ,\{ml\}_q;V;L\right) =
\sum_{i=1}^p c\left(\{l\}_i;V_{1,i};L\right) +
\sum_{j=1}^q c\left(\{ml\}_j;V_{2,j};L\right).
\]
If $c(\{l\};V;L)>0$ then $\{l\}$ \textit{is a singularity in $V$ of $g(T)$ with respect to additivity and with respect to $L$.}

Similarly for $\{cl\}$ and $\{ml\}$.

\vs
\noindent
\textbf{(3.05)}

\noindent
\textbf{(b)}
\textit{If $\{l\}$ has the property (l) with respect to all $\{D_n\}$ over $V$ and in $H$ then}
\[
B(\{l\};\{D_n\}) \leq c(\{l\};V;L),\;\;\;\;\;\;
b(\{l\};V;H) \leq c(\{l\};V;L).
\]
\textbf{(b)}
\textit{If $H$ also satisfies (iv)' then $b(\{l\};V;H) = c(\{l\};V;L) $.}

\noindent
Similarly for $\{ml\}$ and $\{cl\}$. Compare Chapter 1, (3.04a and b), and the proofs.
We now put \begin{eqnarray*}
\mbox{osc}(g;V;H)&=&(H)\overline{\int}_V\, g(T) -(H)\underline{\int}_V\, g(T) \;\;\;\;\mbox{ and }\vt
\mbox{osc}N(g;V;L)&=&
(N;L)\overline{\int}_V\, g(T) -(N;L)\underline{\int}_V\, g(T),
\end{eqnarray*}
whenever the right-hand sides exist, and have the following analogues of (3.04) of Chapter 1.

\vs
\noindent
\textbf{(3.06)}
\[
\begin{array}{rlllll}
B(\{l\};\{D_n\})&\leq & \mbox{osc}(g;V;H),\;\;\;\;\;
& B(\{cl\};\{D_n\})&\leq & \mbox{osc}(g;V;H), \vt
B(\{ml\};\{D_n\})&\leq & \mbox{osc}(g;V;H),&&&\vt
b(\{l\};V;H)&\leq & \mbox{osc}(g;V;H),&b(\{cl\};V;H)&\leq & \mbox{osc}(g;V;H),\vt
b(\{ml\};V;H)&\leq & \mbox{osc}(g;V;H),&&&\vt
c(\{l\};V;L)&\leq & \mbox{osc}N(g;V;L),
& c(\{cl\};V;L)&\leq & \mbox{osc}N(g;V;L),\vt
c(\{ml\};V;L)&\leq & \mbox{osc}N(g;V;L).
\end{array}
\]
We now use (vii)' and prove the following.

\vs
\noindent
\textbf{(3.07)}

\noindent
\textbf{(a)}
\textit{If $U$ satisfies (vii)' and if $V_1, V_2$ are in $\bar U$, so are $V_1-V_2'$, $V_2-V_1'$, $V_1.V_2$}.

\noindent
\textbf{(b)}
\textit{Also, if $V_1 \subset V_2$, then \textit{osc}$(g;V_1;H) \leq$ \textit{osc}$(g;V_2;H)$
whenever both sides exist, and \textit{osc}$N(g;V_1;L) \leq$ \textit{osc}$N(g;V_2;L)$ whenever both sides exist.}

\noindent
\textbf{(c)}
\textit{Also} 
\[
\begin{array}{rlllll}
b(\{l\};V;H)&\equiv & b(\{l\};H),\;\;\;\;\;
& b(\{cl\};V;H)&\equiv & b(\{cl\};H), \vt
b(\{ml\};V;H)&\equiv & b(\{ml\};H),& c(\{l\};V;L)&\equiv & c(\{l\};L),\vt  c(\{ml\};V;L)&\equiv & c(\{ml\};L), &
c(\{cl\};V;L)&\equiv & b(\{cl\};L),
\end{array}
\]
\textit{are all independent of $V$.}

\noindent
For (a) and (b) see Chapter 1 (3.06) (a) and (b). For (c), let $V_1\subset V_2$ where $V_1'$ and $V_2'$ are each connected and in $\bar U$. And let $\{l\}$ be inside $V_1'$, except that the ends of $\{l\}$ are on $F(V_1')$ and on $F(V_2')$. Then $V_2'-V_1'$ consists either of one connected set $V_3$, in $\bar U$ by (vii)', or else of two (or more) connected sets $V_3,V_4$. In the latter case, $V_3+V_4$ is in $\bar U$ by (vii)', and since $V_3$ and $V_4$ are a positive distance apart, $V_3$ and $V_4$ are each in $\bar U$.

We deal with the case $V_2'-V_1'=V_3$. The other case is similar. Then $FV_3.F(V_1')$ consists of $\{l'\}_1, \ldots , \{ml'\}_q$ (or of $\{cl'\}$), which satisfy (a) of property (l). Hence from (vi)', if $\{D_n\}$ is in $H$ and over $V_2$, we can find a $\{D_n'\}$ in $H$ and over $V_1$
such that for $n>n_0$, the cluster round $\{l\}$ is the same from $D_n'$ as from $D_n$.

Again, if $\{D_n'\}$ in $H$ is over $V_1$ there is by (vii)' and (i)' a $\{D_n''\}$ in $H$ and over $V_3$,
 and so by (iii)' a $\{D_n\}$ in $H$ and over $V_2$, such that for $n=1,2, \ldots$, all the rectangles of $D_n'$ are rectangles of $D_n$.
 
 These results imply that $b(\{l\};V_1;H)=b(\{l\};V_2;H)$ when $V_2'-V_1'=V_3$. When \[V_2'-V_1'=\sum_{i=3}^{n+2} V_i,\] apply the previous result $n$ times. 
Now let $V_1,V_2$ be general $T^\sigma$ in $\bar U$, such that $(V_1')^o$ and $V_2')^o$ are each connected, with $\{l\}$ inside $V_1'$ and $V_2'$, except that the ends of $\{l\}$ are on $F(V_1')$. Then $V_1, V_2 \subset V_1+V_2$, so that
\[
b(\{l\};V_1;H) =b(\{l\};V_1+V_2;H) =b(\{l\};V_2;H).
\]
Similarly for $\{cl\}$, $\{ml\}$, and the $c$'s.

\vs
\vs
\noindent
\textbf{(3.08)}
\textit{Let $V_1,V_2$ in $\bar U$ be non-overlapping with sum $V_3$. Then}

\noindent
\textbf{(a)}
\begin{eqnarray*}
(H)\overline{\int}_{V_3} \,g(T) &\leq &
(H)\overline{\int}_{V_1} \,g(T) +
(H) \overline{\int}_{V_2} \,g(T) +
(F(V_1') F(V_2'))\sum b(\{bl\};H),\vt
(H)\underline{\int}_{V_3} \,g(T) &\geq &
(H)\underline{\int}_{V_1} \,g(T) +
(H) \underline{\int}_{V_2} \,g(T) -
(F(V_1') F(V_2'))\sum b(\{bl\};H),
\end{eqnarray*}
\textit{where $\{bl\}$ is a collective notation for $\{l\}, \{cl\}$ or $\{ml\}$. Similar results hold for norm-limits.} Compare the proof of (3.08) of Chapter 1.

\vs
\vs
\vs
\noindent
\textbf{Corollary.}

\noindent
\textbf{(a)}
\textit{If $b(\{bl\};H)=0$ everywhere then the upper and lower integrals of $g(T)$ are additive.}

\noindent
\textbf{(b)}
\textit{The $b(\{bl\};H)$ for $(H)\overline{\int}_{V_1} \,g(T) $ is not greater than the $b(\{bl\};H)$ for $g(T)$.}

\vs
\noindent
\textbf{(3.09)}
\textit{If $H$ satisfies (vi)', (vii)', and if osc$(g;V;H)$ exists and is finite, then: }

\noindent
\textbf{(a)}
\textit{There is an enumerable set $s_1, \ldots ,s_n, \ldots $ of distinct singularities in $(V')^o$ (except that their ends are on $F(V')$) such that if $\{l\}$ (or $\{cl\}$, or $\{ml\}$) is distinct from $\sum_{i=1}^\infty s_i$ and is effective for $V$, then $b(\{l\};H)=0$.} Also

\noindent
\textbf{(b)}
\[
\sum_{n=1}^\infty b(s_n;H) \leq \mbox{osc}(g;V;H).
\]
\textit{A similar result holds for norm-limits.
}

\noindent
Since $\{l\}, s_1, \ldots , s_n, \ldots $ are distinct we can use the proof of (3.07), Chapter 1, but modified to deal with two dimensions.

\vs
\noindent
\textbf{(3.10)}
\textit{Let $V_1, \ldots ,V_n$ be non-overlapping, in $\bar U$, and with $\sum_{i=1}^n V_i' =V'$. Let $s_i'= F(V_1'+\cdots + V_i')$. Then}
\[
\sum_{i=1}^{n-1} b(s_i';H) \leq 2\mbox{osc}(g;V;H).
\]
$F(V_1')F(V_2')$ is distinct from $F(V_i')$ for $i>2$ (except perhaps for a finite number of points) since if $l$ is a line of $F(V_1')F(V_2')$, $V_1'$ and $V_2'$ lie on opposite sides of $l$, and $l$ can only meet $F(V_i')$ ($i>2$) at an end of $l$.

Similarly $F(V_1'+\cdots +V_{i-1}')F(V_i')$ is distinct from $F(V_m')$ for $m>i$ (except perhaps for a finite number of points). Hence pairs of $s_1', \ldots , s_{n-1}'$ can only meet in a finite number of points.

Let $W_i$ in $\bar U$ be in $V'$ and with in $\delta$ ($>0$) of $s_i'$ ($i=1, \ldots ,n-1$) such that  $s_i'$ is inside $W_i$, except that the ends of the components $\{l\},\{ml\}$ of $s_i'$, and special corner points, are on $F(W_i')$.
We put 
\[
f_j = \left(\sum_{i=1}^j W_i'\right) .W_{j+1}'.
\]
Then $f_j$ is in $\bar U$ and can be as small as we please by taking $\delta$ small enough, since $f_j$ lies in the sums of neighbourhoods of a fixed number of points. Thus we can suppose $f_1, \ldots, f_{n-1}$ distinct.

We can also take $f_j$ so that
\[
b\left(Ff_j;H\right) < \frac \ve{2^{j+1}} \;\;\;\;\;\;(j=1, \ldots ,n-2).
\]
(For suppose $j$ fixed. We can take a sequence of distinct $Ff_j$, say $s_1, \ldots ,s_n, \ldots$, and as in (3.09) prove that $\sum_{n=1}^\infty b(s_n:H)$ is convergent and so $b(s_n;H) \rightarrow 0$.)
Then by (3.06),
\[
\sum_{i=1}^{n-1} b(s_i';H) \leq \sum_{i=1}^{n-1} \mbox{osc}(g;W_i;H).
\]
But by (3.08), for arbitrary $V_1, V_2$ in $\bar U$, $\mbox{  osc}(g;V_1;H)+\mbox{osc}(g;V_2;H)\;\;\leq$

\begin{eqnarray*}
&\leq & 2b(F(V_1'V_2');H) + \mbox{osc}(g;V_1;H) +\mbox{osc}(g;V_2-V_1';H)+\mbox{osc}(g;V_1.V_2;H), \vt
&\leq & \mbox{osc}(g;V_1+V_2;H) + \mbox{osc}(g;V_1.V_2;H) 
+2b(F(V_1'V_2');H).
\end{eqnarray*}
Hence
\begin{eqnarray*}
\sum_{i=1}^{n-1} \mbox{osc}(g;W_i;H) & \leq &
 \mbox{osc}(g;\sum_{i=1}^{n-1} W_i;H)
 + \sum_{i=1}^{n-1} \mbox{osc}(g;f_i;H) 
+ \sum_{i=1}^{n-1}2b(Ff_i;H) \vt
&<& \mbox{osc}(g;V;H) + \mbox{osc}(g;V;H) + \ve,
\end{eqnarray*}
since $f_1, \ldots ,f_{n-1}$ are disjoint; i.e.
\[
\sum_{i=1}^{n-1} b(s_i';H)  < 2\mbox{osc}(g;V;H) + \ve,
\]
giving the result.

\vs
\noindent
\textbf{Corollary 1.}
\textit{Let $\bar{b}(V;H)$ be the least upper bound of all such $\sum_{i=1}^{n-1} b(s_i';H)$. Then
\[
\bar{b}(V;H) \leq 2\mbox{osc}(g;V;H).
\]}
\noindent
\textbf{Corollary 2.}
\textit{From (3.08),}

\noindent
\textbf{(a)}
\[
(H)\overline{\int}_V\,g(T) \leq\sum_{i=1}^n (H)\overline{\int}_{V_i}\,g(T) + \bar{b}(V;H),
\]
\textbf{(b)}
\[
(H)\overline{\int}_V(H)\overline{\int}_T\,g(T) \leq(H)\overline{\int}_{V}\,g(T) 
\leq
(H)\overline{\int}_V(H)\overline{\int}_T\,g(T)+ \bar{b}(V;H).
\]
We now consider axiom (viii)' in conjunction with (iii)'.

\vs
\noindent
\textbf{(3.11)}
\textit{Let $L$ satisfy (viii)' and let the multiple broken line $\{ml\}$ be composed of broken lines $\{l\}_1, \ldots ,\{l\}_n$. Then
\[
c(\{ml\};L) \leq \sum_{i=1}^n c(\{l\}_i;L).
\]}  
Let $D$ be a division of $V$, giving $D'$ for $\{ml\}$ by (vi)'. At a special corner point $x$ we have the arrangement of Fig.~8, page \pageref{Figures}, where $T_1$ is from $D$, $FT_1$ includes $x$, and $T_1.\{l\}_i$ is not null. For norm$(D)<\ve$ and $\ve$ small enough, $T_1'$ does not meet $\{l\}_3, \ldots , \{l\}_n$, and $\{l\}_i.T_1'$ is straight  ($i=1,2$). Then by (vi)', either $T_1$ is split by $\{l\}_1$ into $T_2, T_3$ (where $T_2'.\{l\}_2$ is not null), or else $T_1$ has a side along $\{l\}_1$.

In the first case $T_2$ is part of the cluster round $\{l\}_2$ given by $D'$. Let the cluster be formed by $T_2$ and the set $S_2$ of $T$, and let $S_2$ have come from the set $S_1'$ of $T$ of $D$. We put $D_2'$ as the division given by the $T$ of $D'$, except that $S_2$ is replaced by $S_1'$. Let $S_3$ be the cluster around $\{l\}_1$ from $D$ and let the corresponding set of $T$ in $D'$ from $S_3$ be $S_4$. For $D_1'$ we replace $S_3$ in $D$ by $S_4$. Then by (viii)' and (iii)', $D_1, D_2$ are in $G$. And ``near'' $x$, the combined cluster round $\{l\}_1 + \{l\}_2$ is the same for $D_1'$ as for $D_2'$.

In the second case, if $T_1$ is split by $\{l\}_2$ we form $D_1', D_2'$ as above, interchanging $\{l\}_1$ and $\{l\}_2$. If $T_1$ is not split by either then by (vi)' a side of $T_1$ is along $\{l\}_i$ ($i=1,2$), and $T_1$ is in both $D$ and $D'$. Thus we put $D_1'=D=D_2'$.

In either case, $(V;D) \sum g(T)- (V;D') \sum g(T) \;=$
\[ \begin{array}{rllll}
=\;\;(V;D_2') \sum g(T)& - &(V;D') \sum g(T) & + &(V;D_1) \sum g(T)\;\; -\vt  -\;\;(V;D_1') \sum g(T)& +& (V;D_1') \sum g(T)& -& (V;D_2') \sum g(T).
\end{array}
\]
Repeating the construction for all corner points $x_1, \ldots , x_{n-1}$, we obtain $D_1'', \ldots ,D_n'''$, all in $G$, and
\begin{eqnarray*}
(V;D) \sum g(T)- (V;D') \sum g(T) &=&
\sum_{i=1}^n \left(
(V;D_i'') \sum g(T)- (V;D_i''') \sum g(T) \right), \vt
\left|(V;D) \sum g(T)- (V;D') \sum g(T)\right| &\leq &
\sum_{i=1}^n \left|
(V;D_i'') \sum g(T)- (V;D_i''') \sum g(T) \right| \vt
&\leq & \sum_{i=1}^n C\left(\{l\}_i;V;L;e\right). \vt
\mbox{i.e. }\;\; 
C\left(\{ml\};V;L;e\right) & \leq &
\sum_{i=1}^n C\left(\{l\}_i;V;L;e\right), \vt
c\left(\{ml\};L\right) & \leq &
\sum_{i=1}^n c\left(\{l\}_i;L\right).
\end{eqnarray*}  
In a \textit{modified cluster} round a line $l$ we allow the finite non-overlapping rectangles $T_1, \ldots , T_n$ of $U$ to satisfy (a) and (b) of the definition of a cluster, except that at each end $x$ of $l$ there may be a set of rectangles (say $T_1, \ldots , T_p$) which do not all satisfy (a) and (b), but such that $\sum_{i=1}^p T_i'=T'$, another rectangle, where $T$ is in $U$ and $x$ is inside $T$.

Let $l$ be part of a broken line $\{l\}$ satisfying (a) of property $l$ for $V$, so that $l$ also satisfies (a) for some $V_3$ (or let $l$ be part of $\{cl\}$ satisfying (a) of property (l) for $V$). Let $D$ be in $L$, and let $c_1$ be the modified cluster round $l$ formed from $D$. Then property (l) defines another modified cluster, say $c_2$, round $l$, \textit{such that if $E_i$ is the set of points    of $c_i$ ($i=1,2$) then $E_1' = E_2'$.} We put 
\[
C_m(l;V;L;\ve) = \mbox{l.u.b.} \left| (c_1) \sum g(T) -(c_2) \sum g(T)\right|
\]
for every $D$ over $V$ and in $L$, with norm less than $e>0$.
\[
c_m(l;V;L)= \lim_{e\rightarrow 0} 
C_m(l;V;L;e).
\]
By (vii)', as in (3.07c), we may omit $V$, provided that we take $l$ in $V^o$. As is easily seen,
\[
0\leq c(l;L) \leq c_m(l;L)
\]
since we can take a division-line perpendicular to each end of $l$, in which case our modified cluster can be an ordinary cluster.

\vs
\noindent
\textbf{(3.12)}
\textit{Let $\{l\}$ (or $\{cl\}$ consist of lines $l_1, \ldots ,l_n$. Then $c(\{l\};L)\leq \sum_{i=1}^n c_m(l_i;L)$.} See Fig.~9, page \pageref{Figures}. 

\vs
If in $D$ a rectangle $T$ covers a corner $x$ of $\{l\}$ as shown in Fig.~9, page \pageref{Figures},  we put $T$ and its decompositions $T_1, \ldots ,T_n$ for one line meeting at $x$, and for the other line we replace $T$ in $D$ by $T_1, \ldots ,T_n$, and also use $T_1,\ldots ,T_n$ corresponding to the $D_1+D_n$ of property (l). The modified $D, D_1, D_2$ all lie in $L$ by (viii)'. The result then follows.

As before, we now have

\vs
\noindent
\textbf{(3.13)}
\textit{Let $l$ satisfy (a) of property (l) for some $V_3 \subset V$, and let $l$ lie in $V^o$. Then $c_m(l;L) \leq \mbox{osc}N(g;V;L)$.}

\vs
\noindent
\textbf{(3.14)}
\textit{If osc$N(g;W;L)$ exists and is finite then there are an at most enumerable set $l_{1,1}, \ldots ,l_{1,n}, \ldots $ of lines parallel to the $x$-axis, and an at most enumerable set $l_{2,1}, \ldots , l_{2,n}, \ldots $ of lines parallel to the $y$-axis, such that}

\noindent
\textbf{(a)}
\textit{if $\{l\}$ ($\{ml\},\{cl\}$) does not include a segment of any one of 
\[
l_{1,1}, \ldots ,l_{1,n}, \ldots;l_{2,1}, \ldots , l_{2,n}, \ldots,
\] 
then  $c(\{l\};L)=0$,}

\noindent
\textbf{(b)}
\[
\sum_{n=1}^\infty c_m(l_{i,n};L) \leq \mbox{osc}N(g;V;L)\;\;\;\;\;(i=1,2),
\]
\textit{where the clusters all lie in $V$. (In particular, $V=W$.)}

\noindent
The proof follows from that of I(3.07) (i.e.~(3.07) of Chapter 1) and then uses (3.11) and (3.12). We now have analogues of I(3.09) to I(3.12).

\vs
\noindent
\textbf{(3.15)}
\textit{Let the upper and lower norm-limits of $g(T)$ over $W$, with respect to $L$, be finite. Then, given $\ve>0$ we can find $\delta>0$ such that every (finite, non-overlapping) set $S_2$ of rectangles of $U$ in $W$, forming a $T^\sigma$ set $V_2$, and with norm$(S_2)<\delta$, gives
\[
(S_2)\sum g(T) < (N;L) \overline{\int}_{V_2}\, g(T)
+(FV_2)\sum c_m(l_{1,n};L) +(FV_2)\sum c_m(l_{2,n};L) + \ve.
\]
Similarly for lower norm-limits.}

\vs
\noindent
\textbf{(3.16)}
\textit{Let the norm-limit of $g(T)$, over $W$, exist with respect to $L$. Then given $\ve>0$ we can find $\delta>0$ such that every (finite non-overlapping) set $S_2$ of rectangles of $U$ in $W$ with norm$(S_2)<\delta$ gives }
\[
\left|(S_2)\sum \left(g(T) - (N;L)\int_T\,g(T) \right)\right|< \ve,\;\;\;\;\;\;
(S_2)\sum \left|g(T) - (N;L)\int_T\,g(T) \right|< \ve.
\]
\noindent
\textbf{(3.17)}
\textit{Let $\delta$ be as in (3.15). If for a division $D$ in $L$ over $V$ in $\bar U$ we have
\[
(V;D)\sum g(T) >(N;L) \overline{\int}_V\,g(T) -\ve 
\]
and norm$(D)<\delta$, then every partial set $S_2$ of the set $S_3$ of rectangles of $D$, with the corresponding region $V_2$, is such that
\[
(S_2)\sum g(T) >-2\ve-(F(V_1-V_2))\sum \left(c_m(l_1;L)+c_m(l_2;L)\right) +(N;L) \overline{\int}_{V_2}\, g(T).
\]
}

\noindent
\textbf{(3.18)}
\textit{Let $x$ be a corner-point of the rectangle $T$ in $U$. Then as $\delta(T) \rightarrow 0$,}
\[
\overline{\lim} \left(g(T)-(N;L)\overline{\int}_T\,g(T)\right) =0=
\underline{\lim} \left(g(T)-(N;L)\underline{\int}_T\,g(T)\right) .
\]
If $x$ is on a side of $T$, or inside $T$, there are similar (but more complicated) results to the case $J=I$ in I(3.12).

\section{$k$-integration and the $\sigma$-limits.}
\label{Chapter 2 k-integration and the sigma-limit}
As in I, \textsection 4, we can easily see that when the Burkill integrals (or norm-limits) are finite and unequal they are not necessarily additive, and we are thus led to consider the use of permanent lines of division. To do this effectively there must be at most an enumerable number of singularities in $W$, and this is only ensured when we introduce axiom (viii)', the modified clusters, the $c_m(l;L)$, and the lines $l_{1,1}, \ldots, l_{1,n}, \ldots , l_{2,1}, \ldots ,l_{2,n}, \ldots $.

We therefore suppose that $L$ obeys axioms (i)' to (viii)', save possibly (iv)' and (v)', and that $(N;L) \overline{\int}_W\,g(T)$ and $(N;L) \underline{\int}_W\,g(T)$ are finite. In our notation we can drop the dash from $k'$ because there is only one set of singularities.

If in a Riemann succession in $H$, each singular line $l_{i,n}$ is a fixed line of the divisions, with a specified set of bracket conventions (the $k$-convention) along it, after a finite number of divisions depending on $i,n$, then we call the succession a $k$-succession.

The family of all $k$-successions in $H$ is denoted by $H_k$. When $H=H_s$ then $H_k$ is denoted by $H_{sk}$ ($0\leq s<1$, or ``$s=r$''). $H_k$ does not necessarily exist, but we will suppose that $H_k$ exists and satisfies axiom (i)'. Then we have
\[
-\infty <(H)\underline{\int}_V\,g(T) \leq (H_k)\underline{\int}_V\,g(T)\leq (H_k)\overline{\int}_V\,g(T)\leq (H)\overline{\int}_V\,g(T) <+\infty.
\]
Since axioms (ii)', (iii)', (vii)', (viii)' are true for $H$ they are true for $H_k$. Axiom (iv)' in general is never true for an $H_k$, and the corresponding ``closure'' axiom is (ix)'. Let $\{D_n^{(i)}\}$ ($i=1,2, \ldots $) be a set of $k$-successions in $H_k$, of divisions of $V$ in $\bar U$. If $D_1^{(1)},D_2^{(1)},D_1^{(2)},D_3^{(1)},D_2^{(2)},D_1^{(3)},D_1^{(1)},
D_4^{(1)},\ldots$ is a $k$-succession it is in $H_k$.

If (iv)' is true for $H$ then (ix)' is true for $H_k$. If (v)' is true for $H$ it is true for $H_k$. If the $k$-convention is \textbf{all} conventions then (vi)' is true for $H_k$ since it is true for $H$. Otherwise we need in general to postulate that $H_k$ satisfies (vi)'.

We now define the upper and lower $k$-limits. The \textit{upper $k$-limit $(k;L)\overline{\int}_V\, g(T)$ of $g(T)$ over $V$ in $\bar U$, with respect to the set $L$ of divisions} is the greatest lower bound of numbers $dk(e;m;m)$ for all $e>0$ and integers $m,n$ where $dk(\ve;m;n)$ is the least upper bound of $(V;D)\sum g(T)$ for all $D$ in $L$ such that $D$ is a division of $V$ with (at least) the first $m$ singularities $\{l_{1,i}\}$ and the first $n$ singularities $\{l_{2,i}\}$ included with the correct $k$-convention, and with norm$(D)$ less than $e$. As $e \rightarrow 0$ and $m,n \rightarrow \infty$, $dk(e;m;n)$ is monotone decreasing, so that we may take $m=n$, and the lower bound is also the limit. For $m=n$ the division is a $k_{me}$-\textit{division}. Similarly for a $k_{me}$-\textit{set} of a finite number of non-overlapping rectangles.

A similar definition holds for the \textit{lower $k$-limit} $(k;L)\underline{\int}_V\,g(T)$. If for $V$ the upper and lower $k$-limits are equal, we say that the \textit{$k$-limit of $g(T)$ exists over $V$ in $\bar U$, with respect to the family $L$ of divisions}, and we write the common value as $(k;L)\int_V\,g(T)$.
\[
-\infty <(N;L)\underline{\int}_V\,g(T) \leq (k;L)\underline{\int}_V\,g(T)\leq (k;L)\overline{\int}_V\,g(T)\leq (N;L)\overline{\int}_V\,g(T) <+\infty.
\]
We have similar results to I(2.05), $\ldots $, I(2.20), as for the 1-dimensional $k'$-limits, and the following analogue of I(4.03), I(4.04).

\vs
\noindent
\textbf{(4.1)}
\textit{The upper and lower Burkill integrals with respect to $H_k$, and the upper and lower $k$-limits, are each additive.
}

\vs
\noindent
\textbf{(4.2)}
\textit{If keeping the same $k$-convention at each $\{l_{i,j}\}$, we add more permanent lines, each of which satisfy (a) of property (l), we do not alter the $k$-limits.
}

\vs
\noindent
\textbf{(4.3)}
\[
\begin{array}{rllll}
&&H&&\vt
&\nearrow&&\searrow& \vt
H_k &&&&N\;\;\;(<+\infty) \vt
&\searrow&&\nearrow& \vt
&&k&&
\end{array}
\]
There is an analogue of I(4.11), and also results similar to I(4.12), $\ldots$ , I(4.15).

Replacing points by lines, we may easily extend the idea of a $\sigma$-limit to two dimensions and produce results analogous to I(4.16), $\ldots$ , I(4.21).

\vs
\noindent
\textbf{(4.4)}
\textit{We can have $(H_{rk})\overline{\int}_{T_1}\,g(T)
< (H_{k})\overline{\int}_{T_1}\,g(T)$ for $H=H_0$.
}

\noindent
Let $T=[0,1;0,1]$. If $T$ lies in $0\leq x<1/2$ with edges along $y=0$ and $y=2^{-2n}$ put $g(T) =2^{2n}mT$. If $T$ lies in $1/2 \leq x<1$ with edges along $y=0$ and $y=2^{-2n-1}$ put $g(T) = 2^{2n+1}mT$. Otherwise put $g(T)=0$. Then
\[
(H_k)\overline{\int}_{T_1}\,g(T)=1,\;\;\;\;\;\;
(H_{rk})\overline{\int}_{T_1}\,g(T)=\frac 12.
\]
Note that the only singularity is $(y=0,\;0\leq x \leq 1)$.

\chapter{Special Functions.}\label{Chapter 3}
\section{Functions of bounded variation.}
\label{Functions of bounded variation}
In 1-dimension the limit $(N;G) \overline{\int}_{R}\,|g(I)|$ is called the \textit{variation $\mbox{Var}(g;R;G)$ of $g(I)$ in $R$ with respect to the family $G$.} In 2-dimensions the limit $(N;L)\overline{\int}_{V}\,|g(T)|$ is called the \textit{variation $\mbox{Var}(g;V;L)$ of $g(T)$ in $V$ with respect to the family $L$.}

For $G=G_1$, $L=L_0$, these are implicit in Burkill [6]. Kempisty [12] 17, \textsection 3, gives them for $G=G_1$, $L=L_{\frac 12}$. In these cases there is no distinction between the Burkill integrals and the norm-limits.

When $g(I) = S(f;I)$ (\textit{the Stieltjes increment of $f$---P.M.}) then $\var(g;R;G_1)$ is the total variation of $f(x)$ in $R$. This explains the name of $\var(g;R;G)$.

When $\var(g;R;G)$ is finite we say that \textit{$g(I)$ is of bounded variation (b.v.) in $R$ with respect to the family $G$}. Similarly in two dimensions.

\vs
\noindent
\textbf{(1.01)}

\noindent
\textbf{(a)}
\textit{Let $\var(g;R;G) \leq M$ $(<\infty)$. Then if $F$ satisfies (i), (ii), (iii), (vii), (viii), given $\ve>0$ there is a $\delta>0$ such that for each non-overlapping finite or infinite set $\{I_i\}$ in $S$ and in $R$, with norm $<\delta$, we have $\sum |g(I_i)|<M+\ve$.}

\noindent
\textbf{(a)}
\textit{Let $F$ satisfy (i). If for every $\ve>0$ there is a $\delta>0$ such that for norm less than $\delta$ we have 
$\sum |g(I_i)|<M+\ve$, then $\var(g;R;G)\leq M$. Similarly in two dimensions.}

\noindent
(a) Let $\delta>0$ be such that if $D$ in $G$ and over $R$ has norm$(D)<\delta$, then
\[
(R;D)\sum |g(I)| < M + \frac \ve 2.
\]
Now let the non-overlapping $I_1, \ldots ,I_n$ in $S$ and in $R$ have norm less than $\delta$. Then $R_1 = \sum_{i=1}^n I_i'$ is in $\bar S$ so that by (vii), $R-R_1$ is in $\bar S$. Hence there is a $D_1$ in $G$ over $R-R_1$ such that norm$(D_1)<\delta$. Then by (viii) and (iii), $D_1$ with $I_1,\ldots ,I_n$ form a division $D$ in $G$ over $R$, so that
\[
\sum_{i=1}^n |g(I_i)| \leq (R;D)\sum |g(I)| < M + \frac \ve 2.
\]
This is true for $n=1,2, \ldots $, so that
$
\sum_{i=1}^\infty |g(I_i)| \leq  M + \frac \ve 2 <  M + \ve .
$ Hence (a).

\noindent
(b) If $D$ in $G$ over $R$ has norm$(D)<\delta$ (this is possible by (i)) then
\[
(R;D) \sum |g(I) <M+\ve.
\]
Hence $\var(g;R;G) \leq M + \ve$, and this is for all $\ve>0$. Hence (b).

Thus, in some sense the variation is a limiting variation.

\vs
\noindent
\textbf{(1.02)}
\textit{If $\bar S$ satisfies (vii) and if $R_1,R_2$ are in $\bar S$, and $R_1 \subset R_2$, then \[\var(g;R_1;G) \leq \var(g;R_2;G).\] }
For $R_2-R_1'$ is in $\bar S$, so that there is in $G$ a division $D$ over $R_2-R_1'$ with norm$(D)<\delta$. Then if $D_1$ is over $R_1$ and in $G$, $D+D_1$ is over $R_2$ and in $G$ by (iii), and
\[
(R_1;D_1)\sum |g(I)| \leq (R_2; D+D_1) \sum |g(I)|.
\]
Hence the result.

\vs
\noindent
\textbf{Corollary.}
\textit{If $g(I)$ is b.v.~in $R$ in $\bar S$, then $g(I)$ is b.v.~in every $R_1$, in $\bar S$, contained in $R$.}

We now put $A(g;R;G) \equiv A(R)$ as the maximum of
\[
\left|(N;G) \overline{\int}_R\,g(I)\right|\;\;\mbox{ and }\;\;
\left|(N;G) \underline{\int}_R\,g(I)\right|.
\]
\textbf{(1.03)} $A(R) \leq \var(g;R;G)$ (---because $|\sum g(I) \leq \sum |g(I)|$).

\vs
\noindent
\textbf{Corollary.}
\textit{If $g(I)$ is b.v.~in $R$, its upper and lower norm-limits are bounded in $R$.}

\vs
\noindent
\textbf{(1.04)} $\var(A;R;G) \leq \var(g;R;G)$.

\noindent
For by (1.03), and then by I(2.11),
$
(R)\sum A(I) \leq (R) \sum \var (g;I;G) \leq \var (g;R;G).
$

\noindent
\textbf{Corollary.}
\textit{If $g(I)$ is b.v.~in $R$, so are its upper and lower norm-limits.}

\vs
\noindent
\textbf{(1.05)}
\[\var (g;R;G) \leq \var (A;R;G)+ 2 (R^o)\sum c(y;G)\] \textit{whenever $A(W) < \infty$ (taking $W$ in $\bar S$) and $F$ satisfies (vi), (vii), (viii).}

\noindent
Since $A(W)<\infty$ and $R\subset W$, the series $(R^o)\sum c(y;G)$ exists by I(3.07).

In a division $D$ in $G$ of $R$, with norm$(D)<\delta$, take first the set of intervals in which $g \geq 0$, and then the rest. By I(3.09) and then I(2.09a),
\[
\begin{array}{rllll}
(R;D) \sum|g(I)|&\leq & (R;D)\sum A(I) &+& 2(R^o) \sum c(y;G) + \ve \vt
&\leq & \var(A;R;G) &+& 2(R^o) \sum c(y;G) + 2\ve .
\end{array}
\]
Hence the result.

\vs
\noindent
\textbf{Corollary 1.}
\textit{If $A(W)<\infty$ and $\var(g;R;G) = +\infty$ then $\var(A;R;G) = +\infty$. }

\vs
\noindent
\textbf{Corollary 2.}
\textit{If $g(I)$ is continuous then $\var(A;R;G)=\var(g;R;G)$.}

\vs
\noindent
\textbf{(1.06)}
\textit{If $F$ satisfies (vii), (viii) and if $A(R)\leq M$ for every $R$ in $\bar S$ then \[\var(g;R;G) \leq 4M.\]}
Given $\ve>0$ there is a $\delta = \delta(R;\ve)>0$ such that for every division $D$ in $G$ over $R$ with norm$(D)<\delta$, $|(R;D) \sum g(I)|<M+\ve$. Let $S_2$ be a finite set of non-overlapping intervals of $S$ in $R$ with norm less than $\delta$, and covering the set $R_2$ in $\bar S$, and put $R_1=R-R_2'$. By I(2.09b) we take $D_1$ in $G$ over $R_1$ with norm$(D_1) < \delta$ and
\[
(R_1;D_1) \sum g(I) > (N;G) \overline{\int}_{R_1}\, g(I) - \ve.
\]
Then $D=D_1 +S_2$ is over $R$, and in $G$ by (viii) and (iii), so that
\[
\begin{array}{rll}
(S_2)\sum g(I) &=& (R;D)\sum g(I) - (R_1;D_1)\sum g(I) \vt
&<& M+\ve -(N;G) \overline{\int}_{R_1}\,g(I) +\ve \vt
&\leq & 2(M+\ve).
\end{array}
\]
Similarly $(S_2)\sum g(I) >-2(M+\ve)$. Now let $S_2$ be the set of intervals in $(R;D)\sum g(I)$ for which $g\geq 0$, and let $S_3$ be the set of the rest. Then for $i=2,3$,
\[
(S_i)\sum | g(I)| <2(M+\ve),\;\;\;\;\;\;
(R;D)\sum |g(I)| <4(M+\ve).
\] Hence the result.
A similar result is given  by Dienes [14] when $g(I)$ is additive and ``bounded for all $I^\sigma$ in an interval $W$'', without using the norm-limits.

Note that it is not sufficient that $A(I) \leq M$ for every \textbf{interval} $I$ in $S$; for let $g(I) = S(f;I)$ where $f(x)$ is bounded but not of bounded variation.

\vs
\noindent
\textbf{(1.07)}
\textit{If $R_1,R_2$ in $\bar S$ are non-overlapping with sum $R_3$, then for $A(R_3)<\infty$,}
\[
A(R_3) \leq A(R_1)+A(R_2)+(FR_1.FR_2)\sum c(y;G).
\]
If
\[
A(R_3) = \left|(N;G) \overline{\int}_{R_3}\,g(I)\right| \neq \left| (N;G) \underline{\int}_{R_3}\,g(I) \right|
\]
then $(N;G) \overline{\int}_{R_3}\,g(I) >0$ and by I(3.08),
\[
\begin{array}{rll}
A(R_3) &= &-(N;G) \underline{\int}_{R_1}\,g(I)
-  (N;G) \underline{\int}_{R_2}\,g(I)
+ (FR_1.FR_2) \sum c(y;G)
 \vt
&\leq & A(R_1) + A(R_2) +(FR_1.FR_2) \sum c(y;G).
\end{array}
\]
We are left with the case
\[
A(R_3) = \left|(N;G) \overline{\int}_{R_3}\,g(I)\right| = \left| (N;G) \underline{\int}_{R_3}\,g(I) \right|
\]
If $(N;G) \overline{\int}_{R_3}\,g(I)>0$ we use the first argument. Otherwise $(N;G) {\int}_{R_3}\,g(I)$ exists and $c(y;G)=0$, and we can use the second argument.
\vspace{5pt}

\noindent
\textbf{Corollary.}
\textit{If $A(R)<\infty$ then $(N;G) \underline{\int}_{R}\,A(I) \geq A(R) - (R^o)\sum c(y;G)$.}

\vspace{10pt}

\noindent
\textbf{Unless otherwise stated, we shall suppose in the remainder of this section that $g(I)$ is b.v.}
We put 
\[
j(y) = \lim \var(g;(x,z);G)
\]
where $x<y<z$, $(x,z)$ is in $\bar S$, and $x,z \rightarrow y$. Then $j(y)$ is defined for every $y$ inside some interval of $S$. This limit exists by (1.02) if $\bar S$ satisfies (vii), which we will suppose. When $y$ satisfies (a) of property (p) in I \textsection 3 (the only interesting case), and when $j(y)>0$, we say that $y$ \textit{is a singularity of $g(I)$ with respect to its variation and with respect to $G$}, or a \textit{variation singularity.}

Obviously, by (1.03),

\vs
\noindent
\textbf{(1.08)}
\textit{if $j(y)=0$ and if both ends of $I$ tend to $y$, then $A(I) \rightarrow 0$.}

\vs
\noindent
\textbf{(1.09)}
\textit{The variation singularities of $g(I)$ are at most enumerable, and if they are $v_1, \ldots ,v_n, \ldots $ then (for $W$ in $\bar S$)} $\sum_{n=1}^\infty j(v_n) \leq \var(g;W;G)$.

\noindent
Let $I_1, \ldots , I_n$ be non-overlapping intervals of $\bar S$. Then by I(2.11a), and then by (1.02),
\[
\sum_{i=1}^n \var (g;I_i;G) \leq \var (g;\sum_{i=1}^n I_i;G) \leq \var (g;W;G).
\]
Fix $v_i$ in $I_i^o$ and let $mI_i \rightarrow 0$ $(=1, \ldots ,n)$. Then $\sum_{n=1}^n j(v_i) \leq \var(g;W;G)$.
The Cantor argument completes the proof (cf.~I(3.07)).

We will now consider $k'$-integration, and first a result for $A(R)$.

\vs
\noindent
\textbf{(1.10)}
\textit{$(k';G)\int_R \,A(I)$ exists and is not less than $A(R) - (R^o)\sum c(y;G)$.}

\noindent
From (1.07) we have $\;\;\;A(R_3) -((R_3')^o)\sum c(y;G) \;\;\leq$
\[
\leq A(R_1) -((R_1')^o)\sum c(y;G) 
+ A(R_2) -((R_2')^o)\sum c(y;G)
\]
so that
$\leq A(I) -(I^o)\sum c(y;G)$ increases on subdivision, in the sense of \textsection 3.
But if
$
h(I) =(I^o)\sum c(y;G)$, then $(k';G)\int_R \,h(I)=0$
for all bracket conventions. Hence by (3.3) we have the results.

By (1.08) the permanent points of the integration are included in the variation singularities.

\vs
\noindent
\textbf{Corollary 1.}
\[
(k';G)\int_R \,A(I)\geq \var(A;R;G) -(R^o)\sum c(y;G),
\]
so that from (1.05),
\[
 \var(g;R;G)\leq (k';G)\int_R \,A(I) + 3(R^o)\sum c(y;G).
\]
\noindent
\textbf{Corollary 2.}
\textit{If $g(I)$ is continuous then 
\[ (k';G)\int_R \,A(I)= \var(A;R;G)= \var(g;R;G). \] }

\noindent
\textbf{(1.11)}
\textit{Let $x<y<z$, $I=x$---$z$, $I_1=x$---$z$, $I_2=y$---$z$, where $I_1, I_2$ are in $S$. Then as $x,z$ tend independently to $y$,
\[
j(y) = \max\left(\overline{\lim}|g(I)|,\;\overline{\lim}\left(|g(I_1)|+|g(I_2)|\right)\right).
\]
}
Since $\var(g;J;G)$ is finite and tends to a limit when $J=I$ or $I_1$ or $I_2$, the result follows from I(3.12).

\vs
\noindent
\textbf{Corollary 1.}
\textit{$c(y;G) \leq 2j(y)$, so that each singularity of $g(I)$ with regard to additivity is also a variation singularity.}

\vs
\noindent
\textbf{Corollary 2.}
\textit{If $j(y)=0$ and if $y$ is fixed in $I'$ with $mI\rightarrow 0$, then $g(I) \rightarrow 0$.}

\vs
\noindent
\textbf{Corollary 3.}
\textit{If $g(I)$ is continuous then $j(y)=0$.}

\vs
\noindent
The result (1.11 Corollary 1) can be proved without assuming that $G$ satisfies (viii). But we have to use (vi) and (vii). For, using the notation of Chapter 1 \textsection 3,

\begin{eqnarray*}
\left|g(J) -g(J_1) - g(J_2)\right| &=&
\left| (I;D)\sum g(I) -(I;D_1+D_2) \sum g(I)\right| \vt
&\leq & 2 \max \left\{
 (I;D)\sum \left|g(I)\right|  ,\;\;(I;D_1+D_2) \sum\left| g(I)\right| \right\}
 \end{eqnarray*}
since $|A-B| \leq 2 \max\left\{|A|,\;|B|\right\} $.
Hence $C(x;G;e) \leq 2d_1(e)$, where (in the notation of Chapter 1 \textsection 2), $d_1(e)$ is the $d(e)$ for $|g(I)|$. Hence
\[
c(x;G) \leq 2(N;G) \overline{\int}_I\, |g(I)|,
\]
and this for all $I$ with $x$ inside $I$. Hence $c(x;G) \leq 2j(x)$.

When $g(I)$ is not b.v., the equality in (1.11) can be false. For example, let $g(I)=1$ for $I=2^{-n}$---$2^{-n+1}$ and $n=1,2, \ldots$, and otherwise $g(I)=0$. Then $j(0) =+\infty$. But the three upper limits are zero.

From (1.03 Corollary), the upper and lower norm-limits of $g$ over $R$ are finite, so that we may define the $k'$-limits. From (1.11 Corollary 1) and I(4.04) we obtain:

\vs
\noindent
\textbf{(1.12)}
\textit{We do not obtain new limits for $g(I)$ alone if we take the variation sing\-ularities as permanent points.}

\vs
\noindent
\textbf{(1.13)}
\textit{As $z\rightarrow y$, with $z$---$y$ in $\bar S$, $(N;G) \overline{\int}_{(z,y)}\,g(I)$ tends to a limit, say $\bar N(y-)$, which is zero unless $y$ is a variation singularity. Similarly}
\[
(k';G) \overline{\int}_{(z,y)}\,g(I) \rightarrow  {\bar {k}}'(y-),\;\;\;\;\;\;A((z,y)) \rightarrow A(y-).
\]
\textit{Similarly when $z>y$ and $z \rightarrow y$;
when $x<y<z$ and $x,z \rightarrow y$; and for lower limits.
}

Let $z_1< \cdots <z_n<y$ and $z_n \rightarrow y$ as $n \rightarrow \infty$ such that $I_n =(z_n,z_{n+1})$ and $(z_1,y)$ are all in $\bar S$. Then by (vii), all $(z_n,y)$ are in $\bar S$. As $n \rightarrow \infty$,
\[
\sum_{i=1}^n (N;G)\overline{\int}_{I_i} \,g(I) +
(N;G)\overline{\int}_{(z_{n+1},y)} \,g(I)
\] 
is monotone decreasing by I(2.11), and by (1.03) its modulus is $\;\leq$
\[
\leq\; \sum_{i=1}^n \var(g;I_i;G) + \var(g;(z_{n+1},y);G), \;\;\leq \var(g;W;G)
\]
by I(2.11). Thus
$
\sum_{i=1}^n (N;G)\overline{\int}_{I_i} \,g(I) +
(N;G)\overline{\int}_{(z_{n+1},y)} \,g(I)
$
tends to a limit, $P$ say. Also for every $n$ we have by (1.03) and I(2.11),
\[
 \sum_{i=1}^n (N;G)\overline{\int}_{I_i} \,g(I) \leq \var(g;W;G),
 \]
so that the series $\sum_{i=1}^\infty (N;G)\overline{\int}_{I_i} \,g(I)$ is convergent, with sum $Q$ say. Then
\[
(N;G)\overline{\int}_{(z_{n+1},y)} \,g(I) \rightarrow P-Q
\]
as $n \rightarrow \infty$, and this for each sequence $z_1<\cdots <z_n<y$ with $z_n \rightarrow y$ and with $z_n$---$z_{n+1}$ and $z_1$---$y$ all in $\bar S$. Hence $P-Q$ is independent of the choice of $\{z_n\}$, so that we have the first result. By (1.08), $P-Q=0$ unless $j(y)>0$.
Similar proofs hold for the other results.

\vs
\noindent
\textbf{(1.14)}
\textit{If $z\rightarrow y$ and $I=z$---$y$ is in $S$ then
\[
\overline{\lim} \,g(I) = {\overline N}(y-),\;\;\;\;\;\;
\underline{\lim}\, g(I) = {\underline{N}}(y-).
\]
There are many similar results.}

\noindent
By I(3.12) and (3.13)

\vs
\noindent
\textbf{Corollary.}
\textit{When $g(I)$ is also continuous then so are the upper and lower norm-limits.}
\noindent
(Or use (1.11 Corollary 3) and (1.08).)

We can split $g(I)$ up into an additive part, and a part whose upper norm-limit is continuous to the left and to the right. For let $B^-(R)$, $B^+(R)$ be the respective additive functions formed from
\[
B^-(a\mbox{---}b)=\left((a,b)\right)\sum \bar N(v-),\;\;\;\;\;\;B^-(a\mbox{---}b)=\left([a,b)\right) \sum \bar N(v+),
\]
and put $g_1(I) = g(I) -B^-(I)-B^+(I)$. Then

\noindent
\textbf{(1.15)}
\[
(N;G)\overline{\int}_R\,g(I) = (N;G)\overline{\int}_R\,g_1(I) + B^-(R) +B^+(R),
\]
\textit{where the second limit is left and right continuous.}

\noindent
For $(R) \sum g(I) =  B^-(R) +B^+(R)+(R)\sum g_1(I)$, giving the first results. And
\[
\lim_{z \rightarrow y} (N;G)\overline{\int}_{(z,y)}\,g_1(I)
=
\lim_{z \rightarrow y} \left( (N;G)\overline{\int}_{(z,y)}\,g(I) -B^-((z,y))-B^+((z,y))\right) =0.
\]
Similarly $\lim_{z \rightarrow y} (N;G)\overline{\int}_{(y,z)}\,g_1(I)=0$.
Similar decompositions can be given of the lower norm-limit and the $k'$-limits.

Let $g(I)$ be independent of bracket conventions, and let $G$ satisfy (viii). We then define the variations of $g(I)$, as follows.

\textit{The upper positive variation $\overline{p}(J)$ of $g(I)$ in $J$ of $S$}
is l.u.b.$(S_2)\sum g(I)$ for every finite non-overlapping set $S_2$ of intervals of $S$ in $J$, and \textit{the lower negative variation $\underline{n}(J)$ of $g(I)$ in $J$ of $S$ is $-\mbox{g.l.b.}(S_2)\sum g(I)$.} The \textit{upper negative variation $\overline{n}(J)$ of $g(I)$ in $J$ of $S$ is $\overline{p}(J) -g(J)$} and the \textit{lower positive variation $\underline{p}(J)$ of $g(I)$ in $J$ of $S$} is $\underline{n}(J) +g(J)$, so that
\[
g(J) =\overline{p}(J) - \overline{n}(J)
= \underline{p}(J) - \underline{n}(J).
\]
By (1.01) the variations are finite for $mJ<\delta$. And by (1.09) and (1.11 Corollary 2), $\overline{p}(I)\geq 0$ and $\underline{n}(I)\geq 0$.
Obviously $\underline{p}(I)\geq 0$ and $\overline{n}(I)\geq 0$.

\vs
\noindent
\textbf{(1.16)}
\textit{The $k'$-limits of $\overline{p}(I)$ and $\underline{n}(I)$ exist.}

\noindent
Let $\overline{P}(I)$ be the $\overline{p}(I)$ for $h(I) = (k';G) \overline{\int}_I\,g(I)$. Then $h(I)$ is additive. Let $I_1,I_2,I_3=I_1+I_2'$, be intervals of $\bar S$ with $I_1^o.I_2^o$ null. Then obviously 
\[
\overline{P}(I_1) + \overline{P}(I_2) \leq \overline{P}(I_3) \leq\overline{P}(I_1)+\overline{P}(I_2) + \mbox{lu.b.} \left| h(J) -h(J_1) -h(J_2)\right|
\]
where $J$ in $I_3$ splits up into $J_1$ in $I_1$ and $J_2$ in $I_2$, and $J_1, J_2$ are in $\bar S$. But $h(J)=h(J_1)+h(J_2)$, so that $\overline{P}(I)$ is additive.

Given $\ve>0$ there is an $S_2$ (with the corresponding $I^\sigma, R_2$ in $\bar S$) such that 
\[
(k';G)\overline{\int}_{R_2} \,g(I) >\overline{P}(I) -\ve\;\;\;\mbox{ and }\;\;\;R_2 \subset I.
\]
There is then a division $D$ in $G$ of $R_2$ such that
\[
(R_2;D)\sum g(I) >
(k';G)\overline{\int}_{R_2} \,g(I)-\ve.
\]
Hence for each $\ve>0$, $\overline{p}(I) > \overline{P}(I)-2\ve$, so that $\overline{p}(I) \geq \overline{P}(I)$.

Now let $D$ be a $k'_{me}$-division of $J$. Then by I(4.12),
\begin{eqnarray*}
\overline{P}(J)& =&(J;D) \sum \overline{P}(I) 
\leq (J;D) \sum \overline{p}(I) 
\leq (J;D) \sum \overline{P}(I) +\delta, \vt
\overline{P}(J) &\leq& (J;D) \sum \overline{p}(I) \leq 
\overline{P}(J) +\delta.
\end{eqnarray*}
 Hence $(k';G)\overline{\int}_{J} \,\overline{p}(I) =\overline{P}(J)$.
Similarly 
\[(k';G)\overline{\int}_{J} \,\underline{n}(I) =\underline{N}'(J),\] the $\underline{n}(J)$ for 
$(k';G)\overline{\int}_{J} \,g(I)$.

\vs
\noindent
\textbf{Corollary.}
\begin{eqnarray*}
(k';G)\overline{\int}_{J} \,g(I) &=& \overline{P}(J) -
(k';G)\underline{\int}_{J} \,\overline{n}(I), \vt
(k';G)\underline{\int}_{J} \,g(I) &=&
(k';G)\underline{\int}_{J} \,\underline{p}(I)
- \underline{N}'(J).
\end{eqnarray*}
If $g(I)$ depends on the bracket convention of $I$, no such simple decompositions can be given,

It is well-known that a function such as $\overline{P}(I)$ is additive, but the simple proof is included for completeness.

We denote by $E_+'$ the set of $x$ such that there is a point of $E$ in every right-hand neighbourhood $(x,y)$ of $x$, and by $E_-'$ the set of $x$ such that there is a point of $E$ in every left-hand neighbourhood $(y,x)$ of $x$.

For $x$---$v$ and $v$---$x$ in $S$ let $g(x$---$v) \rightarrow g(v-)$ and $g(v$---$x)\rightarrow g(v+)$ as $x \rightarrow v$ with a chosen $k'$-convention at the singular points $v$.  We then put
\[
A^-(R) = (R_-' )\sum g(v-),\;\;\;
A^+(R) = (R_+' )\sum g(v+),\;\;\;
g_2(I) = g(I) - A^-(I) - A^+(I).
\]
As is easily seen, if $R'$ is in $R_1^o$, and if $F$ satisfies (vi), (vii), (viii) (which we suppose),
\[
(R')\sum \left(|g(v-)| + |g(v+)| \right) \leq \var(g;R_1;G).
\]
Hence the series for $A^-(R)$ and $A^+(R)$ are convergent.

Let $h(I)$ be bounded by $K$ in $S$. Then by (1.11), the $c(y;G)$ for $h(I)g(I)$ is not greater than $2Kj(y)$, so that the $\{y_i'\}$ for $h(I)g(I)$ is included in the $\{v_i\}$ for $g(I)$. We put
\begin{eqnarray*}
\overline{\lim}\,h(x\mbox{---}v) &=& \left\{
\begin{array}{rllll}
\overline{h}(v-)&\mbox{when}& g(v-) &\geq & 0, \vt
\underline{h}(v-)&\mbox{when}& g(v-) &<& 0,
\end{array}
\right. \vt 
\underline{\lim}\,h(x\mbox{---}v) &=& \left\{
\begin{array}{rllll}
\underline{h}(v-)&\mbox{when}& g(v-) &\geq & 0, \vt
\overline{h}(v-)&\mbox{when}& g(v-) &<& 0;
\end{array}
\right.
\end{eqnarray*}
and similarly for $\overline{h}(v+),\underline{h}(v+)$. Then

\vs
\noindent
\textbf{(1.17)}
\begin{eqnarray*}
(k';G)\overline{\int}_{R} \,h(I)g(I) &=&
(k';G)\overline{\int}_{R} \,h(I)g_2(I)\;\;+\vt
&&\;\;\;\;\;\;\;\;+\;\;(R_-')\sum \overline{h}(v-)g(v-)
+(R_+')\sum \overline{h}(v+)g(v+),\vt
(k';G)\underline{\int}_{R} \,h(I)g(I) &=&
(k';G)\underline{\int}_{R} \,h(I)g_2(I)\;\;+\vt
&&\;\;\;\;\;\;\;\;+\;\;
(R_-')\sum \underline{h}(v-)g(v-)
+(R_+')\sum \underline{h}(v+)g(v+).
\end{eqnarray*}
We have

\noindent
\textbf{(1)}
\[
\sum h(I)g(I) =\sum h(I)g_2(I) + \sum h(I)A^-(I) +\sum h(I)A^+(I) .
\]
Let $m$ be such that

\noindent
\textbf{(2)}
\[
\sum_{n=m+1}^\infty
\left(|g(v_n-)|+|g(v_n+)|\right) < \frac \ve K.
\]
Then if $x_n$---$v_n$ and $v_n$---$z_n$ are intervals in the $k'_{me}$-division $D$ $(n=1, \ldots ,m)$, and if $S_2$ is the set of the rest,

\noindent
\textbf{(3)}
\begin{eqnarray*}
&&\;\;\;\;\;\;\left| (R;D) \sum h(I) A^-(I) +(R;D) \sum h(I) A^+(I)
\right. \;\;- \vt
&&
\;\;\;\;\;\;\;\;\;\;\;\;\;\;\;\;
\;\;\;\;\;\;\;\;-\;\;
\sum_{n=1}^m \left(h(x_n\mbox{---}v_n)g(v_n-) +
\left. h(v_n\mbox{---}z_n)g(v_n+)  \right) \right| \vt
&&\leq \;
 (S_2) \sum \left|h(I) \left(A^-(I) + A^+(I)\right)\right|
+ \sum_{n=1}^m \left|h(x_n\mbox{---}v_n)
\sum_{\stackrel{i>m}{x_n<v_i<v_n}}g(v_i-)\right| \;\;+\vt
&&
\;\;\;\;\;\;\;\;\;\;\;\;\;\;\;\;
\;\;\;\;\;\;+\;\;\sum_{n=1}^m \left|
h(v_n\mbox{---}z_n)\sum_{\stackrel{i>m}{v_n<v_i<z_n}}g(v_i+)   \right|, \vt
&&< \;K.\frac \ve K\;=\;\ve,\;\;\;\mbox{ by (2)}.
\end{eqnarray*}
Let $D$ be such that

\noindent
\textbf{(4)}
\[
(R;D)\sum h(I) g_2(I) >(k';G)\overline{\int}_{R} \,h(I)g_2(I)-\ve.
\]
\textbf{(5)}
Taking $k'$-conventions into account we know that $g_2(I)$ is continuous by the definitions of $A^-(I)$;
and $g_2(I)$ is b.v.~since $g(I), A^-(I)$ and $A^+(I)$ are b.v. Hence by (1.14) for $k'$-limits, the upper and lower $k'$-limits of $|g_2(I)|$ are continuous to the left and right at every $v$.

We now take $x_n<x_n'< v_n$ so that $x'_n$---$v_n$ is in $S$, and so by (vii), $x_n$---$x_n'$ is in $\bar S$; and

\noindent
\textbf{(6)}
\[
h(x_n'\mbox{---}v_n) g(v_n-) >\bar{h}(v_n-)g(v_n-)-\frac \ve n,
\]
and we replace $x_n$---$v_n$ in $D$ by $x'_n$---$v_n$ and a division $D'$ of $x_n$---$x'_n$; and similarly for $v_n$---$z'_n$ and $z'_n$---$z_n$.

We take $D'$ as a $k'_{me}$-division such that

\noindent
\textbf{(7)}
\[
((x_n,x'_n);D')\sum | g_2(I)| <(k';G)\overline{\int}_{(x_n,x'_n)} \,|g_2(I)|  +\frac \ve{Km}.
\]
By (5) and I(4.03) we can choose $x_n$ so near to $v_n$ that

\noindent
\textbf{(8)}
\[
(k';G)\overline{\int}_{(x_n,x'_n)} \,|g_2(I)|  < \frac \ve{Km}\;\;\;\mbox{ and }\;\;\;
|g_2(x\mbox{---}v_n)| < \frac \ve{Km}\;\;\;\;\;\;(x=x_n, x_n').
\]
\textbf{(9)}
By (7) and (8), replacing $x_n$---$v_n$ in $D$
by $x'_n$---$v_n$ and $D'$ over $x_n$---$x'_n$, and similarly for $v_n$---$z_n$ (for $n-1, \ldots ,m$) to form $D''$, we alter $(R;D)\sum h(I)g_2(I)$ by at most $8\ve$.

By (viii) and (iii), $D''$ is a $k'_{me}$-division of $R$, and we have by (1), (3), (9), that $(R;D'')\sum h(I)g(I) \;\;=$
\begin{eqnarray*}
&=&
(R;D'')\sum h(I)g_2(I) + (R;D'')\sum h(I)A^-(I) +(R;D'')\sum h(I)A^+(I) \vt
&>&
(R;D)\sum h(I)g_2(I) - 9\ve
+\sum_{n=1}^m \left(
h(x_n'\mbox{---}v_n)g(v_n-) + h(v_n\mbox{---}z'_n)g(v_n+)
\right),
\end{eqnarray*} 
and then by (4), (6),
\begin{eqnarray*}
&>&
(k';G)\overline{\int}_{(x_n,x'_n)} \,h(I)g_2(I) - 12\ve 
+
\sum_{n=1}^m \left(
\bar{h}(v_n-)g(v_n-) + \bar{h}(v_n+)g(v_n+)
\right),
\end{eqnarray*} 
i.e. $(R;D'')\sum h(I)g(I)\;\;>$
\[
>\;\;(k';G)\overline{\int}_{R} \,h(I)g_2(I) 
+
\sum_{n=1}^\infty \left(
\bar{h}(v_n-)g(v_n-) + \bar{h}(v_n+)g(v_n+)
\right) -  13\ve ,
\]
where the summation is over every significant $v_n$ in $R$. Hence letting $m \rightarrow \infty$, $e \rightarrow 0$, and then $\ve \rightarrow 0$, $\;\;\;(k';G)\overline{\int}_{R} \,h(I)g(I) \;\;\geq $
\[
\geq \;\;
(k';G)\overline{\int}_{R} \,h(I)g_2(I) 
+
(R'_-)\sum 
\bar{h}(v-)g(v-) + (R'_+)\sum\bar{h}(v+)g(v+)
. 
\]
The reverse inequality is obvious from (1), and hence the result.

Note that $(k';G)\overline{\int}_{R} \,h(I)g_2(I) $ is continuous, so that we have separated the expression $(k';G)\overline{\int}_{R} \,h(I)g(I) $ into its ``continuous and jump functions''.

When $g(I)=D(f;I)$ we can obtain a similar theorem for Riemann-Stieltjes integration. (Dienes [14].)

In two dimensions we can obtain a similar theory. The following result is immediately obtained from the two-dimensional analogue of (1.01).

\vs
\noindent
\textbf{(1.18)}
\textit{For $T=[I_1,I_2]$ put $g(T)=g(I_1,I_2)$. Then when $mI_2<\delta$, $g(I_1,I_2)$ is of bounded variation with respect to $I_1$, with variation not greater than $M+\ve$.}

Of the examples in Chapter 1 \textsection 1, the following are b.v.: (i), (vi), (x), (xiv), (xvi), (xviii). When $f(x)$ is b.v.~then also so are (ii), (iii), (iv), (xi), (xii).

Of the examples in Chapter 2 \textsection 1, the following are b.v.: (i) ($m(T)$, but not $p(T)$ nor $\delta(T)$), (v). When $f(x,y)$ is b.v.~then also so are (ii), (iv) (when $g(T)$ has form (ii)\,), (vi), (vii), (viii).

\section{Absolutely continuous functions.}
\label{Absolutely continuous functions}
Let $S_2$ be a finite sum of non-overlapping intervals of $S$ in $R$, and covering $R_2$. Then $g(I)$ is \textit{absolutely continuous} in $R$ if $(S_2)\sum g(I) \rightarrow 0$ as $mR_2 \rightarrow 0$, for all such $S_2$. (Burkill [6] \textsection 3. See also Kempisty [12] 15.)

Obviously $g(I)$ is continuous, and $|g(I)|$ is absolutely continuous. If $g_1(I), g_2(I)$ are absolutely continuous so are $g_1+g_2$, $g_1-g_2$. If $g_1(I)$ is absolutely continuous and $g_2(I)$ is bounded then $g_1(I)g_2(I)$ is absolutely continuous.
Burkill writes
\[
p(I)=\frac{|g(I)|+g(I)}2,\;\;\;\;\;\;n(I) =\frac{|g(I)|-g(I)}2,
\]
so that $p(I)$ and $n(I)$ are also absolutely continuous, and
\[
g(I)=p(I) - n(I),\;\;\;\;\;\;\;\;|g(I)| = p(I)+n(I).
\]
We now show that

\vs
\noindent
\textbf{(2.1)}
\textit{If $g(I)$ is absolutely continuous in $R$ then each of 
\[
(N;G)\overline{\int}_{R} \,g(I),\;\;\;\;\;\; 
(N;G)\underline{\int}_{R} \,g(I),\;\;\;\;\;\;
(N;G)\overline{\int}_{R} \,|g(I)|
\]
is finite.
}

\noindent
Since
\[
(N;G)\overline{\int}_{R} \,g(I) \leq
(N;G)\overline{\int}_{R} \,p(I),\;\;\;\;\;\;
(N;G)\underline{\int}_{R} \,g(I) \geq -(N;G)\overline{\int}_{R} \,n(I),
\]
and
\[
(N;G)\overline{\int}_{R} \,g(I)\leq
-(N;G)\overline{\int}_{R} \,p(I)+ (N;G)\overline{\int}_{R} \,n(I),
\]
it is sufficient to show that if $g(I) \geq 0$ is absolutely continuous then $(N;G)\overline{\int}_{R} \,g(I)$ is finite.

Suppose not. Given $\ve>0$ there is a $\delta>0$ such that $(S_2)\sum g(I)<\ve$ when $mR_2<\delta$. Choose $N> 1 +(mR) (\ve/\delta)$. Then there is in $G$ a $D$ over $R$ with norm$(D)< \delta/N$ and
\[
(R;D) \sum g(I) >N+\ve.
\]
Group the intervals of $D$ to form sets $S_{2,i}$ with $mR_{2,i}<\delta$ and, as far as possible, with $mR_{2,i} >(N-1) \delta /N$. Thus there may be one group of meshes having measure $\leq (N-1)\delta/N$. The number of these groups $S_{2,i}$ is less than
\[
\frac{N.mR}{(N-1) \delta} +1
\]
since $mR_{2,i}<\delta$, $(S_{2,i})\sum g(I) < \ve$. Hence
\[
(R;D) \sum g(I) < \left(\frac{N.mR}{(N-1)\delta}+1\right) \ve <N+\ve.
\]
This gives a contradiction. Hence the result.

\vs
\noindent
\textbf{Corollary.}
\textit{If $g(I)$ is absolutely continuous it is of bounded variation.}

\noindent
Hence, the results of \textsection 1 for a function which is b.v.~and continuous all apply here. There is also a further result.

\vs
\noindent
\textbf{(2.2)}
\textit{If $g(I)$ is absolutely continuous then so are}
$(N;G)\overline{\int}_{R} \,g(I)$, $(N;G)\underline{\int}_{R} \,g(I)$.
This is obvious, and can be used to prove (2.1). These results are easily extended to two or more dimensions. Note the theorem (3.5 Corollary).

Kempisty [12] 15 also defines absolutely semi-continuous functions:
\begin{itemize}
\item
$g(I)$ is \textit{SAC (upper absolutely semi-continuous)}
if given $\ve>0$, then all $(S_2)\sum g(I)>-\ve$ for $mR_2<\delta = \delta(\ve)$.
\item
$g(I)$ is \textit{SAC (upper absolutely semi-continuous}
if given $\ve>0$, then all $(S_2)\sum g(I)<\ve$ for $mR_2<\delta = \delta(\ve)$.
\end{itemize}
\textbf{(2.3)}
\textit{If $g(I)$ is IAC (SAC) then so are 
$(N;G)\overline{\int}_{R} \,g(I)$, $(N;G)\underline{\int}_{R} \,g(I)$.
}

\vs
\noindent
\textbf{(2.4)}
\textit{If $g(I)$ is IAC then  $(N;G)\underline{\int}_{R} \,g(I)> -\infty$. If $g(I)$ is SAC then  \[(N;G)\overline{\int}_{R} \,g(I)> -\infty.\]
}

\noindent
We use a similar proof to that of (2.1).

Of the examples in Chapter 1, \textsection 1, the example (i) is absolutely continuous. When $f(x)$ is absolutely continuous then so are (ii), (iii), (iv), (x), (xi), (xii).

Of the examples in Chapter 2, \textsection 1, the example (i) ($m(T)$, but not $p(T)$ nor $\delta(T)$) is absolutely continuous. When $f(x,y)$ is absolutely continuous, then so are (ii), (iv) (when $g(T)$ has the form (ii)), (vi), (vii), (viii).

\section{Functions which are monotone on subdivision.}
\label{Functions which are monotone on subdivision}
Suppose that $g(I_1)+g(I_2) \leq g(I_3)$ whenever $I_1$ and $I_2$ are non-overlapping, $I_3'=I_1'+I_2'$, and $I_1,I_2,I_3$ are in $S$. Then $g(I)$ is said to \textit{decrease on subdivision}. If $g(I_1)+g(I_2) \geq g(I_3)$ always then $g(I)$ \textit{increases on subdivision}. (Burkill [6] \textsection 5.)

\vs
\noindent
\textbf{(3.1)}
\textit{If $g(I)$ decreases on subdivision then for each interval $J$ in $S$, \[(N;G)\overline{\int}_{J} \,g(I) \leq g(J).\] Similarly $(N;G)\underline{\int}_{J} \,g(I) \geq g(J)$ when $g(I)$ increases on subdivision.}

\vs
\noindent
\textbf{(3.2)}
\textit{If $g(I)$ decreases on subdivision, with an at most enumerable number of singularities $y_1', y_2', \ldots $, and $\sum_{i=1}^\infty c(y_i';G)<\infty$, then either
\[
(N;G)\overline{\int}_{R} \,g(I)  =-\infty\;\;\;\;\;\mbox{ or }\;\;\;\;\;\mbox{osc}N(g;R;G) = (R^o)\sum c(y';G).
\]}
(Burkill [6] theorem 5.1, when $g(I)$ is continuous and $G=G_1$.) Let the division $D$ with division points $FR$ and $x_1, \ldots ,x_n$ be such that $(R;D)\sum g(I) <A$ where $A>(N;G)\underline{\int}_{R} \,g(I) $. Then by (vi) and a proof similar to I(3.08),
\begin{eqnarray*}
(N;G)\overline{\int}_{R} \,g(I)&\leq &
(R;D)\sum (N;G)\overline{\int}_{I} \,g(I) +\sum_{i=1}^n c(x_i;G), \vt
&\leq & (R;D)\sum g(I) + \sum_{i=1}^n c(x_i;G) \vt
&<& A + \sum_{i=1}^n c(x_i;G) \;\;<\;\; A+(R^o)\sum c(y';G).
\end{eqnarray*}
If $(N;G)\underline{\int}_{R} \,g(I)=-\infty$ then we can take $A \rightarrow -\infty$, and since $\sum_{i=1}^\infty c(y'_i;G)$ is convergent we then have 
\[
(N;G)\overline{\int}_{R} \,g(I) \leq(N;G)\underline{\int}_{R} \,g(I) + (R^o)\sum c(y';G),
\]
i.e.~osc$N(g;R;G) \leq (R^o)\sum c(y';G)$. The equality follows from I(3.07).

\vs
\noindent
\textbf{Corollary.}
\textit{When $g(I)$ is also continuous then either $(N;G)\overline{\int}_{R} \,g(I)=-\infty$ or 
osc$N(g;R;G) \leq (R^o)\sum c(y';G)=0$, i.e.
$(N;G){\int}_{R} \,g(I)$ exists.
}

\vs
Similarly we may prove

\vs
\noindent
\textbf{(3.3)}
\textit{If $g(I)$ increases on subdivision, with an at most enumerable number of singularities $y_1', y_2', \ldots $, and $\sum_{i=1}^\infty c(y'_i;G)<\infty$, then either 
\[
(N;G)\underline{\int}_{R} \,g(I)  =-\infty\;\;\;\;\;\mbox{ or }\;\;\;\;\;\mbox{osc}N(g;R;G) = (R^o)\sum c(y';G).
\]
}
\textbf{(3.4)}
\textit{If $g$ decreases on subdivision and 
$(N;G)\underline{\int}_{W} \,g(I)  >-\infty$, then the $k'$-limit of all $g(I)$ exists, where $k'$ is \textbf{all} conventions, and $(k';G)\underline{\int}_{J} \,g(I) \leq g(J)$ for $J$ in $S$.
}

\noindent
By (3.1), $(N;G)\overline{\int}_{W} \,g(I)  <+\infty$ (assuming $W$ in $\bar S$) so that we can define the $k'$-limit. Then $(k';G)\overline{\int}_{J} \,g(I)  \leq g(J)$.
Let the division $D$ in $G$ over $R$ in $\bar S$ be such that
$
(R;D) \sum g(I) < (k';G)\underline{\int}_{R} \,g(I) +\ve.
$
Then
\begin{eqnarray*}
(k';G)\overline{\int}_{R} \,g(I) &=&
(R;D) \sum (k';G)\overline{\int}_{R} \,g(I) \;\;\mbox{ by I(4.05)},\vt
&\leq & (R;D) \sum g(I) \;\;<\;\;(k';G)\underline{\int}_{R} \,g(I) +\ve.
\end{eqnarray*}
Hence $(k';G){\int}_{R} \,g(I)$ exists, and then is not greater than $g(J)$ when $R=J$ in $S$.
Similarly when $g$ increases on subdivision and $(N;G)\overline{\int}_{W} \,g(I)<+\infty$. 

We may generalise the property of increasing on subdivision, or decreasing on subdivision, to two dimensions. 

Suppose that $(J;D)\sum g(T) \leq g(J)$ whenever $J$ is in $U$, and $D$ in $L$ is a division of $J$. Then $g(T)$ \textit{decreases on subdivision with respect to} $L$. Suppose that $(J;D)\sum g(T) \geq g(J)$ whenever $J$ is in $U$, and $D$ in $L$ is a division of $J$. Then $g(T)$ \textit{increases on subdivision with respect to} $L$.

\vs
\noindent
\textbf{(3.5)}
\textit{If $g(T)$ decreases on subdivision, and if $\bar{c}(V;L)$ is finite, where  $\bar{c}(V;L)$ corresponds to the  $\bar{b}(V;H)$  of II(3.10) for $(H)\overline{\int}_{V} \,g(T)$, then either 
\[
(N;L)\overline{\int}_{V} \,g(T) =-\infty\;\;\;{ or }\;\;\;
\mbox{osc}N(g;V;L) = \bar{c}(V;L).
\]
}
(Burkill [6] theorem 5.1, when $g(T)$ is absolutely continuous and $L=L_0$.) For the proof, follow that of (3.2), using the analogue of II(3.10, Corollary 1,2).

\vs
\noindent
\textbf{Corollary.}
\textit{When $g(T)$ is also absolutely continuous, the norm-limit of $g(T)$ exists.}

\noindent
For $g(T)$ is then b.v., so that osc$N(g;V;L)=\bar{c}(V;L), =0$ by absolute continuity.

Of the examples in Chapter 1, \textsection 1, the examples (i), (ii), (iv), (vi), (ix), (xiv), (xvi) are increasing on subdivision, and (i), (ii), (x), (xvi) are decreasing on subdivision.

Of the examples in Chapter 2, \textsection 1, the examples (i) ($m(T)$, $\delta(T)$; and $p(T)$ at least for restricted dimensions), (ii), (iv) (when $g(T)$ has the form (ii)), (v), (vii), (viii) are increasing on subdivision. The examples (i) ($m(T)$), (ii), (v) are additive and so decrease on subdivision also. 

\section{$g(T)=g_1(I_x).g_2(I_y)$ when $T=[I_x,I_y]$ in 2-dimensions.}\label{product}
\textbf{(4.1)}
\textit{Let the family $H_r$ produce the families $F_x$ and $F_y$ by projection on the axes. Then for a rectangle $T=[I_x, I_y]$, we have
\begin{eqnarray*}
(H_r) \overline{\int}_T\,g(T)& \leq & (F_x) \overline{\int}_{I_x}\,g_1(I).(F_y)\overline{\int}_{I_y}\,g_2(I), \;\;\;\mbox{ and} \vt
(N;L_r) \overline{\int}_T\,g(T)& \leq & (N;G_x) \overline{\int}_{I_x}\,g_1(I).(N;G)\overline{\int}_{I_y}\,g_2(I).
\end{eqnarray*}}
For if $D$ is in $L_r$, and $D_x,D_y$ are the projections on the axes of $x$ and $y$,
\[
(I_x;D_x) \sum g_1(I).(I_y;D_y) \sum g_2(I)=
(T;D) \sum g(T).
\]
Hence the results.

\noindent
\textbf{(4.2)}
\textit{For each Riemann succession $\{D_x^{(n)}\}$ in $F_x$ and over $I_x$, and for each Riemann succession $\{D_y^{(n)}\}$ in $F_y$ and over $I_y$,
let there be a Riemann succession $\{D^{(n)}\}$ in $H_r$ and over $T=[I_x,I_y]$ such that $D^{(n)}$ projects into $D^{(n)}_x, D^{(n)}_y$ on the $x$-axis and $y$-axis. Then
\begin{eqnarray*}
(H_r) \overline{\int}_T\,g(T)& \geq & (F_x) \overline{\int}_{I_x}\,g_1(I).(F_y)\overline{\int}_{I_y}\,g_2(I), \;\;\;\mbox{ and} \vt
(N;L_r) \overline{\int}_T\,g(T)& \geq & (N;G_x) \overline{\int}_{I_x}\,g_1(I).(N;G_y)\overline{\int}_{I_y}\,g_2(I).
\end{eqnarray*}}
\textbf{Corollary.}
\textit{If the hypothesis of (4.1) holds also, i.e.~if $H=\mbox{``}F_x.F_y\mbox{''}$, then}
\begin{eqnarray*}
(H_r) \overline{\int}_T\,g(T)& = & (F_x) \overline{\int}_{I_x}\,g_1(I).(F_y)\overline{\int}_{I_y}\,g_2(I), \;\;\;\mbox{ and} \vt
(N;L_r) \overline{\int}_T\,g(T)& = & (N;G_x) \overline{\int}_{I_x}\,g_1(I).(N;G_y)\overline{\int}_{I_y}\,g_2(I).
\end{eqnarray*}
Given $D_x$ in $G_x$, $D_y$ in $G_y$, we obtain $D$ in $L_r$
such that
\[
(T;D) \sum g(T) =(I_x;D_x)\sum g_1(I).(I_y;D_y)\sum g_2(I).
\] Hence the results.

\vs
\noindent
\textbf{(4.3)}
\textit{Let $H_r$ give $F_x,F_y$ by projection. If $g_1(I)$ is b.v.~with respect to $G_x$ and if $g_2(I)$ is b.v.~with respect to $G_y$, then $g(T)$ is b.v.~with respect to $L_r$.} (---by (4.1)).

\section{Integration around a set $E$ of points.}
\label{Integration around a set E of points}
For $g(I)$ defined in $S$ we consider an auxiliary function $g_E(I)$, equal to $g(I)$ when $I.E$ is not null, and equal to 0 otherwise.

Kempisty [12] and [18] gives the theory for a $k$-dimensional Cartesian space and (the equivalent of) $H=H_{\frac 12}$. But by using the general one-dimensional theory developed in Chapter 1 we obtain results which may easily be generalised to a $k$-dimensional theory. 

By an extension of Kempisty's definitions we call $(F)\overline{\int}_R\,g_E(I)$ (written as $(F)\overline{\int}_{R,E}\,g(I)$), and $(F)\underline{\int}_R\,g_E(I)$
(written as $(F)\underline{\int}_{R,E}\,g(I)$) respectively the \textit{upper and lower integrals}, and the integral of $g(I)$ around $R,E$ and with respect to $F$.

Similar definitions can be given for the \textit{upper and lower norm-limits, and the norm-limit, of $g(I)$ around $R,E$ and with respect to $G$}, written respectively as
\[
(N;G)\overline{\int}_{R,E}\,g(I),\;\;\;\;\;
(N;G)\underline{\int}_{R,E}\,g(I),\;\;\;\;\;
(N;G){\int}_{R,E}\,g(I).
\]
When $g(I)$ is additive and absolutely continuous, Lebesgue [16] 159, defines a set function $g(E)$ by using certain outer covers of $E$, and Kempisty [12] 19, has noted that

\vs
\noindent
\textbf{(5.01)}
\textit{when also $E$ is closed, $g(E) =(F){\int}_{I,E}\,g(I)$ for $I \supset E$.}

\vs

If $g(I)$ is absolutely continuous and $(F_1)\int_R\,g(I)$ exists, and if $g_1(E)$ is the set function defined by Burkill [6] theorem 4.2, then

\vs
\noindent
\textbf{(5.02)}
\textit{for closed $E$, $g_1(E) = (F){\int}_{I,E}\,g(I)$ for $I\supset E$.}

\vs
\noindent
Similar results hold in $k$-dimensions.

Let $I$ be an interval of $\bar S$ containing $E$. Then $g(I)$ is \textit{integrable around} $E$ when for each $\ve>0$ there is an $I$ so that osc$(g_E;I;F)<\ve$. Similarly for the norm-limit. In general there is no smallest interval  of $\bar S$ containing $E$. Denote by $I_E'$ the product set of all intervals of $\bar S$ containing $E$. The \textit{integral around} $E$ is then the upper limit of $(F)\overline{\int}_{I,E}\,g(I)$, when $g(I)$ is integrable around $E$. (cf.~Kempisty [12] 19.)

\vs
\noindent
\textbf{(5.03)}

\noindent
\textbf{(a)}
\textit{Let $F$ satisfy (vi), (vii). Then if $g(I)$ is integrable around $E$, and $y$ is in $I'_E$, the $b(y;F)$ for $g_E(I)$ is zero.}

\noindent
\textbf{(b)}
\textit{Also osc$(g;R;F)=0$ for $R\subset I'_E$ and $R$ in $\bar S$.}

\noindent
By I(3.07), $(I'_E)\sum b(y;F) \leq \mbox{osc}(g_E;I_1;F)$ for all $I_1 \supset E$, and this can be made arb\-itrarily small. Hence (a). For (b), osc$(g_E;R;F) \leq \mbox{osc}(g_E;R;F) \leq \mbox{osc}(g_E;I_1;F)$.
\textbf{(5.04)}
\[
(F)\underline{\int}_{R,E}\,g(I)
\leq (F)\underline{\int}_{R,E}(F)\underline{\int}_{I,E}\,g(I)
\leq (F)\overline{\int}_{R,E}(F)\overline{\int}_{I,E}\,g(I)
\leq (F)\overline{\int}_{R,E}\,g(I).
\]
(cf.~Kempisty [12] 19, Theorem 1.)
By I(2.11a),
\begin{eqnarray*}
(F)\overline{\int}_{R,E}(F)\overline{\int}_{I,E}\,g(I)
&\leq &(F)\overline{\int}_{R}\left((F)\overline{\int}_{I}\,g_E(I)\right)_E \vt
&=& (F)\overline{\int}_{R}(F)\overline{\int}_{I}\,g_E(I)\vt
&\leq &(F)\overline{\int}_{R}\,g_E(I)
\;\;\;= \;\;\;(F)\overline{\int}_{R,E}\,g(I).
\end{eqnarray*}
Similarly for the first inequality of the result.

Kempisty puts $g^E(I) = g(I) - g_E(I)$ and then states that
\[(F){\int}_{R,E}(F){\int}_{I}\,g^E(I) =0\]
when $g(I)$ is additive. This is not true in general. For example, let $E=(-1,0)$ and $g(I)-S(f;I)$ where $f(x)=0$ for $x\leq 0$ and $f(x) =1$ for $x>0$. Then 
\[
\begin{array}{l}
g(I)=0\;\;\mbox{ unless }\;\;I=a\mbox{---}b,\;\;\; a \leq 0<b,\;\;\mbox{ when }\;\;g(I)=1, \vt
g^E(I)=0\;\;\mbox{ unless }\;\;I=a\mbox{---}b,\;\;\; a < 0<b,\;\;\mbox{ when }\;\;g_E(I)=1.
\end{array}
\] 
Hence
$g^E(I)=0\;\;\mbox{ unless }\;\;I=0\mbox{---}b,\;\;\mbox{ when }\;\;g^E(I)=1$. Hence $(F)\overline{\int}_I g^E(I)=1$, $(F)\underline{\int}_I g^E(I)=0$ when $I=a\mbox{---}b$, $a<0<b$. However, we have

\noindent
\textbf{(5.05)}
\[
-\mbox{osc}(g_E;R;F) \leq 
(F)\underline{\int}_{R,E}(F)\underline{\int}_{I}\,g^E(I)
\leq 
(F)\overline{\int}_{R,E}(F)\overline{\int}_{I}\,g^E(I)
\leq \mbox{osc}(g_E;R;F)
\]
\textit{whenever $g(I)$ is a Stieltjes difference, i.e.~additive and independent of the bracket conventions at the ends of $I$.}

\noindent
For $(F)\overline{\int}_I g^E(I)=g(I) -(F)\underline{\int}_I g^E(I)$ so that
\[
(F)\overline{\int}_{R,E}(F)\overline{\int}_I g^E(I)
\leq (F)\overline{\int}_{R,E} g(I) -
 (F)\underline{\int}_{R,E}(F)\underline{\int}_I g^E(I) 
\leq \mbox{osc}(g_E;R;F)
\]
by (5.04).

\vs
\noindent
\textbf{Corollary.}
\textit{If $(F)\overline{\int}_{R,E} \,g(I)$ exists and $g(I)$ is a Stieltjes difference then 
\[
 (F){\int}_{R,E}(F){\int}_I g^E(I) =0,
\]
i.e.~each Stieltjes difference may be replaced in this case by two functions, one zero on every interval not containing points of $E$, and the other having its integral zero around $E$.
} (Kempisty [12] 19, theorem 2.)

\vs
\noindent
\textbf{(5.06)}
\textit{If $(F)\int_I \,g(I)$ exists for $I \supset E$, then
\[
(F)\overline{\int}_{R,E} \,g(I) = (F)\overline{\int}_{R,E} \,F(g;I),\;\;\;\;\;\;
(F)\underline{\int}_{R,E} \,g(I) = (F)\underline{\int}_{R,E} \,F(g;I),
\]
where for $R\subset I$ and $R$ in $\bar S$, $F(g;R) =(F)\int_R \,g(I)$.} For proof use I(3.10).

\vs
\noindent
\textbf{Corollary.}
\textit{In (5.05) we may replace ``$g(I)=S(f;I)$'' by ``$F(g;I)$ exists''.}

\vs
\noindent
\textbf{(5.07)}
\textit{Let $I_\ve$ in $\bar S$ be such that $I_\ve \supset E$ and osc$(g_E;I_\ve;F) < \ve$ for given $\ve>0$. If for all $\ve>0$ and all $R$ in $\bar S$ with $R\subset I_\ve$ we have $\left|(F)\overline{\int}_{R,E} \,g(I)\right| < \ve$, then (a) the integral of $|g(I)|$ around $E$, and around $E_1\subset E$, is zero.} (Kempisty [12] 19, theorem 3, for $F-F_1$.) \textit{(b) The integral of $g(I)$ around $E_1$ is zero also.} (Kempisty [12] 20, theorem 4, for $F=F_1$.) 

By I(3.09),

\vs
\noindent
\textbf{(5.08)}
\textit{If $F(g;I)$ exists and is integrable around $E$ then for $E_1\subset E$ and $R\subset I'_\ve$,}
\[
(F){\int}_{R,E_1}(F){\int}_{I} \,g^E(I)=0,\;\;\;\;\;\;
(F){\int}_{R,E_1}(F){\int}_{I,E} \,g(I)  =
(F){\int}_{R,E_1}\, g(I) .
\] 
(Kempisty [12] 20, theorem 5, for $F=F_1$.) From (5.06 Corollary), and (5.07b).

\vs
\noindent
\textbf{(5.09)}
\textit{If $F(g;I)$ exists and is integrable around $E$ then if $R\subset E_1$,}
\[
(F)\overline{\int}_{R}\left|(F)\int_{I,E} \,g(I)\right|
= (F)\overline{\int}_{R,E}\left| \,g(I)\right|
\]
(Kempisty [12] 20, theorem 6, for $F=F_1$.) For (1)
\begin{eqnarray*}
(F)\overline{\int}_{R}\left|(F)\int_{I,E} \,g(I)\right|
&=&
(F)\overline{\int}_{R}\left|(F)\int_{I} \,g_E(I)\right|,\vt& \leq &
(F)\overline{\int}_{R}\left| \,g_E(I)\right|
\;\;\;=\;\;\; (F)\overline{\int}_{R,E}\left| \,g(I)\right|
\end{eqnarray*}
by (1.04) for Burkill integrals. But we have
\begin{eqnarray*}
(F)\int_{I,E}\, g(I) + (F)\int_{I}\, g^E(I)
&=& (F)\int_{I,E}\, \left(g(I)+g^E(I)\right) \;\;\;=\;\;\; F(g;I), \vt
|F(g;I)| &\leq & 
\left|(F){\int}_{I,E}\, g(I)\right|
+ \left|(F){\int}_{I}\, g^E(I)\right|.
\end{eqnarray*}
But by (5.06 Corollary), (5.05), and (5.07a); and then by I(3.10),
\begin{eqnarray*}
(F)\overline{\int}_{R,E}\left|(F)\int_{I} \,g^E(I)\right|
&=& 0, \vt
(F)\overline{\int}_{R,E}\left|g(I)\right|
&=& (F)\overline{\int}_{R,E}\left|F(g;I)\right| \vt
&\leq &
(F)\overline{\int}_{R,E}\left|(F)\int_{I,E} \,g(I)\right| \;\;+\vt
&&\;\;\;\;\;\;\;\;+\;\;
(F)\overline{\int}_{R,E}\left|(F)\int_{I} \,g^E(I)\right| \vt
&=& (F)\overline{\int}_{R}\left|(F)\int_{I,E} \,g(I)\right|
\end{eqnarray*}
From (1) we have the result.

\vs
\noindent
\textbf{Corollary.}
\textit{When $N(I)$ exists and its norm-limit around $E$ exists,}
\[
\var\left(\left((N;G) \int_{I,E}\,g(I)\right);R;G\right)
= \var\left(g_E(I);R;G\right)\;\mbox{ for }\;R \subset I'_E.
\]

\vs
\noindent
\textbf{(5.10)}
\textit{In order that $g(I)$, whose $N(I)$ exists, has a norm-limit around a set $E$ it is necessary and sufficient that given $\ve>0$ there is $R$ in $\bar S$ enclosing $E$ and a $\delta - \delta(\ve)>0$, such that $\left|(R;D)\sum g^E(I)\right| < \ve$ for every division $D$ in $G$ over $R$ with norm$(D)<\delta$
} (Kempisty [12] 20, theorem 7, for $G=G_1$.)

\vs
\noindent
\textbf{Sufficiency.}
$\left|(R;D)\sum g(I) - (R;D)\sum g_E(I)\right| < \ve$ for each $D$ in $G$ over $R$ with norm$(D)<\delta$. Hence letting $\delta \rightarrow 0$,
\[
\left|N(R) -(N;G) \overline{\int}_{R,E}\,g(I)\right| \leq \ve\;\;\;\mbox{ and }\;\;\;
\left|N(R) -(N;G) \underline{\int}_{R,E}\,g(I)\right| \leq \ve.
\]
Hence osc$N(g_E;R;G) \leq 2\ve$. Adding intervals (necessarily in $\bar S$) to fill in the gaps in $R$, without altering the two extreme end-points, we obtain an interval $I$ in $\bar S$ such that $I \supset I'_E$. Hence $g(I)$ has a norm-limit around $E$. (We have assumed that $W$ is in $\bar S$.)

\vs
\noindent
\textbf{Necessity.}
\begin{enumerate}
\item[(1)]
Let $I \supset I'_E$ with osc$N(g_E;I;G)< \ve/5$ and take a division $D$ in $G$ of $I$ such that
\item[(2)]
\begin{eqnarray*}
(I;D) \sum g_E(I) &<& (N;G) \overline{\int}_{I,E}\,g(I)+ \frac \ve 4\;\;\mbox{ and }\vt
(I;D) \sum g_E(I)& > &(N;G) \underline{\int}_{I,E}\,g(I)- \frac \ve 4.
\end{eqnarray*}
Let $R$ be the part of $D$ made up of all intervals containing points of $E$. Then
\item[(3)]
$(I;D) \sum g_E(I) = (R) \sum g_E(I)$. Now let $\delta = \delta(\ve)$ be such that
\item[(4)]
\[
(N;G) \underline{\int}_{I,E}\,g(I)-\frac \ve 4 <
(R;D') \sum g_E(I) <
 \overline{\int}_{I,E}\,g(I)+ \frac \ve 4
 \]
for every subdivision $D'$ of $R$ and in $G$, with norm$(D') < \delta$. 
\end{enumerate}
Since
\[
(R;D') \sum g^E(I) = (R;D') \sum g(I) - (R;D') \sum g_E(I)
\]
we have
\begin{eqnarray*}
\left|(R;D') \sum g^E(I)\right|& =& \left|(R;D') \sum g(I) - (R;D') \sum g_E(I)\right|\vt
&<& \frac \ve 4 + \left|(R) \sum g(I) - (R;D') \sum g_E(I)\right|
\end{eqnarray*}
by I(3.09) for norm$(D)< \delta_1 = \delta_1(\ve)$, since osc$N(g_E;I;G)<\ve/4$; and then by (3) the latter
\[
\;\;\;\;=\;\;\;\;\frac \ve 4 + \left|(I;D) \sum g_E(I) - (R;D') \sum g_E(I)\right|.
\]
Thus from (1), (2), and (4), we now have 
$\left|(R;D') \sum g^E(I)\right| <\ve$ as required.

Kempisty [12] \textsection\textsection 5,6 gives the following definitions (for $S=S_1$). 
\begin{itemize}
\item
\textit{$g(I)$ is absolutely continuous around $E$} when $g_E(I)$ is absolutely continuous in $I'_E$ (i.e.~in every interval $I$ of $\bar S$ which is in $I'_E$).  
\item
\textit{$g(I)$ is of bounded variation (b.v.) around $E$ when} for all intervals $I$ in $\bar S$ with $I \subset I'_E$, $\var(g_E;I;G)<\infty$.
\end{itemize}
\textbf{(5.11)}

\noindent
\textbf{(a)} 
\textit{If $g_1$ and $g_2$ are absolutely continuous (or b.v.) around $E$, so is $g_1+g_2$.}

\noindent
\textbf{(b)}
\textit{If $g$ is absolutely continuous (b.v.) around $E$, it is absolutely continuous (b.v.) around around $e \subset E$.} 

\noindent
(Kempisty [12] \textsection 5, theorems 1,2; \textsection 6, theorems 1,2.)

\vs
\noindent
\textbf{(5.12)}
\textit{If $g(I)$ has a norm-limit around $E$, the necessary and sufficient condition that $g(I)$ is absolutely continuous around $E$, is that its norm-limit around $E$ is absolutely continuous in $I_E'$.}

\noindent
(Kempisty [12 \textsection 5, theorem 5.) Use I(3.09).

\vs
\noindent
\textbf{(5.13)}
\textit{$(N;G)\overline{\int}_R\,g^E(I)$ and $(N;G)\underline{\int}_R\,g^E(I)$ are absolutely continuous around $E$, when $N(I)$ exists and $g(I)$ is integrable around $E$.} 

\noindent
Use (5.06 Corollary) and (5.05), and then I(3.10).

\vs
\noindent
\textbf{(5.14)}
\textit{If $g(I)$ is absolutely continuous around $E$ it is b.v.~around $E$.} 

\noindent
(Kempisty [12] \textsection 6, theorem 3.)

\noindent
$g_E(I)$ is absolutely continuous in $I_E'$, and so is b.v.~in $I_E'$ by (2.1). Hence $g(I)$ is b.v.~around $E$.

\vs
\noindent
\textbf{(5.15)}
\textit{If $g(I)$ is continuous, non-negative, and decreasing on sub-division, then $g(I)$ has a norm-limit around $E$.}

\noindent
(Kempisty [12] \textsection 6, theorem 4.) For $g_E(I)$ is then continuous, non-negative, and decreasing on sub-division, and by (3.2) we have the result.

\vs
In (5.13) we must in general suppose that $N(I)$ exists and is integrable around $E$. For take $F=F_1$, $W=[0,1/2]$, $g(I) = S(f;I)$ where $f(0)=-1=f\left(\frac 12\right)$ and
\[
f(x) =-n\;\;\mbox{ for }\;\;\frac 1{2(n+1)}\leq x<\frac 1{2n},
\]
and $E = \sum_{n=1}^\infty E_n$ where $E_n = \left( \frac 1{2n+2}, \frac 1{2n+1}\right)$. Then as before (5.05), $g^E(I)=0$  unless for some $n>1$, $I=a$---$\frac 1{2n}$ where $a \geq \frac 1{2n+1}$, when $g_E(I)=1$. And by I(2.11),
\[
(N;G) \overline{\int}_{\left(0,\frac 13\right)}\,g^E(I) \geq 1+1+\cdots =+\infty.
\]
In reality, when $S$ is not $S_1$, the absolutely continuous and b.v.~functions around $E$ are usually generalised absolutely continuous and generalised b.v.~functions around $E$, respectively, in some sense. Further investigation in this direction would bring in the derivatives of $g(I)$, which are outside the scope of this thesis.

\section{Density integration.}\label{Density integration}
In this section we consider an interval function (depending on $g(I)$) whose integral over the fundamental interval $W$ gives the completely additive extension of $g(I)$ over sets $E$ when $g(I)$ is additive and absolutely continuous. But the interval function is interesting in itself, and provides an exercise for the methods which have been built up in the previous sections. In appendix 1 is given an application of this interval function to Hilbert space.

Let $E$ be a (Lebesgue) measurable set in the fundamental interval $W$, and let $g(I)$ be an interval function. We put
\[
E(I) = K(g;E;I) = \frac{g(I)m(EI)}{mI}
\]
where $mE$ is the Lebesgue measure of $E$, and consider integration of $K(I)$ over $W$. The \textit{upper density integral $(D;G) \overline{\int}_E\,g(I)$ of $g(I)$ over $E$} is defined to be
$(N;G) \overline{\int}_W\,K(I)$, and the \textit{lower density integral 
$(D;G) \underline{\int}_E\,g(I)$ of $g(I)$ over $E$} is defined to be
$(N;G) \underline{\int}_W\,K(I)$.

\vs
\noindent
\textbf{(6.01)}
\textit{If $g(I)$ is additive and absolutely continuous, then
\[
(D;G) \overline{\int}_E\,g(I)=(D;G) \underline{\int}_E\,g(I)=g(E),
\]
where $g(E)$ is the Lebesgue integral of $g(I)$ over $E$.
}

\noindent
Since $g(I)$ is additive and absolutely continuous, the derivative of $g(I)$ at a point $x$, say $g'(x)$, exists almost everywhere in $W$, and is integrable over $W$. And $g(E) = \int_E\,g'(x)\,dx$ for every measurable set $E$.

Let $D'$ in $G$ be a division of $W$. Then 
\[
(W;D') \sum K(I) =(W;D') \sum g(I) \frac{m(EI)}{mI}
= \int_W (W;D') \sum c(I;x) \frac{m(EI)}{mI}g'(x) dx,
\]
where $c(I;x)$ is the characteristic function of $I$.
\begin{enumerate}
\item[(1)]
If $x$ is inside $J$, an interval of $D'$,
\[
\left|(W;D') \sum c(I;x) \frac{m(EI)}{mI}g'(x) \right|
=\left| \frac{m(EJ)}{mJ}g'(x)\right| \leq |g'(x)|.
\] 
\item[(2)]
The division points of a sequence of $D'$ are enumerable and so of measure zero, so that as norm$(D') \rightarrow 0$, 
$(W;D') \sum c(I;x) \frac{m(EI)}{mI}$ tends to 1 almost everywhere in $E$, and to 0 almost everywhere in $CE$, by Lebesgue's density theorem. Further,
\item[(3)]
$\int_W\,|g'(x)| dx < +\infty$.
\end{enumerate}
From (1), (2), (3), and Legesgue's convergence theorem, 
\[
(W;D')\sum K(I) \rightarrow \int_E\,g'(x) dx =g(E)
\]
when norm$(D') \rightarrow 0$. This is the result required.

Note that $F$ need only contain one Riemann succession, and that over $W$, for the result to be true.

Once the set function $mE$ has been obtained, this evaluation of $g(E)$ is particularly useful when
\begin{enumerate}
\item[(a)]
$E$ is fixed and we have many $g(I)$.
\item[(b)]
$g(I)$ is fixed and we have many sets $E$.
\end{enumerate}
\textbf{(6.02)}
\begin{enumerate}
\item[\textbf{(a)}]
\textit{If $E$ is of measure zero then}
\[
(D;G)\overline{\int}_E\,g(I)=0=(D;G)\underline{\int}_E\,g(I).
\]
\item[(b)]]
\textit{For $R$ an $I^\sigma$ in $\bar S$,}
\[
(D;G)\underline{\int}_R\,g(I) \leq (N;G)\underline{\int}_R\,g(I) \leq(N;G)\overline{\int}_R\,g(I) \leq (D;G)\overline{\int}_R\,g(I).
\]
\item[(c)]
If also, $g(I)$ is continuous, then
\[
(D;G)\underline{\int}_R\,g(I) = (N;G)\underline{\int}_R\,g(I),\;\;\;\;\;(D;G)\overline{\int}_R\,g(I) = (N;G)\overline{\int}_R\,g(I).
\]
\end{enumerate}
By (vii), $W-R'$ is in $\bar S$, so that by (iii) we have (b). For (c), let $R$ consist of $m$ intervals. Then in any division of $W$ there are at most $2m$ intervals $I$ which each contain points of $FR$. For these intervals, $g(I)\rightarrow 0$ by hypothesis.

\vs
\noindent
\textbf{(6.03)}
\textit{In order that $(D;G)\underline{\int}_R\,g(I)$ and $(D;G)\overline{\int}_R\,g(I)$ should be finite for every measurable set in $W$ it is necessary and sufficient that $g(I)$ should be b.v.
}

\noindent
\textbf{Sufficiency.}
\[
\left|(W;D') \sum g(I) \frac{m(EI)}{mI}\right|
\leq (W;D') \sum| g(I) |,
\]
whence the result.

\noindent
\textbf{Necessity.}
We first prove three Lemmas, supposing that (vi), (vii), (viii) hold.

\vs
\noindent
\textbf{Lemma 1}
\textit{$g(I)$ is bounded when a fixed point $x$ is in $I'$ and $mI \rightarrow 0$.}

\noindent
Suppose that $|g(I)|\rightarrow +\infty$, and let $E$ have been defined as an $I^\sigma$ in $W-J$ of $\bar S$, where $x$ is in $J^o$. Let the rest of $E$ lie in $J_1$ of $S$, which is a positive distance from $W-J$, and such that $x$ is in $J_1^o$. Them for norm$(D')$ small enough,
\[
\left|(W;D') \sum g(I) \frac{m(EI)}{mI}\right| \geq\left| g(J_1) \frac{m(EJ_1)}{mJ_1}\right| -
\left|(W-J;D') \sum K(I)\right|,
\]
where $D'$ includes $J_1$. We take $D'$ in $W-J$ so that the last term is arbitrarily near to $(N;G)\overline{\int}_{W-J}\,K(I)$, which is finite by I92.12b) and by hypothesis. By suitable choice of $J_1$, and then an interval $I_1$ in $J_1$ for $E$, we have
\[
\left| g(J_1) \frac{m(EJ_1)}{mJ_1}\right| -
\left|(W-J;D') \sum K(I)\right| >
\frac 13\left|g(J_1)\right| - \left|(N;G)\overline{\int}_{W-J}\,K(I)\right| -1>n.
\]
We can obviously take $I_1$ so that $x$ is not in $I_1'$, enabling us to continue the construction of $E$. Then by induction we find an $E$ such that one of $(D;G)\overline{\int}_E\,g(I)$, $(D;G)\underline{\int}_E\,g(I)$ is not finite, contrary to hypothesis.

Similarly if $x$ is a permanent point of division of all $\{D_n\}$ in $F^*$ (the closure of $F$). We consider each side of $x$ separately. Hence $|g(I)|$ is bounded.

\vs
\noindent
\textbf{Lemma 2.}
\textit{For a sequence $J_1, \ldots , J_n, \ldots$ of intervals of $\bar S$ let the length of an interval $I_n$ of $S$ tend to zero as $n\rightarrow \infty$, where $J_1, \ldots , J_n$ are outside $I_n$ and $J_{n+1}, J_{n+2}, \ldots$ are inside $I_n$. Then
\[
\sum_{n=1}^\infty \left|(N;G)\overline{\int}_{J_n}\,g(I)\right|\;\;\;\mbox{ and }\;\;\;
\sum_{n=1}^\infty \left|(N;G)\underline{\int}_{J_n}\,g(I)\right|
\]
are convergent.}

\noindent
Let $E$ be $\sum_{n=1}^\infty J_n$. There is a division $D_n$ in  $G$ over $W$with division -points the ends of $J_1, \ldots , J_n$ (by (vii) and (iii)) such that if $S_i$ is the set of intervals over $J_i$ ($i=1, \ldots ,n$), with norm$(S_i)<1/n$,
\[
(S_i) \sum g(I) >(N;G) \overline{\int}_{J_i}\, g(I) - \frac 1{2^{n+1}}.
\]
By (viii) we can also suppose that $I_n$ is in $D_n$. Then
\begin{eqnarray*}
(W;D_n) \sum K(I) &=&
\sum_{i=1}^n (S_i) \sum g(I) + g(I_n)\frac{m(EI_n)}{mI_n} \vt
&>& 
\sum_{i=1}^n  (N;G) \overline{\int}_{J_i}\, g(I) + g(I_n)\frac{m(EI_n)}{mI_n}.
\end{eqnarray*}
The last term is bounded by Lemma 1 since $mI_n \rightarrow 0$ by hypothesis. Hence
\[
\sum_{i=1}^n  (N;G) \overline{\int}_{J_i}\, g(I) \leq M
\]
independent of $n$, for some $M<+\infty$. Similarly
$\sum_{i=1}^n  (N;G) \underline{\int}_{J_i}\, g(I))\geq M'$
independent of $n$, for some $M'>-\infty$. Take $|M'|\leq M$. Then we have
\[
\left|\sum_{i=1}^n  (N;G) \overline{\int}_{J_i}\, g(I) \right| \leq M.
\]
A subsequence of $J_1, \ldots ,J_n, \ldots $ satisfies the same conditions as the main series, so that taking the positive terms in $\sum_{i=1}^n  (N;G) \overline{\int}_{J_i}\, g(I)$, and then the negative terms, we see that
$\sum_{i=1}^n  \left|(N;G) \overline{\int}_{J_i}\, g(I))\right| \leq 2M$, and similarly
\[
\sum_{i=1}^n  \left|(N;G) \underline{\int}_{J_i}\, g(I))\right| \leq 2M .
\]
\textbf{Lemma 3.}
\textit{$\sum_{i=1}^\infty |g(J_n)|$is convergent when the $J_n$ are also in $S$. }

\noindent
By hypothesis,
\[
(N;G) \overline{\int}_W\, g(I) = (D;G) \overline{\int}_W\, g(I) \;\;\;\mbox{ and }\;\;\;
(N;G) \underline{\int}_W\, g(I) = (D;G) \underline{\int}_W\, g(I)
\]
are finite, so that there is an at most enumerable number of singularities $\{y_i'\}$ and $\sum_{i=1}^\infty c(y'_I;G)$ is finite (see I(3.07)). Then by I(3.09),
\[
\sum_{n=N}^A (N;G) \underline{\int}_{J_n}\, g(I) - \sum_{i=1}^\infty c(y_i';G) - \ve \;\;<\;\;
\sum_{n=N}^A (N;G) \overline{\int}_{J_n}\, g(I) + \sum_{i=1}^\infty c(y_i';G) + \ve
\]
for $A\geq N >N_0=N_0(\ve)$. Hence by Lemma 2, for some $M_1$ ($<\infty$) independent of $n$, 
$\left| \sum_{i=1}^n g(J_i)\right| \leq M_1$. Then as in Lemma 2 we have the result.

We are now in a position to prove the main result. Suppose that $\var(g;W;G)=+\infty$ and let $D'$ in $G$ over $W$ contain at least three intervals. Let $W_1$ be the interval covered by the first two, and let $W_2$ be the interval covered by all the intervals of $D'$ except the first. Then
\[
\var(g;W;G) \leq \var(g;W_1;G) + \var(g;W_2;G)
\]
as may easily be seen when we use (viii) and (vii). (Each interval of $D'$ can occur in $D_1$ in $G$ over $W_1$ together with $D_2$ in $G$ over $W_2$.) Hence $\var(g;W_i;G) = +\infty$ for $i=1$ or 2 (or both).

Take $i$ the first such. Repeating this construction we obtain a sequence $I_1, I_2, \ldots $ of intervals such that $I_1 \supset I_2 \supset \cdots$, $mI_n \rightarrow 0$ as $n \rightarrow \infty$, and $\var(g;I_n;G) =+\infty$ for $n=1,2, \ldots$. There is a fixed point, say $x$, in $I_n'$ of $\bar S$.

By Lemma 1, $|g(I)$ is bounded, say by $M_2$,  if $x$ is in $I'$ and $mI \rightarrow 0$ with $I$ in $S$. Then since $\var(g;I_n;G)=+\infty$, there is a division $D_n$ of $I_n$ and in $G$, with sum greater than $2M_2+2$. Omitting the interval (or two intervals) $K_n$ in $D_n$ which has $x$ in $K_n'$, we have $(D_n-K_n)\sum |g(I)| >2$. Taking the odd or even terms we have $\sum |g(I)| >1$.

There is an $I_m'$ in $K_n^o$ ($m>n$), so that we may repeat the construction, and so obtain a series
$
\sum |g(J_n)| \geq 1+1+\cdots = +\infty
$
of the type considered in Lemma 3, which is therefore contradicted. Hence $\var(g;W;G) < \infty$, which was to be proved.

\vs
\noindent
\textbf{Corollary.}
\textit{If $E$ is restricted to be open, it is necessary that $g(I)$ should be b.v.}

\vs
\noindent
\textbf{(6.04)}
\[
\left.
\begin{array}{r}
\left|(D;G) \overline{\int}_E\, g(I)\right| \vt
\left|(D;G) \underline{\int}_E\, g(I)\right|
\end{array}
\right\}
 \leq 
(D;G) \overline{\int}_E\, |g(I)| \leq \var(g;W;G).
\]
For $|\sum K(I)| \leq \sum |K(I)| \leq \sum | g(I)| $.

\vs
\noindent
\textbf{(6.05)}
\textit{If $E_1 \subset E_2$ then }
\[
(D;G) \overline{\int}_{E_1}\, |g(I)| \leq
(D;G) \overline{\int}_{E_2}\, |g(I)|, \;\;\;\;\;\;
(D;G) \underline{\int}_{E_1}\, |g(I)|\leq
(D;G) \underline{\int}_{E_2}\, |g(I)|.
\]
For $\sum |g(I)| \frac{m(E_1I)}{mI} \leq \sum |g(I)| \frac{m(E_2I)}{mI}$.

\vs
\noindent
\textbf{(6.06)}
\textit{
Let $g(I)$ have a norm-limit $N(R)$ over every $R$ of $\bar S$. Then}
\[
(D;G) \overline{\int}_{E}\, g(I) =
(D;G) \overline{\int}_{E}\, N(I), \;\;\;\;\;\;
(D;G) \underline{\int}_{E}\, g(I)=
(D;G) \underline{\int}_{E}\, N(I).
\]
By I(3.10), we have for norm$(D')<\delta = \delta(\ve)$,
\[
\left|(W;D')\sum g(I) \frac{m(EI)}{mI} -
(W;D')\sum N(I) \frac{m(EI)}{mI}\right| \leq
(W;D')\sum |g(I)- N(I)| <\ve.
\]
Hence the two sums give the same limit-points, and in particular we have the results.

\vs
\noindent
\textbf{(6.07)}
\textit{In order that $(D;G_1) \overline{\int}_{E}\, g(I)$ and $(D;G_1) \underline{\int}_{E}\, g(I)$ should be finite and equal for every (measurable) set $E$ in $W$ it is necessary and sufficient that the norm-limit of $g(I)$ should exist and be absolutely continuous.}

\vs
\noindent
\textbf{Sufficiency.}
By (6.06) and (6.01).

\noindent
\textbf{Necessity.}
By hypothesis $N(W)$ exists (take $E=W$) so that $N(I)$ exists for every $I$ in $\bar S$, and by (6.06) we can replace $g(I)$ by $N(I)$. Then by a simplified form of proof as in (6.03), $N(I)$ must be b.v.

By (1.13), $N(I)$ tends to a limit, say $N(v-)$, if $I=x$---$v$ in $\bar S$ and $x \rightarrow v$, and similarly for $N(v+)$. Then by (1.17)
\begin{eqnarray*}
&&(k';G)\underline{\int}_{W}\, N(I)\frac{m(EI)}{mI} \vt
&=&(k';G)\underline{\int}_{W}\, N_2(I)\frac{m(EI)}{mI}
+(W'_-)\sum  N(v-) +
(W'_+)\sum \underline{m}(v+)N(v+)
\end{eqnarray*}
where $\overline{m}(v-), \ldots$ are the $\overline{h}(v-), \ldots$
for $h(I) = \frac{m(EI)}{mI}$, and $N_2(I)$ is the $g_2(I)$ for $N(I)$, and $G=G_1$.

But the norm-limit of $N(I)\frac{m(EI)}{mI}$ exists by hypothesis, so that the $k'$-limits on the left are equal. Hence
\[
(W'_-)\sum \overline{m}(v-) N(v-) =
(W'_-)\sum \underline{m}(v-) N(v-).
\]
But we can easily find an $E$ for each $v$ in $W_-'$ such that $\overline{m}(v-)>\underline{m}(v-)$, and $E$ can be open. Hence $N(v-)=0$, and similarly $N(v+)=0$. Hence $N(I) = N_2(I)$, i.e., $N(I)$ is the Stieltjes difference of a function $f(x)$ which is continuous and b.v. We put $f(x)=g(x) +h(x)$ where $g(x)$ is absolutely continuous and $h(x)$ is the singular function.

Using (6.01), where $g(x)$ produces $g(E)$, we have
\begin{eqnarray*}
(D;G) \overline{\int}_{E}\, N(I)|&=&
g(E) +(D;G) \overline{\int}_{E}\, S(h;I)\vt
& =&
(D;G) \underline{\int}_{E}\, N(I) \vt
&=&
g(E) +(D;G) \underline{\int}_{E}\, S(h;I)\;\;\;\mbox{ so that}\vt
(D;G) \overline{\int}_{E}\, S(h;I)|&=&
(D;G) \underline{\int}_{E}\, S(h;I).
\end{eqnarray*}
Since $h(x)$ is singular we can find a set $\mathcal E$ with $m{\mathcal E}=0$, and containing all the variations of $S(h;I)$. We can then find a sequence $G_1 \supset \cdots \supset G_n \supset \cdots $ of open sets such that $\mathcal{E} \subset G_n$ for $n=1,2, \ldots$, and, for each component interval $I$ of $G_n$,
\[
m(I.G_{n+1}) < \frac{mI}{n}.
\]
We construct a set $E$ as follows. Let $H_n$ be the sum of those component intervals $I$ of $G_n$ with $S(h;I)\geq 0$, and then put
\[
E=\sum_{n=1}^\infty \left(H_{2n-1}-G_{2n}\right).
\]
\begin{enumerate}
\item[(1)]
Let $I$ be a component interval of $G_{2n}$. Then
\begin{eqnarray*}
I.E &=& \sum_{i=n+1}^\infty \left(H_{2i-1}-G_{2i}\right)I\;\;\;\mbox{ so that}\vt
\frac{m(EI)}{mI} &=& \frac{\sum_{i=n+1}^\infty \left(H_{2i-1}-G_{2i}\right).I}{mI} \vt
&\leq & 
\frac{\sum_{i=n+1}^\infty m(G_{2i-1}.I)}{mI} \vt
&&\mbox{and by the construction of $\{G_n\}$ this is} \vt
&<& \sum_{i=n+1}^\infty \frac 1{(2i-2)\ldots (2n)} \vt
\mbox{i.e. }\;\;\frac{m(EI)}{mI} &<& \frac e{2n}.
\end{eqnarray*}
\item[(2)]
Let $I$ be a component interval of $H_{2n-1}$. Then
\begin{eqnarray*}
I.E &=& \left(I-G_{2n}.I\right) +\sum_{i=n+1}^\infty \left(H_{2i-1}-G_{2i}\right)I ,
 \vt
\frac{m(EI)}{mI} &=&
\left(1 -\frac{m(G_{2n}.I)}{mI}\right) +
\frac{\sum_{i=n+1}^\infty m\left(\left(H_{2i-1}-G_{2i}\right).I\right)}{mI}
 \vt
&<& 1-\frac 1{2n-1}.
\end{eqnarray*}
\item[(3)]
Let $I$ be a component interval of $G_{2n-1}-H_{2n-1}$. Then
\begin{eqnarray*}
I.E &=& \sum_{i=n+1}^\infty \left(H_{2i-1}-G_{2i}\right).I ,
 \vt
\frac{m(EI)}{mI} &<&  \frac e{2n},
\end{eqnarray*}
as in (1).
\item[(4)]
Now \[
S(h;I)=S(h;I.G_n) = \sum_{i=1}^\infty S(I.I_i^{(n)}), \]
where $I_i^{(n)}$ are the component intervals of $G_n$. Since $h(x)$ is b.v., the series is absolutely convergent. There is therefore a number $i_n$ such that
\item[(5)]
\[
\sum_{i=i_n+1}^\infty\left|S\left(I_i^{(n)}\right)\right| < \frac 1n.
\]
We take $I_1^{(n)}, \ldots ,I_{i_n}^{(n)}$ to be intervals in a division $D_n$ of $W$ and complete $D_n$ with intervals $J_1, \ldots ,J_m$ of norm less than $1/n$. Then $D_n$ is in $G$ since $G=G_1$. Then, by (4) and (5),
\item[(6)]
\[
\left|
(W;D_n)\sum S(I)\frac{m(EI)}{mI} - \sum_{i=1}^{i_n} S\left(I_i^{(n)}\right)
\frac{m\left(EI_i^{(n)}\right)}{mI_i^{(n)}}
\right| < \frac 1n.
\]
Let $n$ be even. Then by (1) and (1.01), we have for $n>n_0(\ve)$,
\[
\left|
\sum_{i=1}^{i_n} g\left(I_i^{(n)}\right)
\frac{m\left(EI_i^{(n)}\right)}{mI_i^{(n)}}
\right| <  (M+\ve)\frac{e}{2n}.
\]
\item[(7)]
Hence from (6), as $n \rightarrow \infty$,
\[
\left|
\sum_{i=1}^{i_{2n-1}} S\left(I_i^{({2n-1})}\right)
\frac{m\left(EI_i^{({2n-1})}\right)}{mI_i^{({2n-1})}}
\right| \rightarrow 0.
\]
\item[(8)]
But by (2) and (3), in an obvious notation,
\begin{eqnarray*}
&&\sum_{i=1}^{i_{2n-1}} S\left(I_i^{({2n-1})}\right)
\frac{m\left(EI_i^{({2n-1})}\right)}{mI_i^{({2n-1})}}\vt
&>&
(H_{2n-1})
\sum_{i=1}^{i_{2n-1}} S\left(I_i^{({2n-1})}\right)
\left( 1 - \frac 1{2n-1}\right)\;\;- \vt
&&\;\;\;\;\;\;-\;\;(G_{2n-1}-H_{2n-1})
\sum_{i=1}^{i_{2n-1}} \left|S\left(I_i^{({2n-1})}\right)\right| \frac e{2n},\;\;\;\mbox{ and by (1.01),}
\end{eqnarray*}
\item[(9)]
\[
(G_{2n-1}-H_{2n-1})
\sum_{i=1}^{i_{2n-1}} \left|S\left(I_i^{({2n-1})}\right)\right| \frac e{2n} <(M+\ve)\frac e{2n}
\]
for $n>n_(\ve)$. 
\item[(10)]
Hence by (7), (8), (9),
$(H_{2n-1}) \sum_{i=1}^{i_{2n-1}} \left|S\left(I_i^{({2n-1})}\right)\right| \rightarrow 0$ as $n \rightarrow \infty$.
\item[(11)]
Similarly, by taking another $E$.
\[
(G_{2n-1}-H_{2n-1}) \sum_{i=1}^{i_{2n-1}} \left|S\left(I_i^{({2n-1})}\right)\right| \rightarrow 0\;\;\;\mbox{  as }\;\;\;n \rightarrow \infty.
\]
\item[(12)]
Now let $D_n'$ consists of $I_1^{(2n-1)},\ldots,I_{2n-1}^{(2n-1)}$ together with $J_1, \ldots ,J_m$ with norm$(D_n')\rightarrow 0$ as $n \rightarrow \infty$, such that each $J_j$ does not cut any of the intervals $I_p^{(2n-1)}$ ($p>i_{2n-1}$). Then by (5), (10), (11), since there is no variation outside $G_{2n-1}$,
\[
(W;D_n')\sum |S(I)| \rightarrow 0 \;\;\;\mbox{ as }\;\;\; n\rightarrow \infty.
\]
\end{enumerate}
But $S(I)$ is increasing on subdivision, with its $c(y';G)=0$. Hence by (3.2) and (12), and then by (1.03), (1.02),
\[
\var(S;W;G) =0,\;\;\;\;\;\;|S(I)| \leq \var(S;W;G)=0.
\]
Hence $N(I)=S(g;I)$, i.e., the norm-limit of $g(I)$ exists and is absolutely continuous, which was to be proved.

\noindent
\textbf{(6.08)}
\textit{Let $E_1, \ldots ,E_n$ be disjoint sets with sum $E$. Then}
\[
\sum_{i=1}^n (D;G)\underline{\int}_{E_i} g(I)
\leq (D;G)\underline{\int}_{E} g(I)
\leq (D;G)\overline{\int}_{E} g(I)
\leq \sum_{i=1}^n (D;G)\overline{\int}_{E_i} g(I).
\]
For $\sum_{i=1}^n (W;D)\sum g(I)
\frac{m\left(EI_i\right)}{mI}
=
(W;D)\sum g(I)
\frac{m\left(EI\right)}{mI}.
$

\vs
\noindent
\textbf{Corollary.}
\textit{If 
\[
(D;G)\underline{\int}_{E_i} g(I) 
=(D;G)\overline{\int}_{E_i} g(I)
\]
($i=1, \ldots ,n$) then $(D;G)\underline{\int}_{E} \,g(I)=(D;G)\overline{\int}_{E} \,g(I)$}

We cannot in general extend these results to infinite sums. For let $g(I)=a>0$ when 0 is in $I'_-$, and otherwise $g(I)=0$ Then
\[(D;G)\underline{\int}_{I} \,g(I)=(D;G)\overline{\int}_{I} \,g(I)=0
\]
when 0 is not in $I'_-$, and otherwise
\[(D;G)\underline{\int}_{I} \,g(I)=0,\;\;\;\;\;\;(D;G)\overline{\int}_{I} \,g(I)=a>0,
\]
and this falsifies the corollary for infinite sums.

In the rest of this section we suppose that $g(I)$ is b.v. We now consider the $k'$-limits of $K(I)$. As for (1.17) we see that the singularities $y'$ of $K(I)$ are included in the singularities $v$ of $g(I)$. Let $g(I)$ satisfy the conditions of (1.17). Then for $h(I) = m(EI)/mI$,

\vs
\noindent
\textbf{(6.09)}
\[
(Dk';G)\overline{\int}_{E} \,g(I)=(Dk';G)\overline{\int}_{E} \,g_2(I)+(W'_-)\sum \bar{h}(v-) g(v-)
+ (W'_+)\sum \bar{h}(v+) g(v+),
\]
\textit{and similarly for the lower $k'$-limit, where we have written 
$(Dk';G)\overline{\int}_{E} \,g(I)$ for 
$(k';G)\overline{\int}_{W} \,K(I)$,
$(Dk';G)\underline{\int}_{E} \,g(I)$ for 
$(k';G)\underline{\int}_{W} \,K(I)$.}

\vs
\noindent
\textbf{Corollary.}
\textit{If 
$(Dk';G)\overline{\int}_{E} \,g(I)=(Dk';G)\underline{\int}_{E} \,g(I)$ then
}
\begin{enumerate}
\item[\textbf{(a)}]
\textit{for every $v$ either $g(v-)=0$ or $\overline{h}(v-)=\underline{h}(v-)$, and similarly for $v+$;}
\item[\textbf{(b)}]
\[(Dk';G)\overline{\int}_{E} \,g_2(I)=(Dk';G)\underline{\int}_{E} \,g_2(I).\]
\end{enumerate}
In connection with (a) it may be pointed out that the $k'$-limits produce additive, nut not necessarily completely additive, interval functions 
\[
g_3(I) = (k';G)
\overline{\int}_{I} \,g(I),\;\;\;\;\;\;
g_4(I) =\underline{\int}_{\,I} \,g(I).
\]
In a sense, when $g(I)$ is b.v., $g_3(I)$ and $g_4(I)$ are completely additive for chains of intervals except in the neighbourhoods of the singularities $v_1, \ldots , v_n, \ldots$, so that except in those neighbourhoods, $g_3(I)$ and $g_4(I)$ can be extended to form set functions $g_3(E), g_4(E)$ respectively by the usual method of outer covers of $E$.

But this method fails ``near'' a $v$, and the extension there is somewhat arbitrary. We have to give to each set $E$ some kind of ``weight'' in the neighbourhood of a $v$, and quite reasonable weights are given by $\overline{h}(v-)$ and $\overline{h}(v+)$ if $\overline{h}(v-)=\underline{h}(v-)$ and $\overline{h}(v+)=\underline{h}(v+)$.
The $Dk'$-integration supplies these weights automatically, though it does not always give useful $g_3(E), g_4(E)$.

\vs
\noindent
\textbf{(6.10)}
\textit{Let the $k'$-limit of $g(I)$ exist as $k'(I)$. Then}
\[
(Dk';G)\underline{\int}_{\,E} \,g(I)=
(Dk';G)\underline{\int}_{\,E} \,k'(I),\;\;\;\;\;\;
(Dk';G)\overline{\int}_{E} \,g(I)=
(Dk';G)\overline{\int}_{E} \,k'(I).
\]
Use I(4.12).

\vs
\noindent
\textbf{(6.11)}
\textit{If $E=R$ in $\bar S$ then}
\[
(Dk';G)\underline{\int}_{\,R} \,g(I)=
(k';G)\underline{\int}_{\,R} \,g(I),\;\;\;\;\;\;
(Dk';G)\overline{\int}_{R} \,g(I)=
(k';G)\overline{\int}_{R} \,g(I).
\]
By I(4.04) we can take $FR$ and $v_1,v_2, \ldots$ as permanent points in the $k'$-limit of $g(I)h(I)$, and then by (vi) we have the result.

\vs
\noindent
\textbf{(6.12)}
\textit{Let the open set $\sum_{i=1}^\infty I_i$ (where $\rho(I_i;I_j)>0$ for $i\neq j$) contain all the variation of $g_2(I)$, i.e.,
\[
\var\left(g_2;\;W-\sum_{i=1}^n I_i;\;G\right) \rightarrow 0\;\;\;\mbox{ as }\;\;\;n \rightarrow \infty.
\]
Then
\[
(Dk';G)\overline{\int}_{E} \,g_2(I) = \sum_{i=1}^\infty
(Dk';G)\overline{\int}_{EI_i} \,g_2(I) ,
\]
and similarly for lower limits.
}

\noindent
By I(4.04) and $\sum_{i=1}^\infty\var\left(g_2;I_i;G\right) <\infty$.

\chapter{Appendices.}

\section{Appendix 1:\\ Function Spaces and Density Integration.}\label{Appendix 1}
Let $K$ be a (restrictedly) additive family of (Lebesgue's) measurable sets $E$ in $W=[0,1]$, and let $g(E)$ be defined for all $E$ in $K$ such that

\vs
\noindent
\textbf{(1)}
$g(E_1)+g(E_2) =g(E_1+E_2)$ \textit{when $E_1,E_2$ are disjoint and in $K$. }

\vs
We then define a functional $P(f)$ to be
\[
\sum_{i=1}^n b_ig(E_i)\;\;\;\mbox{  whenever }\;\;\;
f=f(x) = \sum_{i=1}^n b_i c(E_i;x),
\]
where $b_1, \ldots ,b_n$ are real constants and $E_1, \ldots ,E_n$ are disjoint, non-null, and in $K$.

If $f(x)$ is a characteristic function then each $b_i$ is either 0 or 1, since the $E_1, \ldots , E_n$ are non-null and disjoint. And then by (1) we see that the definition is consistent. Thus by construction,

\vs
\noindent
\textbf{(2)}
\textit{$P(f)$ is a distributive functional of the functions $f$.}

\vs
We denote by $\sigma_2(f)$ the Hilbert function space of functions for which 
\[
||f||^2 = \int_0^1 |f(x)|^2 dx\;\;\;\;\mbox{ is finite,}
\]
and we take $||f||\geq 0$.

\vs
\noindent
\textbf{(3)}
\textit{In order that $P(f)$ should be continuous with respect to $||f||$, for $f$ in $\sigma_2(f)$, it is necessary and sufficient that}
\[
\var\left(\frac{g^2(E)}{mE};W\right) \leq M^2\;\;\;\mbox{ for some }\;\;M>0.
\]
The continuity of $P(f)$, i.e.~$P(f_1) \rightarrow P(f)$ if $||f-f_1|| \rightarrow 0$, is by (2) equivalent, for some $M>0$, to $|P(f)| \leq M||f||$, i.e.,

\vs
\noindent
\textbf{(4)}
\[
\left(\sum_{i=1}^n b_ig(E_i)\right)^2 \leq M^2 \sum_{i=1}^n b_i^2 mE_i \;\;\mbox{ \textit{for all} }b_1, \ldots , b_n, \;\mbox{ \textit{and} }\; n=1,2, \ldots .
\]
For $n=1$ we have $|g(E)|\leq M\sqrt{mE}$, so that $g(E)$ is absolutely continuous. We may obviously suppose that
\[
\left|g(E_1)\ldots g(E_n)\right| >0\;\;\;\mbox{ and so }\;\;\;mE_1\ldots mE_n>0.
\]
From the theory of quadratic forms, (4) is equivalent to 
\[
\det\left(M^2 \,mE_i\, \delta_{ij} -g(E_i)\,g(E_j)\right) \geq 0
\]
where $\delta_{ij}=1$ $(i=j$), $\delta_{ij}=0$ $(i \neq j)$, and $i,j=1,\ldots , n$; and $n=1,2, \ldots$.

To evaluate the determinant we multiply the first row by $g(E_2)$ and subtract from it $g(E_1)$
times the second row. We then multiply the second row by $g(E_3)$ and subtract from it $g(E_2)$ times the third row, and so on; so that we obtain
\[
\frac{M^{2(n-1)}}{g(E_2) \ldots g(E_n)}
\left|
\begin{array}{rrrrrrr}
m_1g_2&-m_2g_1&0&0 &\ldots &0&0\vt
0 & m_2g_3 & -m_3g_2&0&\ldots & 0&0\vt						
\vdots &\vdots &\vdots &\vdots & \ldots &\vdots&\vdots \vt
0&0&0&0&\ldots &m_{n-1}g_n & -m_ng_{n-1} \vt
-g_ng_1&-g_ng_2&-g_ng_3&-g_ng_4&\ldots &-g_ng_{n-1}&M^2m_n-g_n
\end{array}
\right|
\]										
where $g_i = g(E_i)$ and $m_i=mE_i$ for conciseness, and expanding by the bottom row, this equals
\[
M^{2(n-1)} mE_1\ldots mE_n \left(M^2 - \sum_{i=1}^n \frac{g^2(E_i)}{mE_i}\right).
\]
Hence in order that (4) should be true it is necessary and sufficient that for $n=1,2, \ldots$,

\vs
\noindent
\textbf{(5)}
\[
\sum_{i=1}^n \frac{g^2(E_i)}{mE_i} \leq M^2
\]
where $E_1, \ldots ,E_n$ are disjoint and in $K$. Thus, in an obvious notation,
\[
\var\left( \frac{g^2(E)}{mE};W\right) \leq M^2
\]
is necessary; and if $\var\left( \frac{g^2(E)}{mE};W\right) \leq M^2$ then for disjoint sets $E_1, \ldots ,E_n$ in $K$, of measure less than $\delta=\delta(\ve)$,
\[
\sum_{i=1}^n \frac{g^2(E_i)}{mE_i} \leq \left(M+\ve\right)^2
\]
so that $|P(f)| \leq (M+\ve)||f||$ for every $f$ formed from sets $E_1, \ldots ,E_n$ in $K$, or measure less than $\delta$. Since $K$ is additive, this is true for all our special $f$. Hence the results.

\noindent
\textbf{Corollary.}
\textit{From theorems on derivatives, a necessary condition is
\[
\int_0^1 |g'(x)|^2 \leq M^2.
\]
If $\frac{g^2(E)}{mE}$ is absolutely continuous, this integral condition is also sufficient.}

We now construct a complete orthonormal (i.e.~orthogonal and normal) set for this space, by first forming tables of signs.

Inductively, let $a_{11}(1) = +1$, and for $n=1,2, \ldots$, 
\[
\begin{array}{rrll}
a_{i,2j-1}(n+1) &=& a_{i,2j}(n+1) &=\; a_{ij}(n) \; \mbox{for} \; 1\leq i\leq 2^{n-1},\; 1\leq j\leq 2^{n-1};\vt
a_{i,2j-1}(n+1) &=&- a_{i,2j}(n+1) &=\; a_{ij}(n) \; \mbox{for} \; 
i=r+2^{n-1} \mbox{ and }1\leq r,j \leq 2^{n-1},
\end{array}
\]
For fixed $n,i$, the set $a_{i,1}(n), \ldots , a_{i,2^{n-1}}(n)$ is the $i$th row of the $n$th stage.

\vs
\noindent
\textbf{(6)}
\textit{Every row of $2^{n-1}$ terms, of $0$'s except for a single $1$, can be obtained by linear combinations of the rows of the $n$th stage.}

\vs
\noindent
Suppose true for $n$. Then every row of $2^{n-1}$ terms, of $0,0$'s, except for a single pair of the form $+1,+1$ or $+1,-1$ can be got from the rows of the $(n+1)$th stage. Hence the result is true for $n+1$. Being true for 1, it is true generally.

\noindent
\textbf{(7)}
\[
\sum_{k=1}^{2^{n-1}} a_{ik}(n) a_{jk}(n) =0 \;\;\;\;\;(i \neq j).
\]
Suppose true for $n$. Then true for $n+1$, and $1\leq i \leq 2^{n-1}$, $1\leq j\leq 2^{n-1}$; and $2^{n-1}<i\leq 2^{n-1}$, $2^{n-1}<j\leq 2^n$. Hence by induction.

\vs
\noindent
\textbf{(8)}
\[
a_{ij}(n)=a_{ji}(n).
\]
This is true for $n=1$. Suppose true for $1,2, \ldots , n-1$. By construction, 
\[
a_{i,2j-1}(n) = (-1)^r a_{i,2j}(n) \;\;(n>1)\;\;\mbox{ where }\;\;r= 
\left\{
\begin{array}{rll}
0&\mbox{for}& 1\leq i\leq 2^{n-2},\vt
1&\mbox{for}& 2^{n-2}<i\leq 2^{n-1}.
\end{array}
\right.
\]
And by easy induction,
\[
a_{i,2j-1}(n) = (-1)^s a_{2i,j}(n) \;\;(n>1)\;\;\mbox{ where }\;\;s= 
\left\{
\begin{array}{rll}
0&\mbox{for}& 1\leq j\leq 2^{n-2},\vt
1&\mbox{for}& 2^{n-2}<j\leq 2^{n-1}.
\end{array}
\right.
\]
Also, $a_{2i-1,j}(n) =a_{2i-1,j+2^{n-2}}(n)$, ($1\leq j\leq 2^{n-2}$), by easy induction. Thus we need only consider $a_{2i-1,2j-1}(n) $, and putting $m\equiv i$, $p\equiv j$ (mod $2^{n-3}$) for $1\leq m,p \leq 2^{n-3}$, we have
\[
\begin{array}{rlll}
a_{2i-1,2j-1}(n) &=&a_{2m-1,j}(n-1)&\mbox{by construction}, \vt
&=&a_{j,2m-1}(n-1)&\mbox{by induction}, \vt
&=&a_{p,m}(n-2)&\mbox{by construction}, \vt
&=&a_{m,p}(n-2)&\mbox{by induction}, \vt&=&a_{m,2p-1}(n-1)&\mbox{by construction}, \vt&=&a_{2p-1,m}(n-1)&\mbox{by induction}, \vt&=&a_{2j-1,2m-1}(n)&\mbox{by construction}, \vt&=&a_{2j-1,2i-1}(n).&\mbox{Hence the result.}
\end{array}
\]
We now put 
\[
N=2^{n-1},\;\;\;I_{j,n} = \left[\left.\frac{j-1}N;\,\frac jN\right)\right.\;\;(j\leq N),
\]
and take $S$ to be the set of $I_{j,n}$ for all $j,n$. We suppose $S\subset K$.

Let $h_i(x) = a_{ij}(n)$ for $x$ in $I_{jn}$ and $i \leq N$. By construction, $h_i(x)$ is independent of $n$ for $N \geq i$, and $||h_i|| =1$. By (7),
\[
\int_0^1 h_i(x)h_j(x) \,dx =0,\;\;\;\;(i \neq j).
\]
And by (6), every function $f$ constant in intervals $I_{jn}$ ($j=1, \ldots ,N$) for some $n$ depending on $f$, is linearly dependent on a finite number of $h_1,h_2, \ldots$. 
The functions $f$ are everywhere dense in the set of step-functions, which are everywhere dense in $\sigma_2(f)$. Hence we have proved

\vs
\noindent
\textbf{(9)}
\textit{$h_1(x), h_2(x), \ldots $ form a complete orthonormal set for $\sigma_2(f)$.}

\vs
We now have the principal results of this appendix.

\vs
\noindent
\textbf{(10)}
\textit{When $f(x)=\sum_{i=1}^\infty b_ih_i(x)$ and $f_n(x) =\sum_{i=1}^n b_ih_i(x)$ then
\[
P(f_N) = N\sum_{j=1}^N F_{jn}g(I_{jn})\]
where
\[ 
N=2^{n-1},\;\;I_{j,n} = \left[\left.\frac{j-1}N;\,\frac jN\right)\right.\;\;(j\leq N),\;\;F_{jn}=\int_{I_{jn}}f(x)\,dx.
\]}
\textbf{(11)}
\textit{When $f(x)=c(E;x)$ then $F_{jn} =m(EI_{jn})$, and since $mI_{jn}= 1/N$, we have a density integration with a single Riemann succession $\{D_n\}$ over $W$, where $D_n$ is composed of $I_{1,n},I_{2,n}, \ldots , I_{N,n}$.}

\vs
\noindent
\textbf{(12)}
\textit{When $P(f)$ is continuous,
\[
P(f) = \lim_{n \rightarrow \infty} N\sum_{j=1}^N F_{jn} g(I_{jn}).
\]}
For 
\begin{eqnarray*}
P(f_N)& =& \sum_{i=1}^N b_iP(h_i)\;\;\;=\;\;\;
\sum_{i=1}^N \sum_{j=1}^N a_{ij}(n)F_{jn}\sum_{k=1}^N a_{ik}(n)g(I_{kn})\vt
&=& 
\sum_{j=1}^N \sum_{k=1}^N F_{jn}g(I_{kn})\sum_{i=1}^N a_{ij}(n)a_{ik}(n) \;\;\mbox{ and by (8), and then (7)},\vt
&=&
\sum_{j=1}^N \sum_{k=1}^N F_{jn}g(I_{kn})\sum_{i=1}^N a_{ji}(n)a_{ki}(n)\;\;\;
=\;\;\; N\sum_{j=1}^N F_{jn}g(I_{jn}).
\end{eqnarray*}
This gives (10). For (11) it is obvious that $F_{jn}=m(EI_{jn})$. 

If $F_2$ is the family consisting of $\{D_n\}$ with other successions as in axiom (ii), then $F_2$ would not satisfy any other axiom save (vii), so that very little of the theory of Chapters 1 and 3 can be applied. However we can define the corresponding $G_2$ and
\[
(D;G_2)\overline{\int}_E\,g(I),\;\;\;\;\;\;
(D;G_2)\underline{\int}_E\,g(I).
\]
For (12), since $f_N \rightarrow f$ and $P(f)$ is continuous, we have $P(f_N) \rightarrow P(f)$, and this with (10) gives (12).

\vs
\noindent
\textbf{(13)}
\textit{When $P(f)$ is continuous then
\[
P\left(c(E;x)\right) = (D;G_2){\int}_E\,g(I)=g(E).
\]
where $g(E)$ is the Lebesgue extension of $g(I)$.}

\vs
\noindent
For by the argument in (3), taking $K=\bar S$, we see that $g(I)$ is absolutely continuous in $\bar S$, and then we use III(6.01).

\section{Appendix 2. \\
Further Generalisations.}\label{Appendix 2}
As already noted in Chapter 2, the axioms (i) to (viii), save (vi), can at once be generalised to $n$ dimensions, so that those results which depend only on these axioms can at once be stated for the $n$-dimensional theory. It is not possible to generalise (vi) in a straightforward manner to two dimensions, and (vi)' is rather involved.

Let us now take the case of an abstract space $W$ containing  sets $R$, in which there is an operation of \textit{addition}, denoted by $+$, and having the following properties.
\begin{enumerate}
\item[(a)]
Let $R_1,R_2$ be in $W$. If $R_1+R_2$ exists it is in $W$. If $R_1+R_2$ exists we say that $R_1$ and $R_2$ are \textit{disjoint}.
\item[(b)]
$R_1+R_2=R_2+R_1$, in the sense that if one side exists so does the other. Similarly for (c).
\item[(c)]
$R_1+(R_2+R_3)=(R_1+R_2)+R_3$. We write this as $R_1+R_2+R_3$ and suppose that it exists if $R_1+R_2$, $R_1+R_3$, and $R_2+R_3$ all exist.
\item[(d)]
If $R_1+R_2$ does not exist then one and only one of (d1), (d2), (d3), (d4) holds.
\begin{itemize}
\item[(d1)]
For $R_3$ in $W$, $R_1=R_2+R_3$, we then put $R_2 \subset R_1$ ($R_2$ is contained in $R_1$), and $R_3$ as $R_1-R_2$.
\item[(d2)]
For $R_3$ in $W$, $R_2=R_1+R_3$.
\item[(d3)]
$R_1=R_2$, when we also put $R_1 \subset R_2$, $R_2 \subset R_1$.
\item[(d4)]
There is an $R$ in $W$, denoted by $R_1.R_2$ such that 
\[
R_1.R_2 \neq R_1\;\mbox{ or }\;R_2,\;\;\;\;\;\;
R_1.R_2 \subset R_1,\;\;\;\;\;\;R_1.R_2 \subset R_2,
\]
and if $R'$ in $W$ is such that $R'\subset R_1$, $R'\subset R_2$, then $R'\subset R_1.R_2$. We then say that $R_1$ and $R_2$ \textit{overlap} with \textit{common part} $R_1.R_2$, and we denote by $R_1-R_2$ the set $R_3$ in $W$ such that $R_1.R_2 +R_3 =R_1$.
\end{itemize}
\item[(e)]
If (d1), (d2), or (d3) is satisfied, or if there is an $R$ of $W$ contained in both $R_1$ and $R_2$, then $R_1+R_2$ does not exist.
\item[(f)]
If $R_1+R_2 =R_1+R_3$ then $R_2=R_3$. We cannot have $R_1=R_1+R_2$.
\end{enumerate}
The properties (d) and (e) replace axiom (vii) in the theory of Chapter 1.

Given a subspace $S$ of $W$ such that each $R$ of $W$ is a ``finite sum'' of sets in $S$, we denote the general member of $S$ by $I$. Then $W$ takes the place of $\bar S$ in the theory of Chapter 1.

If rules exist which associate a definite set of numbers $\{g(I)\}$ with each set $I$ of $S$, then \textit{$g(I)$ is a many-valued function of the sets $I$ of $S$.}
(The set of numbers $\{g(I)\}$ takes the place of the different numbers for different bracket conventions, and also enables us to deal directly with the function $g(I;\xi)$ for $\xi$ in $\rho(I)$.)

Let $R$ in $W$ be the sum of $I_1, \ldots, I_n$ in $S$, i.e.~$R=\sum_{i=1}^n I_i$. Then $I_1, \ldots ,I_n$, together with a rule (or \textit{value convention}) which picks out a definite value from $\{g(I_i)\}$ for each $I_i$ ($i=1, \ldots ,n$) are said to form a \textit{division} $D$ of $R$. Similarly, if $R$ is in $S$ then $R$ with a value convention forms a \textit{division} $D$ of $R$.
Similarly, if $R$ is in $S$ then $R$ with a value convention forms a \textit{division} $D$ of $R$.

Summation of \textit{numbers} for a division $D$ over a set $E$ of objects is denoted by $(E;D)\sum$. If a division $D$ is not in question, or is assumed known, we put $(E)\sum$ for the convention. If also $E$ is assumed known we put $\sum$.

We suppose that to each $I$ of $S$ there is associated a definite strictly positive number $n(I)$, called the \textit{norm} of $I$, such that if $I_1 \subset I_2$ then $n(I_1)\leq n(I_2)$. Then if $D$ is a division of $R$ consisting of $I_1, \ldots ,I_m$ and value conventions, the \textit{norm} of $D$ is
\[
\mbox{norm}(D) = \max\{n(I_1), \ldots , n(I_m)\}.
\]
If $D_1,D_2, \ldots$ is a sequence of divisions of $R$ in which norm$(D_i)$ tends to 0 as $i \rightarrow \infty$, then $\{D_i\}$ is called a \textit{Riemann succession} of divisions of $R$.

We can now consider, as before, a family $F$ of Riemann successions of divisions of the $R$ in $W$, such that $F$ obeys one or more of the following axioms.
\begin{enumerate}
\item[(i)]
If $R$ is in $W$ there is at least one Riemann succession of divisions of $R$ which is n $F$.
\item[(ii)]
If $\{D_n\}$ is a Riemann succession in $F$ and $n_1, \ldots , n_i, \ldots $ is any sequence of integers tending to infinity then $\{D_{n_i}\}$ is in $F$.
\item[(iii)]
If $R_1,R_2, R_1+R_2$ are all in $W$, and if $\{D_{1,n}\}, \{D_{2,n}\}$ are in $F$, where $\{D_{i,n}\}$ is a Riemann succession of divisions of $R_i$ ($i=1,2$), then $\{D_n\}$ is in $F$, where $D_n$ is the division of $R_1+R_2$ formed from the $I$, and their value conventions, of the divisions $D_{1,n}$ and $D_{2,n}$, i.e., $D_n = D_{1,n}+D_{2,n}$.
\item[(iv)]
Let $\{D_n^{(i)}\}$ ($i=1,2, \ldots $) be a set of Riemann successions in $F$, of divisions of $R$ in $W$. Then if
\[
D_1^{(1)},\;\;\;D_2^{(1)},\;\;\;D_1^{(2)},\;\;\;D_3^{(1)},\;\;\;D_2^{(2)},\;\;\;D_1^{(3)},\;\;\;D_4^{(1)},\ldots
\]
is a Riemann succession it is in $F$.
\item[(v)]
Let $R=\sum_{i=1}^m R_i$ with the property that \textit{every} division $D$ of $R$ can be divided up to form  $D^{(1)}, \ldots ,D^{(m)}$ such that $D^{(i)}$ is a division of $R_i$ ($i=1, \ldots ,n$). Then if the Riemann succession $\{D_n\}$ of successions of $R$ is in $F$, and $D_n$ gives $D_n^{(i)}$ over $R_i$ ($i=1, \ldots ,m$), we have 
$\{D_n^{(1)}\}, \ldots ,\{D_n^{(m)}\}$ also in $F$.

Let $G$ be the family of divisions from all the Riemann successions in $F$.
\item[(viii)]
If $I$ is in $S$, and we take any value convention for $I$, then there is a division $D$ of $G$ and over $I$ such that $D$ is $I$ with the chosen value convention.
\end{enumerate}
We are now left with an axiom like (vi).

Let $\{D_n\}$ in $F$ be over $R=R_1+R_2$. Then $R_1,R_2$ \textit{have the property} (r) \textit{with respect to }$\{D_n\}$ if the following are satisfied.
\begin{enumerate}
\item[(A)]
$\{D_{1,n}\}$ over $R_1$ and $\{D_{2,n}\}$ over $R_2$ are both in $F$, where $\{D_{1,n}\}$  and $\{D_{2,n}\}$ are defined in (B), (C).
\item[(B)]
If $J$ of $S$ and in $D_n$ is contained in $R_i$ then $J$ occurs in $D_{i,n}$ with the same value convention there, ($i=1,2$).
\end{enumerate}
Now if $J \subset R_1$ there is $R_3$ in $W$ such that $R_1=R_3+J$, or else $J=R_1$. In the latter case, $J+R_2 =R_1+R_2$ exists. In the former case,
\[
R_1+R_2 = (R_3+J) +R_2 =R_3 +(J+R_2)
\]
so that again, $J+R_2$ exists. Thus $J$ cannot be in $R_2$, nor be equal to $R_2$, nor overlap with $R_2$. Similarly if $J \subset R_2$.
\begin{enumerate}
\item[(C)]
We now suppose it false that $J \subset R_1$.
\textit{If $R_1 \subset J$ and $R_1 \neq J$, then there is an $R_3$ such that $J=R_1+R_3$. Let $D_n$ be $J, I_1, \ldots ,I_m$, with value conventions. Then
\[
R_1+R_2 = J+I_1+ \cdots +I_m,
\]
so that by (f) we have $R_3 \subset R_2$.} The $I_1, \ldots ,I_m$ are in $R_2$ and are dealt with in (B).
\end{enumerate}
If there are no $I_1, \ldots,I_m$, then obviously $R_1+R_2 =J=R_1+R_3$, $R_2=R_3$, by (f). For the $R_3$ we suppose
\begin{enumerate}
\item[(C1)]
In $D_{1,n}$ there are some sets $I$ with sum $R_1$; and in $D_{2,n}$ there are some sets $I$ with sum $R_3$. \textit{Let $J_1, \ldots ,J_p$ be the sets $I$ of $D_n$ which each overlap with $R_1$.}
\end{enumerate}
If at least $J_1$ exists then there is no $I$ of $D_n$ such that $R_1 \subset I$. For then \[J_1.R_1 \subset R_1 \subset I\]
so that $I$ and $I_1$ have a common part in $W$. This is impossible since $J_1+I$ exists. Hence we do not have case (C1). And obviously (B) cannot hold for $J_1,\ldots ,J_p$, or else $R_1$ and $R_2$ would have a common part.

Now $R_{3,i}=J_i -R_i$ can be added to $R_1$. If $R_{3,i}+R_2$ does not exist, and $R_{3,i}\neq R_2$, then $R_{3,i}$
contains or overlaps with $R_2$, or $R_{3,i}\subset R_2$. Except in the last case,
$ R_{3,i}= R_{4,i}+R_{5,i}$ where $R_{4,i} \subset R_2$, and $R_{5,i}$ can be added to $R_2$. But
\[
R_{3,i}+R_{1}=\left(R_{4,i}+R_{5,i}\right) +R_{1}
=R_{4,i} + \left(R_{5,i}+R_{1}\right)
\]
so that $R_{5,i}$ can also be added to $R_{1}$. But
$R_{5,i} \subset J_i$ so that $\sum_{i=1}^p \sum'R_{5,i}$ exists over all the $R_{5,i}$, which exist (by (c)). If $R_{3,i} +R_2$ does exist we put $R_{5,i} = R_{3,i}$.

These are for $i=1, \ldots ,p$. Naturally, if $R_{3,i} \subset R_2$, $R_{5,i}$ does not exist.

\noindent
\textit{Out of the sets $I$ in $D_n$ let $J_1', \ldots ,J_q'$ be the sets which can be added to $R_1$.}
Then $J_j'=J_{1,j}'$ or $J_{2,j}'$ or $J_{1,j}'+J_{2,j}'$ (as  for $J_i$), where $J_{1,j}' \subset R_2$, and $J_{2,j}'+R_2$ exists ($j=1, \ldots ,q$). Then we have
\[
R_1+R_2 =R=J_1+\cdots +J_p + J_1'+ \cdots + J_q' + I_1 +\cdots +I_r
\]
where $I_k \subset R_1$ ($k=1, \ldots ,r$) and so satisfy (B), so this $=$
\[
=\left(
\sum_{i=1}^p J_i.R_1 +\sum_{k=1}^r I_k \right)
+\left(
\sum_{i=1}^p R_{4,i} 
+\mbox{$\sum'$}_{\;i=1}^{\;p} R_{5,i}
+\sum_{j=1}^q J_{1,j}'
+\sum_{j=1}^q J_{2,j}'\right) = R_6 + R_7
\]
where $R_6 \subset R_1$ and $R_7 +R_1$ exists. Let $R_1=R_6 +R_8$. Then by (f),
\[
R_8 + R_2 = R_7;\;\;\;\;\;\;\;\;\;
R_1+R_7 = R_1 +(R_8+R_2)=(R_1+R_8) +R_2;
\]
i.e.~$R_1+R_8$ exists. This contradicts (d). Hence $R_1=R_6$, so that by (f), $R_2=R_7$, i.e.,
\[
R_2 =
\left(
\sum_{i=1}^p R_{4,i} +\sum_{j=1}^q IJ'_{1,j}\right)
+\left(\mbox{$\sum'$}_{\;i=1}^{\;p} R_{5,i}
+\sum_{j=1}^q J_{2,j}'\right), = R_9 +R_{10}
\]
where $R_9 \subset R_2$ and $R_{10}+R_2$ exists. Hence, as before, $R_9=R_2$ and $R_2 = R_2 +R_{10}$. Hence by (f), $R_{10}$ cannot exist, so that
\[
J_i =J_i.R_1 +R_{4,i}\;\;\mbox{ where }\;\;J_i.R_1 \subset R_1,\;\;\;R_{4,i} \subset R_2 \;\;\;(1 \leq i \leq p).
\]
Also $J_j' = J'_{1,j}$ ($1 \leq j \leq q$) which therefore satisfy (B). Thus there only remains
\begin{enumerate}
\item[(C2)]
For $i=1, \ldots ,p$ there are some sets $I$ of $S$ in $D_{1,n}$ with sum $J_i.R_1$, and some sets $I$ in $D_{2,n}$ with sum $R_{4,i}$.
\end{enumerate}
We can now set up axiom (vi) and define $b(R_1;R_2;F)$, $c(R_1;R_2;G)$.
\begin{enumerate}
\item[(vi)]
If $R_1+R_2$ exists then $R_1,R_2$ have the property (r) with respect to every Riemann succession $\{D_n\}$ in $F$ and over $R_1+R_2$.
\end{enumerate}
We define $B\left(R_1;R_2;\{D_n\}\right)\;\;=$
\[
=\;\;
{\overline{\lim}}_{n \rightarrow \infty}
\left| (R_1+R_2;D_n) \sum g(I) -(R_1+R_2;D_{1,n}+D_{2,n}) \sum g(I) \right|
\]
and
$
b(R_1;R_2;F) = \mbox{l.u.b.}B\left(R_1;R_2;\{D_n\}\right)
$ for all $\{D_n\}$ in $F$ and over $R_1+R_2$.

Similar definitions, using norm-limits, may be given for
$C\left(R_1;R_2;G;e\right)$ and $c\left(R_1;R_2;G\right)$, and we can obtain analogues of I(3.03), I(3.04), and I(3.06b). The analogue of I(3.08) is as follows.

\vspace{10pt}
\noindent
\textit{Let $R=\sum_{i=1}^n R_i$. Then}

\noindent
\textbf{(a)}
\[
(F) \overline{\int}_R\,g(I) \leq\sum_{i=1}^n
(F) \overline{\int}_{R_i}\,g(I) 
+ \sum_{i=1}^{n-1} b\left(R_1+\cdots +R_i;R_{i+1};F\right),
\]

\noindent
\textbf{(b)}
\[
(F) \underline{\int}_R\,g(I) \geq\sum_{i=1}^n
(F) \underline{\int}_{R_i}\,g(I) 
- \sum_{i=1}^{n-1} b\left(R_1+\cdots +R_i;R_{i+1};F\right),
\]
\textit{Similar results hold for norm-limits.}

\vspace{10pt}
Analogues of I(3.09), \ldots ,I(3.11) now follow, and in this way an abstract theory can be built up. But since the thesis should only deal with functions of ($n$-dimensional) intervals, the abstract theory is really beyond the scope of the thesis.

We might have supposed that the values of $g(I)$ were objects in some space, for example, points in Banach space. But that degree of generality does not really seem necessary.

For a generalisation in another direction we can consider the relative differentiation with respect to an $h(I)$, of $g(I)$ and its integrals (e.g.~see Saks [10], 214 \textit{et seq.}) But this again is beyond the terms of reference of the thesis, which deals with the integration alone.

\newpage

\noindent
{\large{\textbf{Bibliography and References.}}}
\label{Bibliography}
\begin{itemize}
\item[{[1] 1910.}]
Z.~de Ge\"{o}cze, \textit{Quadrature des surfaces courbes}, Math.~Nat\-ur\-wiss.~Ber.~Ungarn. 26, 1--88.
\item[{[2] 1913.}]
W.H.~Young, \textit{On integrability with respect to a function of bounded variation,} Proc.~London Math.~Soc.~(2) 13 109--150.
\item[{1917.}]
W.H.~Young, \textit{On integrals and derivatives with respect to a function,} Proc.London Math.Soc.~(2) 15 35 \textit{et seq}.
\item[{[3] 1920.}]
E.W.~Hobson, \textit{On Hellinger's Integrals}, Proc.~London Math.~Soc.~(2) 18 249--265.
\item[{[4] 1922.}] Moore and Smith, \textit{A general theory of limits,} American J.~of Math.~44 102--121.
\item[{[5] 1923.}]
S.~Pollard, \textit{The Stieltjes integral and its generalisations,} Quart.~J.~of Math.~49, 73--138.
\item[{[6] 1924.}]J.C.~Burkill, \textit{Functions of intervals,} Proc.~London Math.~Soc.~(2) 22 275--310.
\item[{[7] 1924.}]J.C.~Burkill, \textit{The expression of area as an integral,} Proc.~London Math.~Soc.~(2) 22 311--336.
\item[{1924.}]J.C.~Burkill, \textit{The derivates of functions of intervals,} Fundamenta Math.~5 321--327.
\item[{[8] 1925.}]H.L.~Smith, \textit{On the existence of Stieltjes integral,} Trans.~American Math.~Soc.~27, 491--515.
\item[{[9] 1926.}]L.~Tonelli, \textit{Sur la quadrature des surfaces,} Comptes rendus Acad.~Sci.~Paris 182 1198--1200.
\item[{[10] 1927.}]S.~Saks, \textit{Fonctions d'intervalle,} Fundamenta Math.~10, 211 \textit{et seq.}
\item[{1928.}]R.C.~Young, \textit{Functions of  $\Sigma$ defined by addition, or functions of intervals in $n$-dimensional formulation,} Math.~Zeit\-schrift 29 171--216.
\item[{1929.}]R.C.~Young, \textit{On Riemann integration with respect to an additive function of sets,} Proc.~London Math.~Soc.~29 479--489.
\item[{1930.}]M.D.~Kennedy, \textit{Determinate functions of intervals and their rate of increase,} Proc.~London Math.~Soc.~(2) 30 58 \textit{et seq.}
\item[{1930.}]J.~Ridder, \textit{\"{U}ber approximative Ableitungen bei Punktund Intervallfunktionen,} Fundamenta Math.~15 324 \textit{et seq.}
\item[{1932.}]
S.~Kempisty, \textit{Sur les d\'eriv\'ees des fonctions des syst\`emes simples d'intervalles,} Bull.~Soc.~Math.~France 60 106--126.
\item[{1932.}]
F.~Riesz, \textit{Sur l'existence de la d\'eriv\'ee des fonctions d'une variable r\'eelle et des fonctions d'intervalle,} Verhandlungen des internationalen Mathematiker-Kongress Zurich I 258--269.
\item[{1934.}]
S.~Saks and A.~Zygmund, \textit{On functions of rectangles and their applications to the analytic functions,} Ann.~Scuola norm.~super.~Pisa 3 1--6.
\item[{[11] 1935.}] 
B.C.~Getchell, \textit{On the equivalence of two methods of defining Stieltjes integrals,} Bull.~American Math.~Soc.~41 413--418.
\item[{1935.}]
S.~Saks, \textit{On the strong derivatives of functions of intervals,} Fundamenta Math.~25, 235--262.
\item[{1936.}]
A.J.~Ward, \textit{On the differentiation of additive functions of rectangles,} Fundamenta Math.~26 167--182.
\item[{[12] 1936.}] 
S.~Kempisty, \textit{Sur les fonctions absolument continues d'inter\-valle,} Fundamenta Math.~27 10--37.
\item[{1936.}]
S.~Saks, \textit{On derivates of functions of rectangles,} Fundamenta Math.~27 72--76.
\item[{1936.}]
A.J.~Ward, \textit{On the derivation of additive functions of intervals in $m$-dimensional space,} Fundamenta Math.~28 265--279.
\item[{1937.}]
A.J.~Ward, \textit{A sufficient condition for a function of intervals to be monotone,} Fundamenta Math.~29 22--25.
\item[{[13] 1938.}]
T.H.~Hildebrandt, \textit{Definitions of Stieltjes integrals of Riemann type,} American Math.~Monthly 45 265--278.
\item[{1938.}]
A.J.~Ward, \textit{Remark on the symmetrical derivatives of additive functions of intervals,} Fundamenta Math.~30 100--103.
\item[{1941.}]
A.~Gleyzal, \textit{Interval-functions,} Duke Math.~Jour.~8 223--230.
\item[{1941.}]
P.~Reichelderfer and L.~Ringenberg, \textit{The extension of rectangle functions,} Duke Math.~Jour.~8 231--242.
\item[{[14] 1943.}]
P.~Dienes, \textit{On Riemann-Stieltjes integration,} lecture notes.
\item[{[15] 1946.}]
R.~Henstock, \textit{On interval functions and their integrals,} Jour.~London Math.~Soc.~21 204--209.
\item[{Jan.~1947}]
L.A.~Ringenberg, \textit{On the extension of interval functions,} Trans.~American Math.~Soc.~61 (Part 1) 134--146.
\item[{[16] 1928.}]
H.~Lebesgue, \textit{Le\c{c}ons sur l'Int\'egration et la Recherche des Fonctions Primitives,} Paris, 2nd edition.
\item[{[17] 1937.}]
S.~Saks, \textit{Theory of the Integral,} Warsaw, 2nd edition.
\item[{[18] 1939.}]
S.~Kempisty, \textit{Fonctions d'Intervalle non Additives,}
Actualit\'es Scientifiques et Industrielles, 824, Paris.
\item[{[19] 1939}]
E.C.~Titchmarsh, \textit{The Theory of Functions,} Oxford, 2nd edition.
\end{itemize}

\end{document}